\documentclass[11pt]{book}
\usepackage{latexsym,amssymb,amsmath,amscd,amsthm,amsxtra}

\title{\bf ANALYTIC TOPOLOGY \\
	of groups, actions, strings and varietes
	}
\author{Alexander Reznikov}

\date{August 11, 1999}

\makeindex

\begin{document}
\maketitle

\frontmatter

\setcounter{page}{1}
\pagenumbering{roman}
\tableofcontents

\renewcommand{\baselinestretch}{1.2}
\normalsize

\chapter*{Introduction}
This paper is devoted to an application of Analysis to Topology. The latter
is very broadly understood and includes geometric theory of finitely generated
groups, group cohomology, Kazhdan groups, actions of groups on manifolds,
superrigidity, fundamental groups of K\"ahler and quaternionic K\"ahler
manifolds and conformal field theory. The motivation and philosophy which has
led to the present research will be reflected upon in [Reznikov 7] and here
we will merely  say that we believe Analysis to be a major tool in studying
finitely generated groups. An alternative look is provided by arithmetic
method, notably by passing to a pro-p completion and using Galois cohomology.
This will be described in [Reznikov 8].

Each of six chapters which constitute this paper opens with a short overview;
a global picture is as follows. Chapter I and III treat analytic aspects of
geometry of finitely generated groups. Given an
immersion $M\hookrightarrow N$ of negatively curved manifolds ($M$ compact)
there is a boundary map $\partial\tilde M\to \partial\tilde N $, and it has
remarkable regularity properties. Invoking the Thurston theory, we show that
the actions $A$ of pseudoAnosovs on $W_p^{1/p}(S^1)/const $ have striking
properties from the viewpoint of functional analysis, namely,
$$ \sum_{n\in \mathbb{Z}}\|A^nv\|^p<\infty $$
for some $v\not= 0$. We apply this to a classical problem: when a surface
fibration is negatively curved and derive a strong necessary condition.

We then develop a theory of quantization for the mapping class group. A
classical work on $\mathcal{D}iff^{\infty}(S^1)$ suggests a two-step
quantization: first, obtaining a symplectic representation in
$Sp(W_2^{1/2}(S^1)/const)$ with image in the restricted symplectic group
[Pressley-Segal 1] and then using the Shale-Weil representation. The first  step
meets obstacles and the second step breaks down completely  for the mapping class group: first, because $\mathcal{M}ap_{g,1}$
does not act smoothly on $S^1$, so it's unclear why it can be represented in
$Sp(W_2^{1/2}(S^1)/const)$, second, if even it can (this happens to be the
case), there is no way to show that the image lies in the restricted
symplectic group (it almost certainly does not). The solution comes at the
price of abandoning the classical scheme and developing a theory of a new
object which we call bicohomology spaces $\mathcal{H}_{p,g}$. The mapping
class group $\mathcal{M}ap_g $ act in $\mathcal{H}_{p,g}$ and the latter shows
remarkable properties, like duality and existence of vacuum. The last
property is translated into the fact that $H^1(\mathcal{M}ap_{g,1},
W_p^{1/p}(S^1)/const ) $ is not zero. Finally, we find
$\mathcal{M}ap_g$-equivariant maps of the space of all discrete
representations of the surface group into $PSL_2(\mathbb{C})$ to our spaces
$\mathcal{H}_{p,g}$.

Chapter II uses Analysis to study groups, acting on the circle (we need
$\mathcal{D}iff^{1,\alpha}$ regularity, so $\mathcal{M}ap_{g,1}$ is not
included). Our first main theorem says roughly that Kazhdan groups do not
act on the circle. Very special cases of this result, for lattices in Lie
groups, were recently found (see the references in Chapter II). The Hilbert
transform, which played a major role in Chapter I, is crucial for the proof
of this result as well. We then develop a theory of higher characteristic
classes for subgroups of $\mathcal{D}iff^{1,\alpha}$ (the first being
classically known as an integrated Godbillon-Vey class). All this classes
vanish on $\mathcal{D}iff^{\infty}(S^1)$. It is safe to say that the less is
smoothness, the more interesting is the geometry "of the circle".

Chapter III brings us back to asymptotic geometry of finitely generated
groups. We propose, for a non-Kazhdan group, to study the asymptotic behaviour
of unitary cocycles. We prove a general convexity result which shows
that an embedding of $G$ in the Hilbert space, given by a unitary cocycle,
is "uniform". We then prove a growth estimate for unitary cocycles of a
surface group, using very heavy machinery from complex analysis, adjusted
for our situation. Similar result for cocycles in $H^1(G,l^p(G))$ has already
been given in Chapter I.

Chapter IV studies symplectomorphism groups. There is a misterious
similarity between groups acting on the circle and groups acting
symplectically on a compact symplectic manifold. In parallel with the 
above mentioned
result in Chapter II we show, roughly, that transformations of a Kazhdan
group acting on a symplectic manifold must satisfy a partial differential
equation. An example is $Sp(2n,\mathbb{Z})$ acting linearly on $T^{2n}$
and, very probably, $\mathcal{M}ap_g$ acting on the space of stable bundles
over a Riemann surface. (I don't know for sure if $\mathcal{M}ap_g$ is
Kazhdan). In dimension 2 the result is very easy and was known before. We
also introduce new charateristic classes for symplectomorphism groups, in
addition to the two series of classes defined in our previous papers, and
use them to express a volume of a negatively curved manifold through the
Busemann function on the universal cover.

Chapter V studies volume-preserving actions. We introduce a new technique
into the subject, that of (infinite-dimensional, non-positively curved)
spaces of metrics. We define a invariant of an action which is an infinum of
a displacement in the space of metrics and show that for an action of a
Kazhdan group which does not fix a $\log L^2$-metric, 
this invariant is positive (a
weak version of this result for the special case of lattices was known
before). We then turn to a major open problem, that of non-linear
superrigidity and prove what seems to be first serious breakthrough after
many years of effort.

Chapter VI deals with fundamental groups of K\"ahler and quaternionic
K\"ahler manifolds. The situation is exactly the opposite to the studied in
Chapters II and IV, namely, these groups tend to be Kazhdan. We first extend
our rationality theorem for secondary classes of flat bundles over
projective varietes to the case of quasiprojective varietes, answering a
question posed to us by P.Deligne. We then prove that a fundamental group
of a compact quaternionic K\"ahler manifold is Kazhdan, therefore providing
a very strong restriction on its topology. We also discuss polynomial growth of
the group cohomology classes for K\"ahler groups, proved nontrivial in a
previous paper.

The paper uses many different analytic techniques.
Within each chapter, there is a certain coherence in the point of view
adopted for study.

I started this project on a chilly evening of November, 1998 in an African
caf\'e in Leipzig and finished it on a hot afternoon of July, 1999 in
Jerusalem. The manuscript has been written up by August, 11, 1999; I would
appreciate any mentioning of a possible overlap with any paper/preprint
which appeared before this date. During the long time when the paper was being typed and then polished, I found a proof of several statements which had been conjectured in the paper, in particular  a construction of a cocycle for the group of quasisymmetric homeomorphisms valued in $W_p^{1/p}(S^1)/const$, which was conjectured in Chapter I. The proofs will appear in a sequel to this paper.

\mainmatter

\setcounter{page}{1}
\setcounter{chapter}{0}
\pagenumbering{arabic}

\chapter{Analytic topology of  negatively curved manifolds, quantum strings 
and mapping class groups}

Chapter I opens with simple observations concerning the cohomology
$H^1(G,l^{\infty}(G))$ for a finitely-generated group. If $G$ is amenable
we produce plenty of polynomial cohomology classes in $H^*(G,\mathbb{R})$
given by an explicit formula (Theorem 1.2.1). Then we prove a convexity theorem
1.2.2 saying that if there are Euclidean-type quasigeodesics in the Cayley
graph of $G$, then $G/[G,G]$ is infinite. 

We then review some standard facts on $l^p$-cohomology in sections 2--4. One defines  an 
asymptotic  invariant of a finitely generated group $G$,  called a
constant of coarse structure $\alpha(G)$, as an infinum of $p$, $1\le
p\le \infty $ such that $H^1(G,l^p(G))\not=0$. For all noncusp
discrete groups of motions of complete manifolds of pinched negative
curvature, $\alpha(G)<\infty$. For discrete subgroups $G$ of $SO^{+}(1,n)$,
$\alpha(G)\le \delta(G)$, where $\delta$ is the
exponent of the group.
In section 4 we review function spaces. A classical result in weighted Sobolev spaces may be reformulated as an  identification of the $l^p$-cohomology of cocompact real hyperbolic lattices:
$H^1(G,l^p(G))=W_p^{(n-1)/p}(S^{n-1})/const $ .
It follows  that $\alpha(G)=n-1$. 

In section 5 we prove a first result within a program to classify groups according to the cocycle growth. We show for surface groups, that if $L_g\mathcal{F}-\mathcal{F}\in l^p(G) $ for all
$g\in G$, then $|\mathcal{F}(g)|\le const\cdot [length(g)]^{1/p'} $. Here
$\mathcal{F}:G\to \mathbb{R} $ is any function (Theorem 5.6). This result with no doubt generalizes to higher-dimensional cocompact lattices in simple Lie groups of rank one.

In section 6 we present a new theory for boundary maps of negatively curved
spaces, associated with immersions of closed manifolds. The most striking is a
partial regularity result (Theorem 6.1, part 4).

As is well known, the group of quasisymmetric ($n=2$) or
quasiconformal ($n\ge 3$) homeomorphisms of $S^{n-1}$ act on
$W_p^{\frac{n-1}{p}}(S^{n-1}) $ for $p>n-1$. The action of $\mathcal{G}_1$
on $W_2^{1/2}(S^1)/const $ is in fact symplectic. We give application to the
regularity of quasisymmetric homeomorphisms (Theorem 8.2).
In Corollary 9.2 we prove that the unitary representation of a subgroup $G$ of $
SO(1,2)$ in $ W_2^{1/2}(S^1)/const $ is an invariant of a component of the
Teichm\"uller space $\mathbf{T}(G)$. 

In Theorem 9.3 we show striking properties of
invertible operators $A$ in Banach spaces $W_p^{1/p}(S^1)/const,\ p>2$,
induced by quasiAnosov maps in $\mathcal{M}ap_{g,1}$, namely
$$ \sum_{k\in\mathbb{Z}} \| A^k v \|^p < \infty $$
for some $0\not= v $. In Theorem  9.5 we find a new inequality in topological
Arakelov theory, based on the work of [Matsumoto-Morita 1]. In Theorem 9.6 we find very strong restrictions on a
subgroup $G\subset \mathcal{M}ap_g $, such that an induced group extension
$\widetilde G $:
$$ 1\to \pi_1(\Sigma_g)\to \widetilde G \to G \to 1 $$
is a fundamental group of a negatively curved compact manifold (this is a
classical problem). In section 10
we extend the theory to the limit case $p=1$, introducing an $L^1$-analogue
of Zigmund spaces, which we call $\mathcal{L}_{k,\alpha}$.

In section 11 we start a new theory of secondary quantization of Teichm\"uller
spaces. We introduce the bicohomology spaces $\mathcal{H}_{g,p}$ and show
that $\mathcal{M}ap_g$ acts on these spaces. We show (difficult!) that
$\mathcal{H}_{2,p} $ is an infinite-dimensional Hilbert space and there
is a symmetric bilinear nondegenerate form of signature $(\infty,m)$ which
is $\mathcal{M}ap_g$-invariant. What is the value of $m$, we don't know at
the time of writing of this introduction (August,1999). So does the secondary
quantization lead to ghosts? We provide a holomorphic realization in the
space of $L^2$-holomorphic 2-forms on $\mathcal{H}^2\times\mathcal{H}^2/G $
and $\mathcal{H}^2\times \overline{\mathcal{H}^2}/G $ (Theorem  11.12). In section 12 we interpret $\mathcal{H}_{p,g}$ as
operator spaces (proposition 12.2), and prove the existence of vacuum (Theorem
12.5). We prove that
$H^1(\mathcal{M}ap_{g,1}, W_p^{1/p}(S^1)/const )\not= 0 $ for $p\ge 2$. It
still may be true that $\mathcal{M}ap_{g,1}$ is Kazhdan, because the action is
not orthogonal.

In section 13 we construct $\mathcal{M}ap_g$-equivariant maps of the space
of discrete representations of the surface group in
$SO^{+}(1,3)=PSL_2(\mathbb{C})$ to our spaces $\mathcal{H}_{p,g}$ (Theorem 13.1). In
Theorem 13.2 we summarize our knowledge of the functional-analytic structure coming from
hyperbolic 3-manifolds which fiber over the circle.

\section{Metric cohomology}

\noindent\\{\bf 1.1.1.}
Let $G$ be a finitely generated group. Let
$\mathbb{K}=\mathbb{R},\mathbb{C}$. Let $V$ be a locally convex topological
$\mathbb{K}$-vector space which is a $G$-module, that is, there is a
homomorphism $G\to Aut(V)$. If $\{ g_i \},\quad i=1,\cdots,n $ is a finite
set of generators of $G$, then the evaluation map $f\mapsto \{ f(g_i) \} $
establishes an injective homomorphism $Z^1 (G,V)\to \Pi_{i=1}^n V $ of
the space of 1-cocycles of $G$ in $V$. One calls the induced topology in
$Z^1(G,V) $ the cocycle topology; it does not depend on the choice of
generators. A coboundary map $V\to Z^1 (G,V) $ may have an image 
$\overline {B^1(G,V)}$
which is not closed in $Z^1(G,V)$; the quotient $Z^1(G,V)/\overline{B^1(G,V)} $ is
called {\it reduced } first cohomology space. One way to produce nontrivial
cohomology classes is to consider limits of coboundaries, that is, elements
of $\overline{B^1(G,V)} / B^1(G,V) $. That amounts to considering
nets $\{ v_\alpha\in V \} $ such that $g_i v_\alpha -v_\alpha \to l(g_i) $ for
$i=1,\cdots,n$. If $V$ is a Banach space and $G$ acts isometrically without
invariant vectors, then $B^1(G,V)$ is closed in $Z^1(G,V) $ if and only if
there are no almost invariant vectors, that is, sequences $v_j, \|v_j\|=1 $,
such that $\| g_iv_j -v_j \|\to 0 $ for all $ i=1,\cdots,n$. This statement
is an immediate consequence of the Banach theorem and is called Guichardet's
lemma [Guichardet 1]. So if there are almost invaiant vectors, then
$H^1(G,V)\not= 0 $, though the reduced cohomology $H_{red}^1(G,V) $ may be
zero.

If $V$ is Banach and $G$ acts isometrically, let $l\in Z^1(G,V) $ be a
cocycle. Then $$\| l(g) \|\le \max_{i=1}^n \| l(g_i) \|\cdot length(g) $$,
where $length(g)$ is the length of the element $g$ in the word metric,
induced by $\{ g_i \} $. The proof is immediate by induction, using the
cocycle equation $l(gh)=gl(h)+l(g)$.

Now let $V_j,\quad j=1,\cdots,m$ be a collection of Banach spaces on which
$G$ acts isometrically and let $\varphi: \otimes_{j=1}^m V\to \mathbb K $ be a map
continuous in a sense that $\varphi (\otimes_{j=1}^m v_j ) \le const\cdot
\Pi_{j=1}^m \| v_j \| $. Let $ l_j \in Z^1(G,V_j) $ and let 
$l\in Z^m(G,\mathbb K) $
be the cup product $l(g_1,\cdots,g_m)=\varphi ( \otimes_{j=1}^m l_j(g_j) ) $.

\noindent\\
{\it Lemma 1.1.}--- $l\in Z^m(G,\mathbb{K}) $ is of polynomial growth,
more precisely
$$ |l(g_1,\cdots,g_m)|\le const\cdot \Pi_{i=1}^m length (g_i) .$$

\noindent\\
{\it Proof.}--- is immediate from the remarks made above.

A general definition of polynomial cohomology is to be found in
[Connes-Moscovici 1]. As we will see, Lemma 1.1 is a very powerful tool for
constructing cocycles of polynomial growth in concrete situations.

\noindent\\
{\it Proposition 1.1.}--- Let $G$ be an infinite finitely generated group.
Consider a left action of $G$ on $l^{\infty}(G)$. Then
$H^1(G,l^{\infty}(G))\not= 0 $. Moreover, $H^1(G,l_0^{\infty}(G))\not= 0$.

\noindent\\
{\it Proof.}--- Let $\{g_i \} $ be a finite set of generators of $G$, and let
$length(g)$ be a word length of an element $g$. Define a right-invariant
word metric by $\rho (x,y) =length (xy^{-1})$. Let $x_0\in G$ and let
$F(x)=\rho(x_0,x)$. Obviously, $F$ is unbounded. Now let $l(g)=L_g F-F $
where $L_g$ is a left action on functions, that is,
$l(g)(x)=F(g^{-1}x)-F(x)$. We find $$|l(g)(x)|=|\rho
(x_0,g^{-1}x)-\rho(x_0,x) |\le |\rho (g^{-1} x,x)|=\rho(g^{-1},1).$$  So $l$
is a cocycle of $G$ in $l^{\infty}(G)$. If it were trivial, we would have a
bounded function $f$ such that $L_g F-F=L_g f-f $ that is, $F-f$ would be
invariant, therefore constant, a contradiction. The second statement of the
Proposition  will be proved later in section 1.3.

\noindent\\
{\bf 1.1.2.}
Now let $G$ be amenable. In this  case we have a continuous map $\varphi:
\Pi_{j=1}^m l^{\infty}(G) \to \mathbb{K} $ given by $(f_1,\cdots,f_m)\mapsto
\int_G f_1\cdots f_m $. By an integral we mean a left-invariant normalized
mean of bounded functions. We obtain

\noindent\\
{\it Theorem 1.2.1.}--- Let $G$ be a finitely generated amenable group, let
$\rho_j,\ j=1,\cdots,m $ be a collection of right-invariant word metrics
on $G$. A formula
$$ l(g_1,\cdots,g_m)=\int_G \Pi_{j=1}^m [\rho_j (x_0,
g_j^{-1}x)-\rho_j(x_0,x)] $$
defines a real-valued $m$-cocycle on $G$ of polynomial growth: $$|
l(g_1,\cdots,g_m)|\le const\cdot \Pi_{j=1}^m length(g_j) $$ for any word
length $length(\cdot)$.

\noindent\\
{\it Examples.}--- Let $G=\mathbb{Z}$. If we choose generators $\{ -1,1\}$, then
$ length (g)=|g| ,$ and $$\rho (x_0,g^{-1} x)-\rho(x_0,x) =
|x_0-x+g|-|x_0-x|\to \pm |g| $$ as $ x\to \pm\infty $ and
$$\int_{\mathbb{Z}}(|x_0-x+g|-|x_0-x|)=0 .$$ However, if we choose generators
$\{ -1,2\} $, then
\renewcommand{\baselinestretch}{1.5}
\large
$$length(g)=\left\{ \begin{array}{cc}
    |g|,& g\le 0 \\
    \frac{g}{2}, & g\ge 0 \mbox{ and even }\\
    \frac{g+1}{2}, & g\ge 0 \mbox{ and odd }
    \end{array}
    \right.
$$
\renewcommand{\baselinestretch}{1.2}
\normalsize
Then $$length(x_0+g-x)-length(x_0-x) $$ for $g>0 $ and even will have limits
$\frac{g}{2}$ when $x\to -\infty $ and $-g$ when $x\to \infty $, so
$$\int_{\mathbb{Z}}[length(x_0+g-x)-length(x_0-x)]=-\frac{g}{4} .$$ So we
obtain a cocycle $l: \mathbb{Z}\to \mathbb{R} $ given by $g\mapsto
-\frac{g}{4} $. Now, if $G=\mathbb{Z}^k $, $k\ge 2$, let $\rho_j,
j=1,\cdots,k$ be a word metric defined by a set of generators
 $$\{ e_1^{\pm
1}, e_2^{\pm 1},\cdots,e_j^{-1},e_j^2,e_{j+1}^{\pm 1},\cdots,e_k^{\pm 1} \}$$
where $e_s$ is a generator of the s-th factor. If $1\le
j_1<j_2<\cdots<j_m\le k $ is a set of indices, then Theorem 2 provides a
cocycle
$$ l(g_1,\cdots,g_m)=\int_{\mathbb{Z}^k}\Pi_{r=1}^m [\rho_{j_r}(x_0,
g_j^{-1}x)-\rho_{j_r}(x_0,x)] .$$
If $\pi_i :\mathbb{Z}^k \to \mathbb{Z} $ is a projection to i-th factor,
then $l(g_1,\cdots,g_m)=(-\frac{1}{4})^m\cdot \Pi_{r=1}^m \pi_{j_r} (g_r) $.
It follows that classes of cocycles, given by Theorem 1.2.1, generate the real
cohomology space of $\mathbb{Z}^k$.

\noindent\\
{\it Remark.}--- If $G$ is amenable, $\rho$ is a right-invariant word metrics
and for some $x_0, g\in G$, 
$$\int_G [\rho(x_0,g^{-1}x)-\rho(x_0,x)]\not= 0,$$
then $H_1(G,\mathbb{R})\not= 0 $ and in fact $g\notin [G,G] $ for all
$s\not= 0$. This is a direct corollary of Theorem 1.2.1. A more interesting
structure theorem is given below.

\noindent\\
{\it Theorem 1.2.2.}--- Let $G$ be a finitely generated amenable group, $\rho $
a right-invariant word metric. Let $g\in G$, assume a following convexity
condition: there is some $C>0$, such that for any $x\in G$ there exists
$N\ge 0$ such that $\rho (g^k,g^{-1}x)-\rho(g^k,x)\ge C $ for $k\ge N$. Then
$H_1(G,\mathbb{R})\not= 0 $ and moreover, $g^s\notin [G,G] $ for all $s\not=
0 $.

\noindent\\
{\it Corollary 1.2.3.}--- Let $G$ be a Heisenberg group $\{ x,y,z |
[x,y]=z,[x,z]=[y,z]=1 \} $. Then for any right-invariant word metric $\rho
$, there exists $a\in G$ such that $\liminf_{k\to \infty}
[\rho(z^k,z^{-1}a)-\rho(z^k,a)]\le 0 $.

\noindent\\
{\it Proof of the Corollary.}--- Since $z\in [G,G] $, the result follows from
Theorem 1.2.2. Indeed $G$ is nilpotent, therefore amenable.

\noindent\\
{\it Proof of the Theorem.}--- Consider a 1-cocycle
$l(\gamma)(x)=\rho(x_0,\gamma^{-1}x)-\rho(x_0,x) $, $l\in
Z^1(G,l^{\infty}(G)).$ Set $x_0=g^n$, so

$$l_n(g)(x)=\rho(g^n,g^{-1}x)-\rho(g^n,x) .$$
If for any $x$ and sufficiently
big $n$, $\rho(g^n,g^{-1}x)-\rho(g^n,x)>C $ then a pointwise limit
$\lim_{n\to \infty} l_n(g)(x) $ exists and is $\ge C$. Since $| l_n(z)(x)
|\le \rho(z^{-1},1)$, there is a subsequence $n_k$ such that $l_{n_k}(z) $
converges pointwise for any $z$ to a bounded function $l(z)$. One sees
immediately that $l: G\to l^{\infty}(G) $ is a cocycle, so $z\mapsto \int_G
l(z) $ is a homomorphism from $G$ to $\mathbb{R}$. Since $l(g)(x)\ge C>0 $
for all $x$, $\int_G l(g)\ge C>0$, so $ H_1(G,\mathbb{R})\not= 0 $ and $
g^s\notin [G,G] $, as desired.

\noindent\\
{\bf 1.1.3.}
Let $\varphi: R_{+}\to R_+ $ be a smooth function such that $\varphi(x)\to
\infty $ as $x\to \infty $ and $\varphi'(x)\to 0 $. Let $G$ be a finitely
generated group and let $\rho$ be a right-invariant word metric. Consider
$F(x)=\varphi (\rho(x_0,x) ) $ where $x_0\in G$ is a fixed element. Since
\renewcommand{\baselinestretch}{1.5}
\large
$$\begin{array}{rl}
	|(L_g F-F)(x) |&=|F(g^{-1}x)-F(x) | \\
    	&=|\varphi(\rho(x_0,g^{-1}x))-\varphi(\rho(x_0,x))|   \\
        &\le \sup_{t\in I}|\varphi'(t)|\cdot |\rho(g^{-1}x,x)|  \\
        &\le \sup_{t\in I}|\varphi'(t)|\rho(g^{-1},1)
    \end{array}
$$
\renewcommand{\baselinestretch}{1.2}
\normalsize
where $I=[\min (\rho(x_0,x),\rho(x_0,g^{-1}x)),
\max(\rho(x_0,x),\rho(x_0,g^{-1}x))] $, we see that $L_g F-F \in
l_0^{\infty} $. Therefore $H^1(G, l_0^{\infty}(G))\not= 0 $, because the
cocycle $ L_g F-F $ cannot be trivial as a cocycle valued in $l_0^{\infty} $
(by the same reasons as in the proof of the first statement of Proposition
1.1 ) . The proof of Proposition 1.1 is now complete.

Notice that, since $\rho(u,v)=length (u\cdot v^{-1})$,
$$\rho(x_0,x)-length(g)\le \rho(x_0,g^{-1}x)\le \rho(x_0,x)+length(g),$$ so
that $$ |(L_g F-F)(x)|\le \sup_{|t-\rho(x_0,x)|\le length(g)}
|\varphi'(t)|\times \rho(g^{-1},1).$$

\noindent\\
{\it Remark.}--- Let $S(N)=\{ g| length(g)=N \} $. If $S(N)/S(N-1)\to 1 $ and
$\sum_{k=1}^N S(k)/S(N)\to \infty $ as $N\to \infty $, then for $p>1$ there
is a radial function $F(x)=\varphi(\rho(x))$ such that $L_gF-F\in l^p(G) $
and the cocycle $l: G\to l^p(G) $ defined by $g\mapsto L_g F-F $ is not a
coboundary. Note that $G$ is automatically amenable. On the other hand, if
$S(N)\sim e^{cN} $, then such radial function does not exist. This follows
at once from Hardy's inequality. To produce classes in $H^1(G,l^p(G))$, one
needs to use some more elaborate geometry than just distance function. In the
next section we produce such classes for negatively curved groups/manifolds,
using the visibility angles.

\section{Constants of coarse structure for negatively curved groups}

\noindent\\{\bf 1.2.1.}
Throughout this section we assume that $G$ is a finitely generated,
non-amenable group, therefore $B^1(G,l^p(G))$ is closed in $Z^1(G,l^p(G))$
for $p\ge 1$.

\noindent\\
{\it Definition 2.1.}--- A number $\alpha (G)=\inf_{1\le p\le \infty} \{ p|
H^1(G,l^p(G))\not= 0 \} $ is called a constant of coarse structure of $G$.

\noindent\\
{\it Remark.}--- The definition makes sense since by Proposition 1.1,
$H^1(G,l^{\infty}(G))\not= 0$. We will need a proof of the following well-known fact (see, for example [Pansu 1]). The argument below is a slightly modified, from nonpositive curvature to negative curvature, version of a classical argument of [Mishchenko 1,2].

\noindent\\
{\it Proposition 2.1.}--- Let $M^n$ be a complete Riemannian manifold of negative
curvature, not a cusp, satisfying $K(M)\le -1, Ric(M)\ge -(n-1)K $. Let
$G=\pi_1 (M)$. Then $\alpha (G)\le (n-1)\sqrt{K}$.

\noindent\\
{\it Proof.}--- Let $q_0\in \widetilde M $. Consider a map of $G$ onto an orbit
$\mathcal{O}$ of $q_0: g\mapsto gq_0 $; it is equivariant with respect to
the left action of $G$ on itself. Let $q\notin \mathcal{O} $ and let $v_q(s)
$ be an outward pointing vector from $q$ to $s$, that is, a unit vector in
$T_s\widetilde M$, tangent to geodesic segment joining $q$ and $s$. Consider for
$x\in G$, $F(x)=v_q(xq_0)$. Notice that $F(x)$ takes values in
$T_{xq_0}\widetilde M $. We can consider the restriction of $T\widetilde M$ on
$\mathcal{O} $ as an equivariant vector bundle over $\mathcal{O}$. 
Pulling back to $G$, we obtain an left-equivariant vector bundle over
$G$, equipped
with an equivariant Euclidean structure. Then $F$ is a section of this
bundle. Now consider $(L_gF-F)(x)$. Since the action of $G$ on sections is
given by $L_g F(x)=g_* F(g^{-1}x) $, where $g_*$ is the derivative map (
$ g_* :T_s \widetilde M \to T_{gs} \widetilde M$), we get $ (L_g F-F)(x)=g_*
F(g^{-1}x)-F(x)=g_* v_q (g^{-1}xq_0)-v_q(xq_0)=v_{gq}(xq_0)-v_q(xq_0)$. So
$\| (L_g F-F)(x)\|=|2\sin \frac{1}{2}\sphericalangle
(gq,xq_0,q)|\le \sphericalangle (gq,xq_0,q) $.

Let $E|G$ be the equivariant Euclidean vector bundle considered above (the
pullback of $G$ of $ T\widetilde M |\mathcal{O} $). Let $L^p(E)$ be a Banach
space of $L^p$-sections of $E$. We claim $L_g F-F \in L^p(E) $ for
$p>(n-1)\sqrt{K}$. Let $r(x)=dist_{\tilde M} (q_0,xq_0 ) $. For $ g,q_0,q$ fixed we
have $\sphericalangle (gq,xq_0,q)\le const_1 \cdot e^{-r(x)} $ by a standard comparison
theorem, since $K(M)\le -1 $. On the other
hand, for fixed $\delta >0 $, $\#(x|r-\delta \le r(x)\le r+\delta )\le
const_2 e^{(n-1)\sqrt{K}r} $ by the Bishop's theorem.
Therefore $L_g F-F \in L^p(E) $ for $p>(n-1)\sqrt{K} $. Note we only need
that $G$ acts discretely in $\widetilde M$.

A map $l: G\to L^p(E) $ defined by $l(g)=L_gF-F $ is obviously a cocycle. If
it were trivial, we would have an $L^p$-section $s\in L^p(E)$, such that
$F-s$ is invariant. That means $g_*
(v_q(g^{-1}xq_0)-s(g^{-1}x))=v_q(xq_0)-s(x)$, or $v_{gq}(xq_0)-g_*
s(g^{-1}x)=v_q(xq_0)-s(x)$. Notice that since $\| F(x)\|=1 $, $F-s$ is
invariant and $\| s(g)\|\to 0 $ as $length(g)\to \infty $, $\| (F-s)(x)\|=1
$ for all $x$. In particular, $w=v_q(xq_0)-s(x) $ has norm one. Fix $x$ and
let $g$ vary. We get $\| v_{gq}(xq_0)-w\|=\| g_* s(g^{-1}x) \|\to 0 $ as $
length(g)\to \infty $. Let $P_+, P_- $ be an attractive and repelling fixed
points of $g$ on the sphere at infinity of $\widetilde M$.

Let $w_+, w_- $ be unit vectors in $T_{xq_0}(\widetilde M)$, tangent to
geodesics, joining $xq_0$ with $P_+,P_-$. Then $\|
v_{g^nq}(xq_0)-w_{\pm}\|\to 0 $ if $n\to \pm\infty$. It follows that
$w_{\pm}=w$. Therefore all elements of $G$ are parabolic and have a common
fixed point at infinity. So $M$ is a cusp, a contradiction. So
$H^1(G,l^p(E))\not= 0$. However, $l^p(E) $ is equivariantly isometric to
$l^P(G)\otimes T_{p_0}(\widetilde M) $. So $H^1(G,l^p(E))\simeq
H^1(G,l^p(G))\otimes T_{p_0}(\widetilde M)$. We deduce that $H^1(G,l^p(G))\not=
0$.

\noindent\\
The estimate of the Proposition  is sharp. We will see later that if $G$ is a
cocompact lattice in $SO^{+}(1,n)$, i.e. $K(M)=-1$, then $\alpha(G) $ is exactly
$(n-1)$.
Let now $G$ be a discrete nonamenable subgroup of $SO^{+}(1,n)$, or,
equivalently, $K(M)=-1$. Recall that the exponent $\delta(G)$ is defined by
$\delta(G)=\inf \{ \lambda | \sum_{g\in G} e^{-\lambda r(g)}<\infty
\}$ where $r(g)=dist_{\widetilde M} (p_0, gp_0 ) $ for some fixed $p_0\in
\widetilde M$. If $G$ is geometrically finite, then by a well-known theorem
[Nicholls 1] $\delta(G) $ is equal to the Hausdorff dimension of the limit
set $\dim(\Lambda (G) )\subset S^{n-1}$. Note that if $\Lambda (G)\not=
S^{n-1} $, then $\dim \Lambda (G)<n-1 $ by [Sullivan 1] and [Tukia 1].We now have

\noindent\\
{\it Proposition 2.2.}--- Let $G$ be a discrete subgroup of $SO^{+}(1,n)$, not a cusp
group. Then $\alpha(G)\le \delta (G) $.

\noindent\\
{\it Proof.}--- The Proposition  follows from the proof of the
Proposition 2.1. Indeed, we only need that $\sum_{g\in G} e^{-p r(g) }
<\infty $ to conclude that one has a cocycle $l: G\to l^p(G)$. It has been proven
already that this cocycle is not a coboundary.

\noindent\\
{\it Remark.}--- The relation of the constant of coarse structure to ``conformal dimension at infinity'' is discussed in [Pansu 2].

\noindent\\
\section{Function spaces: an overview}

For $s\ge 0$ an integer and fractional part of $s$ are denoted $[s]$ and$\{
s\} $ respectively. A Sobolev-Slobode\u{c}ky space
$W_p^s(\mathbb{R}^n)$,$(p>1)$ consists of measurable locally integrable
functions $f$ on $\mathbb{R}^n $ such that $D^{\alpha} f\in L^p(R^n) $ for
$|\alpha |\le [s] $ and
$$ \sum_{|\alpha |=[s] } \int\int \frac{|D^{\alpha}f(x)-D^{\alpha}f(y)
|^p}{|x-y|^{n+\{ s\} p } }\ dxdy <\infty $$
A space of Bessel potentials $H_p^s $ consists of functions $f$ for which
a Liouville-type operator
$$ \mathcal{D}^sf=((1+|\xi|^2)^{s/2} \hat f(\xi) )^{\wedge} $$
satisfies $\mathcal{D}^s f\in L^p $. Warning: $H_p^s\not= W_p^s $ if $s$ is
not an integer. For $p=2$ the condition is equivalent to $$(1+\triangle
)^{s/2} \hat f \in L^2(\mathbb{R}^n) .$$ Here $f(x)\to \hat f(\xi) $ is the
Fourier transform and $\triangle=-\sum \frac{\partial^2}{\partial x_i^2 } $.

A space of BMO functions BMO$(\mathbb{R}^n)$ is defined as a space of
functions $f$ for which $$\sup_{Q}\frac{1}{|Q|} \int_Q |f(x)-f_Q |\ dx<\infty
,$$ where $Q$ runs over all cubes in $\mathbb{R}^n$ and
$$f_Q=\frac{1}{|Q|}\int_Q f(x)\ dx ,$$ $ |Q|=\int_Q1\ dx $. One has
$W_p^{n/p}\subset BMO $ for all $1<p<\infty$, and moreover $H_p^{n/p}\subset
H_{p_1}^{n/p_1} $ for $1<p<p_1<\infty $ (this follows from Theorem 2.7.1 of
[Triebel 1]. In some sense $BMO$ is a limit of $H_p^{n/p}$ as $p\to \infty $.

If $f\in W_p^1 $ the restriction of $f$ on hyperplanes $\{ x_n=\epsilon \}\subset
R^n $ (where $(x_1,\cdots,x_n)$ are Euclidean coordinates ) have both $L^p$
and nontangential limits a.e. on $\mathbb{R}^{n-1}=\{ x| x_n=0 \} $ and the
limit function $f|_{\mathbb{R}^{n-1}}$, called trace of $f$, satisfies
$f|_{\mathbb{R}^{n-1}}\in W_p^{1-1/p} $. By a nontangential limit we
mean the following. Let $y\in\mathbb{R}^{n-1} $ and let $C_\delta $ be a
Stolz angle centered at $y$, that is, a set $\{ z,x_n |x_n\ge \delta\cdot
|z| \} $ for $\delta>0$. Then a function $f$ defined in $\mathbb{R}_+^n=\{
x_n>0\} $ has a nontangential limit $f(y)$ at $y$ if $$f(x)
\underset{\underset{x\in C_\delta}{x\to y}} \to f(y) $$ for all $\delta$.
Note that the points in
$C_\delta$ are within a bounded distance from any geodesic of a hyperbolic
metric $$\frac{\sum_{i=1}^n dx_i^2}{x_n^2} ,$$ which has $y$ as a point at
infinity. The trace theorem mentioned above may be found in [Triebel 1],
section 2.7.2. Notice that functions in $W_2^1(\mathbb{R}^2) $ have traces
in $W_2^{1/2}(\mathbb{R}^1) $.

Now let $\Omega\subset \mathbb{R}^n $ be a bounded domain with a smooth
boundary. We define $W_p^s (\Omega) $ as a space of locally integrable
functions with $D^\alpha f\in L^p $ for $|\alpha |\le s $ and such that
$$\sum_{|\alpha |=[s] }\int\int \frac{|D^{\alpha}f(x)-D^{\alpha}f(y) |^p
}{|x-y|^{n+\{ s\}p}} <\infty .$$ Equivalently, $W_p^s(\Omega) $ is a space of
restrictions of function from $W_p^s(\mathbb{R}^n) $ on $\Omega$. [Triebel 1,
Chapter 3]. One also defines $H_p^s(\Omega) $ as a space of restrictions of
$H_p^s(\mathbb{R}^n) $ on $\Omega $. For a compact smooth manifold $M$
without boundary (in particular, for the boundary $\partial \Omega $) one
easily defines the spaces $W_p^s(M) $ and $H_p^s(M)$ [Triebel 1, Chapter
3] ($H_p^s $ is $F_{p,2}^s$ in Triebel's notations ).

If $M$ is compact and $g$ a Riemannian metric on $M$, let $\triangle_g$
be a corresponding Laplace-Beltrami operator. One can construct a space of
Bessel potentials $(1+\triangle )^{-s/2} (L_p(M)) $. It is known
[Rempel-Schulze 1, Theorem 1, section 2.3.2.5], [H\"ormander 1], that this
space coincides with $W_p^s$ (and not $H_p^s$). Warning: our $W_p^s $ is
called $H^{p,s}$ in [Rempel-Schulze 1] and in many other sources. In
particular, $W_2^s(S^1)$ consists of functions $f=\sum_{n\in \mathbb{Z}}a_n
e^{in\theta} $ , such that $\sum |n|^{2s} |a_n|^2 <\infty $. We will see
that $W_2^{1/2}(S^1) $ is especially important in topology.

If $f\in W_p^1(\Omega) $ then $f$ has an $L^p$ and nontangential limit a.e.
on $\partial \Omega $ and $f|_{\partial \Omega}\in W_p^{1-1/p} (\partial
\Omega) $. In particular, for a unit disc $D\subset \mathbb{R}^2 $, and a
function $f\in W_2^1(D)$, $f|_{S^1}\in W_2^{1/2}(S^1)$.

We will need trace theorems for weighted Sobolev-Lorentz spaces [Kudryavcev
1,2], [Vasharin 1], [Lions 1], [Lizorkin 1,2], [Uspenski 1]. Let $\Omega$ be
as above and let $\rho(x)=dist (x,\partial\Omega) $. Consider
$L_p^1(\Omega,\rho^{\alpha}) $ as a space of functions $f$ such that
$\int_{\Omega}|\nabla f|^p\cdot \rho^\alpha \ dx< \infty $. Then $f$ has a
nontangential limit a.e. on $\partial\Omega $ and \\\\
$ 1)\quad f|_{\partial\Omega}=0\mbox{ if }\alpha\le -1 $ \\
$ 2)\quad f|_{\partial\Omega}\in W_p^{\frac{p-1-\alpha}{p}}
(\partial\Omega),\quad\alpha>-1 $.\\\\
Moreover, $$\| f\|_{W_p^{\frac{p-1-\alpha}{p}}}\le const\cdot \int_\Omega
|\nabla f|^p \rho^\alpha \ dx $$ and for any $f\in
W_p^{\frac{p-1-\alpha}{p}}(\partial\Omega) $ and harmonic $h$,
$h|_{\partial\Omega}=f$, one has $$\int_\Omega |\nabla h|^p \rho^\alpha\ dx
\le const\cdot \|f\|_{W_p^{\frac{p-1-\alpha}{p}}} .$$

\section{$l^p$-cohomology of cocompact real hyperbolic lattices}

The following result is an immediate corollary of the Poincar\'e inequality in hyperbolic space, which is equivalent to Hardy inequality, and the classical results on traces of functions in weighted Sobolev spaces, reviewed in the previous section. It first appeared in print, with a different proof, in [Pansu 1]. We include a proof here, as many parts of it will be used in the theory later. 

{\it Theorem 4.1, part 1.}--- Let $G\subset SO^{+}(1,n) $ be a cocompact (uniform)
lattice. Then there is a $G$-equivariant isomorphism of Banach spaces
$$ H^1(G,l^p(G))\simeq W_p^{\frac{n-1}{p}}(S^{n-1})/const $$
for $p>n-1$. For $1<p\le n-1$, $H^1(G,l^p(G))=0$.

\noindent\\
{\it Corollary 4.1.}--- The constant of fine structure $\alpha(G)$ is equal $n-1.$

\noindent\\
{\it Remarks}.\\
1)Theorem 4.1 is a first step in the program of linearization of
3-dimensional topology, which we will develop below in this chapter.
A crucial fact is that $W_2^{1/2}(S^1)$ admits a natural action of the
extended mapping class group $\mathcal{M}ap_{g,1}$. This will be proved in
section 7 below. \\
2)Let $\mathcal{H}^n$ be a hyperbolic space. Since
$G=\pi_1(\mathcal{H}^n/G)$, by the work of 
[Gold\u{s}tein-Kuzminov-Shvedov 1 ] we know
that $H^1(G,l^p(G)) $ equals $L^p$-cohomology of $\mathcal{H}^n$. So Theorem
4.1 computes the $L^p$-cohomology of the hyperbolic space.

Recall that any class $l$ in $H^1(G,l^p(G))$ has a primitive function
$\mathcal{F}:G\to \mathbb{R} $ defined up to a constant, such that
$l(g)=L_gF-F$. This follows from the fact that a module of all functions
$\mathbb{R}^G $ is coinduced from the trivial subgroup and therefore
cohomologically trivial [Brown 1].

\noindent\\
{\it Theorem 4.1, part 2.}--- Let $G$ be a cocompact lattice in $SO^{+}(1,n)$ and
let $l\in H^1(G,l^p(G))$, let $\mathcal{F}:G\to \mathbb{R}$ be a primitive
function for $l$ (unique up to a constant). Let $\partial G\approx S^{n-1} $
be the boundary of $G$ as a word-hyperbolic group. Then for almost all
points $x\in\partial G$, $\mathcal{F}(g)$ has nontangential limits as $g\to
x$, and the limit function $\mathcal{F}|_{S^{n-1}}\in
W_p^{\frac{n-1}{p}}(S^{n-1}) .$

\noindent\\
{\it Corollary 4.2.}--- If $1<p<p_1<\infty $, then a natural map
$H^1(G,l^p(G))\to H^1(G,l^{p_1}(G)) $ is injective. In fact, for
$n-1<p<p_1<\infty $ one has a commutative diagram
\renewcommand{\baselinestretch}{1.5}
\large
$$\begin{array}{rcl}
    H^1(G,l^p(G))& \overset{\sim}{\longrightarrow} &
							W_p^{\frac{n-1}{p}}(S^{n-1})/const \\
    \downarrow & & \downarrow \\
    H^1(G,l^{p_1}(G)) & \overset{\sim}{\longrightarrow} &
                            W_{p_1}^{\frac{n-1}{p_1}}(S^{n-1})/const
    \end{array}
$$
\renewcommand{\baselinestretch}{1.2}
\normalsize
where the right vertical arrows exists by an embedding theorem of
Sobolev-Slobode\u{c}ki space [Triebel 1, 2.7.1 ].

\noindent\\
{\it Proof of the Corollary 4.2.}--- The commutative diagram is
implied by the proof of the Theorem 4.1. The injectivity follows
immediately.

\noindent\\
{\it Proof of the Theorem.}--- Though a shorter proof of part 1 of the
Theorem can be given by using [Goldshtein-Kuzminov-Shvedov 1], in order to
prove part 2 we need to make an isomorphism $H^1(G,l^p(G))\simeq
L^pH^1(\mathcal{H}^n) $ explicit. Here $L^pH^1(V)$ is the $L^p$-cohomology
of a complete Riemannian manifold $V$.

Let $l$ be a cocycle in $Z^1(G,l^p(G))$. We have then an affine isometric
action $ g\overset{\pi}{\mapsto} (v\mapsto L_gv+l(g) )$ of $G$ on $l^p(G)$. To
it corresponds a smooth locally trivial affine Banach bundle over
$M=\mathcal{H}^n/G: E=[\tilde M\times l^p(G) ]/\mbox{\rm{diagonal action}} $.
By local triviality, smooth partition of unity and affine structure on
fibers one constructs a smooth section $s$ of this bundle. It can be
interpreted as an equivalent smooth map $s: \widetilde M\to l^p(G) $, that is, $
s(g^{-1}x)=L_gs(x)+l(g) $. We note that there is some sonstant $C>0$ such
that $\| \nabla s(x)\|<C $ for all $x\in \widetilde M\simeq \mathcal{H}^n $
($\nabla s\in T_x^* \widetilde M\otimes l^p(G)$ ). This is because $M$ is
compact. Now let $\mathcal{F}\in R^G $ be a primitive for $l$, i.e.
$l(g)=L_g\mathcal{F}-\mathcal{F} $. Put $\sigma(x)=s(x)+\mathcal{F} $: this
is a function $\sigma :\widetilde M\to \mathbb{R}^G $ with the same derivative
as $s$ in the sense that for all $g\in G$, $\nabla\sigma_g=\nabla s_g $
where $\sigma_g$, $s_g$ means $g$-th coordinate. Next, we claim that 
$\sigma$ is invariant, i.e. $\sigma(g^{-1}x)=L_g\sigma (x) $. In fact,
$l(g)=L_g\mathcal{F}-\mathcal{F}$, so
$s(g^{-1}x)=L_gs(x)+L_g\mathcal{F}-\mathcal{F}$, so
$\sigma(g^{-1}x)=L_g\sigma(x)$. So for $x\in \tilde M$ and $g,h\in G$ we
have $\sigma(g^{-1}x)(h)=\sigma(x)(g^{-1}h)$. Let $f(x)=\sigma(x)(1)$, then
$\sigma(x)(g)=f(gx)$. Since $\nabla\sigma(x)=\nabla s(x)\in l^p $ and is
bounded in norm, we have for all $x\in\widetilde M$ that $\sum_{g\in G}|\nabla
f(gx)|^p <C $. In particular,
$$ \int_{\widetilde M}|\nabla f|^p =\int_{\widetilde M/G} \sum_{g\in G} |\nabla
f(gx)|^p <C\cdot Vol(M).$$
In other words, $|\nabla f|\in L^p(\mathcal{H}^n) $. Now, we can use a
Poincar\'e model for the hyperbolic space, that is, the unit ball
$B^n\subset \mathbb{R}^n $ with a hyperbolic metric
$$g_h=\frac{g_e}{(1-r^2)^2} .$$ If $\mu_e,\mu_h$ denote a Euclidean and a
hyperbolic measure respectively, $|\nabla f|_e,|\nabla f|_h $ denote a norm
of a gradient of a function in the Euclidean and hyperbolic metric
respectively, $\rho(z)=1-r(z)$ denote a Euclidean distance to the boundary
$\partial B^n\approx S^{n-1}$, then
$$ const_2\cdot\rho^{p-n}\cdot
 |\nabla f|_e^p\cdot \mu_e \le |\nabla f|_h^p\cdot \mu_h \le
const_1\cdot \rho^{p-n} |\nabla f|_e^p \mu_e ,$$ 
so we have $\int_{B^n} \rho^{p-n} |\nabla f|_e^p \mu_e <\infty $.

By a theorem of Kudryavcev-Vasharin-Lizorkin-Uspenski-Lions mentioned above,
we find that $f|_{(1-\epsilon)S^{n-1}} $ has an $L^p$-limit $f|_{S^{n-1}}$
to which it converges nontangentially a.e. , and moreover $f|_{S^{n-1}}\in
W_p^{\frac{n-1}{p}} (S^{n-1}) $ if $p>n-1$ and $f|_{S^{n-1}}=const $ if $p\le
n-1$. We claim that a map $l\mapsto f|_{S^{n-1}} $ is a well-defined bounded
linear operator from $H^1(G,l^p(G))$ to $ W_p^{\frac{n-1}{p}}, p>n-1 $.
First, we notice that since $s:\widetilde M\to l^p(G),
\sigma(x)=s(x)+\mathcal{F} $ and $ \sigma(x)(g)=f(gx) $, we have for almost
all $x\in \widetilde M$, $f(gx)-\mathcal{F}(g)\in l^p(G) $ (as a function of
$g$). In particular, $f(gx)-\mathcal{F}(g)\to 0 $ as $length(g)\to \infty $.
This proves that, identifying $G$ with an orbit of $x$, $\mathcal{F}(g) $
has a nontangential limit a.e. on the boundary $\partial G\approx S^{n-1} $
and $\mathcal{F}|_{S^{n-1}}=f|_{S^{n-1}}\in W_p^{\frac{n-1}{p}}$. In
particular, $f|_{S^{n-1}} \mbox{ mod constants} $ does not depend on the
choice of a section $s$. Since changing $l$ by a coboundary leads to an
isomorphc affine $l^p(G)$-bundle, $f|_{S^{n-1}} \mbox{mod constants}$ depends
only on the class $[l]\in H^1(G,l^p(G)) $. So we get a well-defined operator
$H^1(G,l^p(G))\to W_p^{\frac{n-1}{p}}/const $. We claim it is bounded. An
affine flat bundle $E$ has been defined as $\widetilde M\underset{G}\times
l^p(G) $, where $G$ acts on $l^p(G)$ by $v\mapsto L_gv+l(g) $. It is enough
to show, that there is a constant $C$, depending only on $G$ but not on $l$,
such that $E$ possesses a Lipschitz section $s$ with $ \| \nabla s \| <C\cdot \|
l\| $, where $\| l\|=\sup_{i} \| l(g_i) \| $ for a choice of generators
$g_i, i=1,\cdots,m$. We note that $l$ effectively controls the monodromy of
the flat connection in $E$. A construction of $s$ mentioned above, that is,
a choice of an open covering $\cup U_\alpha =M $, flat sections $s_\alpha $
over $U_\alpha $, a partition of unity $\sum f_\alpha =1 $ with $supp
f_\alpha \subset U_\alpha $, so that $s=\sum f_\alpha s_\alpha $, gives a
bound of $|\nabla s| $ in terms of monodromy, as desired.

We note that by [Gold\u{s}tein-Kuzminov-Shvedov 1],
$H^1(G,l^p(G))=L^pH^1(\mathcal{H}^n )$, so to any class in $H^1(G,l^p(G))$
we have associated a function $f$ such that $df $ is in $L^p$, or,
equivalently, $\int_{\mathcal{H}^n}|\nabla f|_h^p \mu_h <\infty $. What we
in fact did above was an explicit construction of this correspondence between
$l^p$- and $L^p$-cohomology.

So far we have constructed a bounded operator $H^1(G,l^p(G))\to
W_p^{\frac{n-1}{p}}(S^{n-1}) $, $p>n-1$. We wish to show that this operator
is in fact an isomorphism of Banach spaces. To this end, we will need a
Poincar\'e inequality in hyperbolic space.

\noindent\\
{\it Proposition 4.5 }(Poincar\'e inequality in $\mathcal{H}^n$).---
Let $f$ be a locally integrable measurable function with
$\int_{\mathcal{H}^n} |\nabla f|_h^p\ d\mu_h < \infty $. Then\\
1)If $p\le n-1 $, then $\int_{\mathcal{H}^n} |f-c|^p \ d\mu_h <\infty $
for some constant $c$;\\
2)If $p>n-1 $ and $f|_{S^{n-1}} $ as an element of
$W_p^{\frac{n-1}{p}}(S^{n-1}) $ is zero, then $\int_{\mathcal{H}^n}|f|^p\
d\mu_h<\infty $.

\noindent\\
{\it Proof.}--- A much more general theorem is contained in [Strichartz 1].

\noindent\\
We now claim that $H^1(G,l^p(G))=0$ for $p\le n-1 $. This in fact follows
immediately from $H^1(G,l^p(G))=L^pH^1(\mathcal{H}^n) $ 
[Gold\u{s}tein-Kuzminov-Shvedov 1 ] and Proposition 4.5. Now, if $p>n-1$, then we claim that the operator
$H^1(G,l^p(G))\to W_p^{\frac{n-1}{p}}(S^{n-1})/const $ constructed above is
injective. In fact, if $f|_{S^{n-1}}=0$, then by Proposition 4.5, $f\in
L^p(\mathcal{H}^n,\mu_h)$, so $\int_M \Sigma_g |f(gx)|^p\ d\mu_h <\infty $,
so for almost all $x\in \tilde M$, $\sum_g |f(gx)|^p<\infty $. But
$f(gx)-\mathcal{F}(g)\in l^p(G)$, so $\mathcal{F}\in l^p(G)$ and $[l]=0$.
Now, if $h\in W_p^{\frac{n-1}{p}}(S^{n-1})$, we denote by $H$ its harmonic
extension into $B^n$. Then [Uspenski 1], [Lizorkin 2], $\int_{B^n}
\rho^{p-n} |\nabla H|^p\ d\mu_e < \| h\|_{W_p^{\frac{n-1}{p}}(S^{n-1})} $,
so $dH$ is an $L^p$ 1-form on $\mathcal{H}^n$. This shows that the injective
operator $H^1(G,l^p(G))=L^pH^1(\mathcal{H}^n)\to
W_p^{\frac{n-1}{p}}(S^{n-1})/const $ has a bounded right inverse, so it is
an isomorphism by Banach theorem. This proves Theorem 4.1.

\noindent\\
{\it Corollary 4.5.}--- Let $G$ be a cocompact lattice in $SO^{+}(1,n)$
and let $\mathcal{F}:G\to \mathbb{R} $ be such that
$L_g\mathcal{F}-\mathcal{F}\in l^p(G) $, for all $g\in G\quad (p>n-1) $.
Then the limit function $u=\mathcal{F}|_{S^{n-1}} $ belongs to $L^q
(S^{n-1}) $ for all $q>1$. In fact,
$$ \sup_{1<q<\infty} \left(\frac{n-1}{q}\right)^{1/q'} \| u\|_{L^q(S^{n-1})}<\infty $$
Moreover, $u$ is in the linear hull of all functions $f$ satisfying
$$ \int_{S^{n-1}}exp(|f|^{p'})<\infty .$$

\noindent{\it Proof.} is an immediate corollary of Theorem 4.1 and the
properties of the Orlicz space $L_\infty(\log L )_{-a} $ and the fact that
$W_p^{(n-1)/p} (S^{n-1})\subset L_\infty(\log L)_{-a}(S^{n-1}) $ for $a\ge
1/p' $, (see [Edmunds-Triebel 1]).

\noindent\\
We will use this corollary in a sequel to this paper [Reznikov 10] in analyzing the local behaviour of the Cannon-Thurston Peano
curves, corresponding to fibers of the hyperbolic 3-manifolds, fibered over
the circle.
\noindent\\
\section{Growth of primitives for $l^p$-cocycles on the surface group}
{\it Theorem 5.1.}--- Let $G$ be a cocompact lattice in $SO(2,1)$
and let $\mathcal{F}:G\to \mathbb{R} $ be such that for all $g\in G$,
$L_g\mathcal{F}-\mathcal{F} \in l^p(G) $, $p>1$. Then for any word length on
$G$,
$$ |\mathcal{F}(g)|\le const\cdot [length(g)]^{1/p'} .$$
{\it Proof} follows from Theorem 4.1 and a following lemma.

\noindent\\
{\it Lemma 5.2.}--- Let $u$ be a harmonic function in the unit disc
such that $u|_{S^1}\in W_p^{1/p}(S^1) $. Then
$$ |u(z)|\le const\cdot[\log (1-|z|)]^{1/p'} .$$

\noindent\\
{\it Proof of the lemma.}--- Here we only treat the case $p=2$. The
full proof will be given in Section 11. Let
$u(e^{i\theta})=\sum_{n\in\mathbb{Z}}a_ne^{in\theta} $. Since $
(1+\triangle )^{1/4} u\in L^2 $, we have $\{ |n|^{1/2}a_n\}\in
l^2(\mathbb{Z})$, therefore for $|z|<1\quad (b_n=|a_n|+|a_{-n}| ) $.
\renewcommand{\baselinestretch}{1.5}
\large
$$ \begin{array}{rl}
    u(z)-a_0&\le \sum_{n>0} (|a_n|+|a_{-n}|)|z|^n \\
        &=\sum |n|^{1/2} b_n\cdot \frac{1}{|n|^{\frac{1}{2}}} |z|^n \\
        &\le (\sum |n|b_n^2)^{1/2}\cdot(\sum \frac{1}{|n|} |z|^{2n})^{1/2}\\
        &\le const\cdot [\log (1-|z|)]^{1/2}.
    \end{array}
$$
\renewcommand{\baselinestretch}{1.2}
\normalsize

\section{Embedding of negatively curved manifolds and the boundaries of
their universal covers}

A problem of fundamental importance in topology is the following: let
$M^m\overset{\varphi}{\hookrightarrow} N^n $ be a smooth $\pi_1$-injective
embedding of manifolds of nonpositive curvature. Let $\widetilde \varphi: 
\widetilde
M\to \widetilde N$ be a lift of $\widetilde\varphi $.
Is there a limit map $S^{m-1}\approx
\partial\widetilde M\overset{\partial\widetilde\varphi}{\to}
\partial\widetilde N\approx
S^{n-1} $ and if there is, how smooth it is? For instance, let $N^3$ be a
compact hyperbolic 3-manifold, and $M^2$ be an incompressible embedded
surface in $N^3$. Then there always exists a limit continuous map
$S^1\overset{\partial\widetilde\varphi}{\to}S^2 $. Moreover, if $M$ is not a
virtual fiber of a fibration over the circle, then
$\partial\widetilde\varphi(S^1)$ is a quasifuchsian Jordan curve. If $M$ is a
virtual fiber, then $\partial\widetilde\varphi: S^1\to S^2 $ is a Peano
curve in a sense that its image fills $S^2$ [Cannon-Thurston 1]. This deep dichotomy follows from the result of
[Bonahon 1]. We have a following very general theorem 9the embedding condition is superfluous but makes the proof more transparent):

\noindent\\
{\it Theorem 6.1.}--- Let $M^m\overset{\varphi}{\hookrightarrow}N^n $
be a smooth $\pi_1$-injective embedding of complete Riemannian manifolds, of
pinched negative curvature. Suppose $M$ is compact. Let $\widetilde \varphi
:\widetilde M\to \widetilde N $ be a 
lift of $\widetilde \varphi $. Let $p_0\in \widetilde
N$ and $\pi: \widetilde N\backslash \{p_0\} 
\to S^{n-1}(T_{p_0}\widetilde N) $ be  a
radial geodesic projection of $\tilde N\backslash \{ p_0\} $ onto the unit
tangent sphere. Identify $T_{p_0}\widetilde N$ with $\mathbb{R}^{n} $. Let
$q_0\in \widetilde M$. Then: \\
1) For almost all unit tangent vectors $v\in T_{q_0}(\widetilde M) $, the
restriction of $\pi \widetilde \varphi $ on a geodesic $\gamma (q_0,v) $ starting at
$q_0$ with a tangent vector $v$ has an $L^1$-derivative as a map
$\widetilde\varphi|_{\gamma(q_0,v)}: \mathbb{R}_+ \to \mathbb{R}^{n} .$ \\
2) For almost all $v$ there exists a limit $\lim_{t\to \infty}
\pi \widetilde\varphi [\gamma(q_0,v)(t)] $. \\
3) The resulting measurable map $\partial \widetilde M\approx
S^{m-1}\overset{\partial\widetilde\varphi}
{\to}S^{n-1}\approx \partial\widetilde N $
does not depend on the choice of $p_0,q_0$. \\
4) If both $M$,$N$ are (real) hyperbolic, then for any $p>n-1$,
$\partial\widetilde\varphi $ induces a bounded linear operator
$$ \partial\widetilde\varphi_*: W_p^{\frac{n-1}{p}}(S^{n-1})\to
W_p^{\frac{m-1}{p}}(S^{m-1}). $$\\
5) If $M$ is hyperbolic and $-K\le K(N) \le -1$, then for $p>(n-1)\sqrt{K}$, $\partial\widetilde\varphi $ induces a bounded linear operator
$$ \partial\widetilde\varphi_*: C^{\infty} (S^{n-1})\to
W_p^{\frac{m-1}{p}}(S^{m-1}) $$ for $p>(n-1)\sqrt{K}$.

\noindent\\
{\it Theorem 6.2.}---
 Let $N^3$ be a compact oriented hyperbolic 
three-manifold,
let $M^2\overset{\varphi}{\hookrightarrow} N^3 $ be an incompressible
immersed  surface, and let $x_1,x_2,x_3$ be Euclidean coordinates on
$S^2\approx\partial N^3 $. Then \\
1) If $\partial \widetilde \varphi $ is quasifuchsian, 
then $x_i\circ\partial\widetilde\varphi
:S^1\to \mathbb{R} $ are in $W_p^{1/p}$ for $p\ge 2.$ \\
2) If $M^2$ is a virtual fiber then $x_i\circ 
\partial\widetilde\varphi : S^1\to
\mathbb{R} $ are in $W_p^{1/p} $ for $p>2$ (but probably not in $W_2^{1/2} $ ).

\noindent\\
{\it Proof of the Theorem 6.1.}--- We will assume $-k\le K(M) \le -1,
-K\le K(N) \le -1 $. For $x\in \widetilde N $ let $r(x)=\rho(p_0,x)$.

\noindent\\
{\it Lemma 6.3.}--- For $r_0>0$ and $r(x)>r_0 $, $|\nabla\pi(x)|\le
const(r_0)e^{-r(x)} $, where we view $\pi$ as a map $N\backslash \{p_0\} \to
\mathbb{R}^n $.

\noindent\\
{\it Proof} is an immediate application of the comparison theorem,
mentioned above in the proof of Proposition 2.1.

\noindent\\
{\it Lemma 6.4.}---
$$ \int_{\widetilde N\backslash B(p_0,r_0) } |\nabla \pi(x) |^p <\infty \mbox{
for } p>(n-1)\sqrt{K} .$$

\noindent{\it Proof} repeats the argument in the proof of Proposition  2.1.\\

Now consider a tubular neighbourhood of $M$ in $N$. There exists an embedding
of $M\times I\to N$ where $I=[-1,1]$. Moreover, the restriction of the
metric $g_N$ of $N$ onto $M\times I $ is equivalent to the product metric
$g_M+dx^2 $ (we say two metrics are equivalent if each one is bounded above
by another one times a constant). It follows that there is an embedding
$$ \widetilde M\times I \overset{\Phi}{\to}\widetilde N$$
such that $g_{\widetilde N} |\widetilde M\times I $ is 
equivalent to $g_{\widetilde M}
+dx^2 $. Since $\varphi $ is $\pi_1$-injective, for any $r_0>0$ there is
$r_1>0$ such that if $\rho_M(q_0,z)>r_1 $, then $\rho_N(p_0,\Phi (z,t))>r_0
$ for $t\in [-1,1]$. It follows that
$$ {\int\int}_{\widetilde M\backslash B(q_0,r_1)\times I } 
|\nabla \pi\circ\Phi |^p\ dVol(\widetilde M) dt <\infty $$
Therefore for almost all $t_0\in I$,
$$ \int_{\widetilde M\backslash B(q_0,r_1)} |\nabla (\pi\circ\Phi
 (z,t_0))|^p\ dVol(\widetilde M) <\infty $$
Fix such $t_0$ and let 
$f=\pi\circ\Phi (z,t_0) :\widetilde M\backslash B(q_0,r_1)\to
\mathbb{R}^n $. We know that
$$ \int_{\widetilde M\backslash B(q_0,r_1) } 
|\nabla f|^p\ d Vol(\widetilde M) <\infty
.$$
Expressing the integral in polar coordinates and taking into account that
$K(M)\le -1 $ we have
$$ \int_{S^{m-1}(T_{q_0}\widetilde M) } dv 
\int_{r_1}^\infty e^{(m-1)t} |\nabla
f|^p \ dt <\infty .$$
In particular, for almost all $v\in S^{m-1}(T_{q_0}\widetilde M)$,
$$ \int_{r_1}^\infty e^{(m-1)t} |\frac{\partial f}{\partial t}|^p \ dt
<\infty .$$
In other words, for such $v$, $|\frac{\partial f}{\partial t} |\cdot
e^{\frac{(m-1)}{p}t}\in L^p[r_1,\infty] $, therefore $|\frac{\partial
f}{\partial t} |\in L^1[r_1,\infty] $, since $e^{-\frac{m-1}{p}t}\in
L^{p'}[r_1,\infty] $. This proves 1). The statements 2) and 3) follow
directly.

Now suppose $K(M)=K(N)=-1$. Let $u\in W_p^{\frac{n-1}{p}}(S^{n-1}) $,
$p>n-1$. Then a harmonic extension $g$ of $u$ satisfies
$$ \int_{\widetilde N}|\nabla g |^p<\infty $$
as we know from [Lizorkin 1], [Uspenski 1] and the proof of Theorem 5.1. By
the argument above, there is a $t_0\in I$, such that the composite function
$g\circ \Phi(z,t_0) $ satisfies
$$ \int_{\widetilde M} |\nabla (g\circ\Phi(z,t_0)|^p <\infty $$
But then the trace $g\circ\Phi(z,t_0)|_{\partial\widetilde M} $ lies in
$W_p^{\frac{m-1}{p}}(S^{m-1}) $. This proves part 4) of Theorem 6.1. A proof of part 5) is identical. The Theorems 6.1 and 6.2,2) are proved. To prove Theorem 6.2, 1) we notice that a restriction of any function $u\in  W_2^{1}(S^2)$ on a quasicircle belongs to the class  $W_2^{\frac{1}{2}}(S^1 )$. This follows immediately from the invariance of  $W_2^{1}(S^2)$ under quasiconformal homeomorphisms, and a fact that functions from  $W_2^{1}(B^2)$ have traces in  $W_2^{\frac{1}{2}}(S^1 )$(notice that the Dirichlet energy of a function of two variables is an invariant of the conformal class of a metric).

\noindent\\
As the reader has noticed, we could assume $\pi_1(M)=\pi_1(N) $, so that
$\pi_1(M)$ acts discretely in $\widetilde N$ and 
$N=\widetilde N/\pi_1(M) $. On the
other hand, the proof does not use the fact that $M$ is embedded,
 so the Theorem 6.1 stays true for (finite-to-one)immersions in $N$. 

We will outline
now, having in mind the applications in the sequel to this paper, how to study the limit maps from the point of view of ergodic theory. The results thus obtained are weaker then those proved above, but apply to non-discrete representations. Our treatment can be seen as a development of a vague remark of [Thurston 1, 6.4.4].
Let $M^m$ be a compact hyperbolic mainfold, $\widetilde
N=\mathcal{H}^n $ and $\rho: \pi_1(M)\to Iso(\widetilde N) $ a discrete faithful
representation. Let $N=\widetilde N/\rho(\pi_1(M)) $. We would like to study a
boundary map $\partial\widetilde \varphi: 
\widetilde M\to \widetilde N $ where $\varphi$
is a smooth map $M\to N$, inducing $\rho$.

\noindent\\
{\it Lemma 6.5.}--- There exists a $\pi_1(M)$ equivariant measurable
map $\psi$ from $\partial\widetilde M=S^{m-1} $ to the space of probability
measures on $\partial\widetilde N=S^{n-1} $.

\noindent\\
{\it Proof.}--- For any compact Riemannian manifold $M$, any compact
metric space $X$ and any representation $\rho: \pi_1(M)\to Homeo(X) $, there
is a $\pi_1(M)$-equivariant harmonic function from 
$\widetilde M$ to the affine
space of charges on $X$, taking values in probability measures. This simple
fact is various degrees of generality has been shown in [Furstenberg 1], [L.Garnett 1], [Kaimanovich-Vershik 1]. If $M$ is hyperbolic,
then the Poisson boundary of $\widetilde M$ is $\partial\widetilde M$, and the result
follows.\\

Now let $\psi_0+\psi_c$ be the decomposition of $\psi$ into atomic and
non-atomic parts. Obviously, $\psi_c$ is also $\pi_1(M)$-equivariant. We
claim $\psi_c=0$. First, $\int \psi_c $ is a $\pi_1(M)$-invariant function
on $\partial\widetilde M=S^{n-1}$, whence a constant, since $\pi_1(M)$ acts on
$S^{n-1}$ ergodically. So if $\psi_c\not= 0$ we may assume $\psi_c$ is a
probability measure. Second, let $G$ be a center of gravity map from the
nonatomic measures on $\partial\widetilde N$ to $N$
[Furstenberg 2]. Then $G\circ \phi_c $ is a $\pi_1(M)$-equivariant
map from $\partial\widetilde M$ to $N$. In particular, $\rho(G\circ
\psi_c (x), G\circ\psi_c(y))$ is a $\pi_1(M)$-invariant function on
$\partial\widetilde M\times \partial\widetilde M$ 
whence a constant by [Hopf 1] and
[Sullivan 3]. It follows easily that $G\circ\psi_c=const $ which is
impossible since $\rho$ is discrete. So $\psi_c=0 $.

We deduce that $\psi$ is atomic, $\psi(z)=\sum_{i=1}^\infty m_i\delta
(\psi_i(z)) ,m_1\ge m_2\ge \cdots $. 
Though $\psi_i(z): \partial\widetilde M\to
\partial \widetilde N $ are not uniquely defined, 
$m_i:\partial\widetilde M\to
\mathbb{R} $ are. It follows that $m_i$ are $\pi_1(M)$-invariant, whence
constant. Since $\sum m_i=1 $, there is some $i$ such that $m_{i+1}<m_1 $.
Choose first such $i$. Then $m_1=\cdots =m_i $ and we get a measurable
equivariant map
$$\partial\widetilde M=S^{m-1}\to \underbrace{S^{n-1}\times\cdots\times
S^{n-1}}_i /S_i ,$$
where $S_i$ is the symmetric group in $i$ letters.

So far we did not use the fact that $\rho$ is discrete, but only that
$\rho(\pi_1(M))$ does not have fixed points in $\widetilde N=\mathcal{H}^m $.
So:

\noindent\\
{\it Propostion 6.6.}--- Let $M^m$ be a compact hyperbolic manifold
and let $\rho: \pi_1(M)\to SO^{+}(1,n) $ be such that $\rho (\pi_1(M))$ does not
have fixed points in $\mathcal{H}^n$. Then there exists a
$\pi_1(M)$-equivariant measurable map
$$ S^{m-1}=\partial\widetilde M\overset{\psi}{\to} 
(\mbox{\it{subsets of cardinality }}i
\mbox{\it{ of }}S^{n-1}=\partial\mathcal{H}^n ) $$
for some $i\ge 1$.\\

Using cross-ratios and the ergodicity of the action of $\pi_1 (M)$ on 
$\partial\widetilde M\times \partial\widetilde M$ , one can easily show $i=1$. Now to any $x\in\widetilde M$ one associates a Poisson measure $\mu_x$ on
$S^{m-1}$. Its pushforward $\psi_*\mu_x $ is a probability measure on
$S^{n-1} $. The pushforward of a measure by a multivalued map is defined by
$$ \int_{S^{n-1}} f\ d[\psi_*\mu ]=\int_{S^{m-1}}\sum_{y\in\psi(x)}f(y)\
d\mu ,$$
where $f\in C(S^{n-1})$.

Now under some natural conditions $\psi_*\mu_x $ does not have atoms and
using the baricenter map $G$ in $\mathcal{H}^n $ one can define
$s(x)=G(\psi_*\mu_x )$. This can easily be shown to be continuous
equivariant map $\widetilde M\overset{s}{\to} \mathcal{H}^n $, again under some
natural assumption on $\rho$. The multivalued map $\psi$ should be regarded
as a weak radial limit of $s$, but we will not pursue this point any further.

\section{The action of quasisymmetric and quasiconformal homeomorphisms on
$W_p^{(n-1)/p} $}

A well known result [Reimann 1] characterizes quasiconformal maps between
domains $D_1,D_2$ in $\mathbb{R}^n,n>2$ as those which induce an isomorphism
of Banach spaces $BMO(D_1)$ and $BMO(D_2)$. We will see now that this result
in case $D_1=D_2=\mathbb{R}^n $ is a limit as $p\to \infty$ of the following
result which establishes a quasiconformal invariance of fractional Sobolev spaces
$W_p^{n/p} $. Of special importance is the fact that the result holds for
$n=1$ and quasisymmetric homeomorphisms of $S^1$. The proof of the following lemma is ``almost'' contained in remarks made in [Pansu 1--3]. 

\noindent\\
{\it Lemma 7.1.}--- Let $\mathcal{G}_{n-1},n\ge 2 $ be a group of
quasisymmetric ($n=2$) or quasiconformal ($n\ge 3$) homeomorphisms of
$S^{n-1}$. Then for any $p>1 \quad (n=2) $ or $p\ge n-1\quad (n\ge 3) $,
$\mathcal{G}_{n-1}$ leaves invariant a Sobolev-Slobode\u{c}ki space
$W_p^{n-1/p}(S^{n-1})$. For any $\Phi\in \mathcal{G}_{n-1} $, the
corresponding map
$$ \Phi_* :W_p^{n-1/p}(S^n)\to W_p^{n-1/p}(S^{n-1}) $$
is an automorphism of the Banach space $W_p^{n-1/p}(S^{n-1}) $.

\noindent\\
{\it Theorem 7.2.}--- There exists for any $n\ge 2$ a bounded
antisymmetric polylinear map
$$ \underbrace{W_n^{\frac{n-1}{n}}(S^{n-1})/const\times \cdots\times
W_n^{\frac{n-1}{n}}(S^{n-1})/const}_n\to \mathbb{R}\ , $$
defined on the smooth functions by $f_1,\cdots,f_n\to \int_{S^n}f_1\
df_2\cdots df_n $, which is invariant under $\mathcal{G}_{n-1}. $

In particular, we have

\noindent\\
{\it Corollary 7.3.}--- There exists a representation
$$ \mathcal{G}_1\to Sp (W_2^{1/2}(S^1)/const ) ,$$
defined by $\Phi(f)=f\circ \Phi^{-1} $.

\noindent\\
{\it Proof of the Lemma  7.1.}--- We will need a result, proved for
$n=2 $ in [Ahlfors-Beurling 1] for $n=3$ in [Carleson 1] and for $n\ge 4$ in
[Tukia-V\"ais\"al\"a 1]:

\noindent\\
{\it Theorem.}--- Let $\phi :S^{n-1}\to S^{n-1} $ be quasisymmetric
($n=2$) or quasiconformal ($n\ge 3$). Then there exists an extension 
$\widetilde
\phi $ of $\phi $ as a homeomorphism of $B^n$, which is a quasiisometry of
the hyperbolic metric:
$$ const_2\cdot g_h \le \widetilde \phi_* g_h \le const_1\cdot g_h . $$
Now let $f\in W_p^{\frac{n-1}{p}}\quad (p>n-1) $. Let $u$ be a harmonic
function in $B^n$, extending $f$. We know that
$$ \int |\nabla u|_h^p\ d\mu_h \le const_3 \|
f\|_{W_p^{\frac{n-1}{p}}(S^{n-1})} .$$
It follows that
$$ \int |\nabla (u\circ\widetilde\phi )|_h^p\ d\mu_h \le const_4 \|
f\|_{W_p^{\frac{n-1}{p}}} <\infty ,$$
and by the trace theorem,
$$ \| u\circ \widetilde\phi \|_{W_p^{\frac{n-1}{p}}}\le
const_5\|f\|_{W_p^{\frac{n-1}{p}}} ,$$
which proves the theorem for $p>n-1$. For $p=n-1, n\ge 3$, the result is
standard.

\noindent\\
{\it Proof of the Theorem 7.2.}--- Let $f_1,\cdots,f_n\in
W_n^{\frac{n-1}{n}}(S^{n-1}) $. Let $u_i$ be a harmonic extension of $f_i$.
The result follows at once from the formula $$\int_{S^{n-1}} f_1\ df_2 \cdots
df_n =\int_{B^n}du_1du_2\cdots du_n $$. Since $\int |\nabla u_i |_h^n\
du_h <\infty $, the integral $\int_{B^n}\ du_1\cdots du_n $ is finite
by H\"older inequality. The invariance is obvious.

\noindent\\
{\it Proof of the Corollary 7.3.}--- A formula
$<f_1,f_2>=\int_{S^1}f_1 \ df_2 $ gives $W_2^{1/2}/const $ a structure of a
symplectic Hilbert space. This means that a map $$W_2^{1/2}/const\to
(W_2^{1/2}/const)^* $$ given by $f\to <f,\cdot > $ is an isomorphism (not
isometry) of Hilbert spaces. By $Sp (W_2^{1/2}/const) $ we mean a group of
invertible bounded operators which leaves this symplectic form invariant.
The result now follows from Lemma  7.1 and Theorem 7.2 .

\section{Boundary values of quasiconformal maps and regularity of
quasisymmetric homeomorphisms}

\noindent{\it Proposition 8.1.}--- 
Let $\phi$ be a quasiconformal map, defined
in a neighborhood of the unit ball $B^n$. Then $\phi|_{S^n} $ as a map
$S^n\to \mathbb{R}^n $ belongs to a class $W_n^{\frac{n-1}{n}+\delta} $
for some $\delta >0$. In particular if $n=2$ and
$\phi(e^{i\theta})=\sum_{n\in\mathbb{Z}} a_n e^{in\theta} $ then $\sum
|n|^{1+\delta} |a_n|^2 <\infty $. If $\phi$ is just defined in $B^n$ then
for almost all $\alpha\in S^n$ there exists a limit $\lim_{r\to 1} \phi(rx)
$ and $\phi|_{S^{n-1}}\in W^{\frac{n-1}{n}} $. In particular, for $n=2$ and
$\phi(e^{i\theta})=\sum_{n\in \mathbb{Z} } a_n e^{in\theta} $, $\sum |n|
|a_n|^2 <\infty $.

\noindent\\
{\it Remark.}--- The last statement for conformal maps is the
"Flachensatz".

\noindent\\
{\it Proof.}--- Since $\phi$ as a map $B^n\to \mathbb{R}^n $ belongs
to $W_n^1$, the last statement follows immediatedly from the trace theorem.
To prove the first, recall that $\phi$ is locally in
$W_{n+\delta'}^1,\delta'>0 $ [Bojarski 1], [Gehring 2]. Therefore
$\phi|_{S^{n-1}}\in W_n^{\frac{n-1}{n}+\delta } $, again by the trace
theorem.

\noindent\\
{\it Theorem 8.2.}--- Let $\varphi :S^1\to S^1 $ be a quasisymmetric
homeomorphism. Then as a map $S^1\to \mathbb{R}^2 , \varphi\in
W_p^{1/p+\delta(p)}, \delta(p)>0, $ for all $p>1$. If $\varphi
(e^{i\theta})=\sum_{n\in \mathbb{Z}} a_n e^{in\theta} $, then
$\sum_{n\in\mathbb{Z}} |n|^{p'/p+\delta }|a_n|^{p'} <\infty $ for all
$1<p\le 2 $.

\noindent\\
{\it Proof.}--- Let $\Phi: D^2\to D^2 $ be a quasiisometry of the
hyperbolic plane, extending $\varphi $. We know  that
$\Phi,\Phi^{-1} $ are H\"older in Euclidean metric. Let $f$ be a smooth
function defined in a neighbourhood of $D^2$. Then for $p>1$
$$ \int_{D^2} |\nabla f|_h^p \rho_e^{-\epsilon} (x,\partial D^2)\cdot \
d\mu_h <\infty $$
for $\epsilon>0$ small enough (one needs $\epsilon<p-1 $).

Since $\Phi$ is a quasiisometry for the hyperbolic metric and biH\"older for
the Euclidean metric, we have for $g=f\circ \Phi $:
$$ \int_{D^2} |\nabla g|_h^p \rho_e^{-\beta} (y,\partial D^2)\ d\mu_h
<\infty $$
for some $\beta>0$. Rewriting in Euclidean terms, we have
$$ \int_{D^2} |\nabla g|_e^p\cdot [\rho (y,\partial D^2)]^{p-\beta-2
}<\infty ,$$
therefore $g|_{S^1}\in W_p^{\frac{1}{p}+\delta } $ by the trace theorem for
weighted Sobolev spaces. Letting $f$ be an Euclidean coordinate function, we
get $\varphi\in W_p^{\frac{1}{p}+\delta} $. The last statement follows from
Young-Hausdorff theorem.

\noindent\\
{\it Remark 8.3.}--- It had been a famous problem in fifties if $\varphi $
is absolutely continuous ( that is, in $W_1^1$). Though the answer is well-known to be negative, we
see that $\varphi $ is as close to be absolutely continuous  as one wishes.
We will use Theorem 8.2 in a sequel to this paper to prove the existence of the vacuum vector for quantized moduli space for $p>1$. We also notice that the argument above together with the proof of Theorem 6.1 shows the following: if $\varphi : M\to N$ is an $\pi_1$-injective immersion of hyperbolic manifolds, $M$ compact, such that for $g\in \pi_1 (M)$ and some fixed $z_0 \in \tilde N$
$$\rho (z_0 , \varphi_{\ast} (g)z_0 )\ge const\cdot length(g),$$
then $\partial \tilde \varphi$ is of class $W_p^{(m-1)/p + \delta}$ and therefore H\"older continuous. It is not enough, though, to prove a continuity if the Cannon-Thurston curve. See [Reznikov 10] for futher study.
\section{Teichm\"uller spaces and quantization of the mapping class group, I}

We denote $\mathcal{M}ap_g$ the mapping class group of genus $g$ and
$\mathcal{M}ap_{g,1} $ the extended mapping class group. If $\Sigma^g $ is a
closed oriented surface of genus $g$, $\Gamma_g =\pi_1 (\Sigma^g) $ then
$\mathcal{M}ap_{g,1}=Aut (\Gamma_g) $ and one has an exact sequence
$$ 1\to \Gamma_g\to \mathcal{M}ap_{g,1} \to \mathcal{M}ap_g \to 1 .$$
{\it Proposition  9.1}(Quantization of the moduli space)--- For any $p>1$ there
exists a representation
$$ \mathcal{M}ap_{g,1}\overset{\pi_p}{\to} Aut(W_p^{1/p} (S^1)/const ) $$
given by the formula
$$ \pi_p(\varphi)(f)=f\circ\Phi^{-1} ,$$
where $\Phi$ is a quasisymmetric homeomorphism of $S^1$, induced by
$\varphi$ and a choice of a hyperbolic structure in $\Sigma^g $. For $p=2$ the
representation
$$ \pi_2: \mathcal{M}ap_{g,1}\to Aut(W_2^{1/2}(S^1)/const) $$
is symplectic, that is, $\pi_2 (\mathcal{M}ap_{g,1})\subset
Sp(W_2^{1/2}(S^1)/const ) .$

\noindent\\
{\it Proof.}--- Fix a hyperbolic structure on $\Sigma$. Then by a
classical theorem of Nielsen, one gets a representation
$\mathcal{M}ap_{g,1}\to \mathcal{G}_1 $. The theorem now follows from
Theorem 7.1.

Now let $G \overset{\pi_0}{\longrightarrow} 
PSL_2 (\mathbb{R}) $ be a Fuchsian group,
possibly infinitely generated. We recall that a Teichm\"uller space $\mathbf{T}(G)$ is
defined as follows : points of $\mathbf{T}(G) $ are discrete representation
$G\overset{\pi}{\longrightarrow}PSL_2 (\mathbb{R}) $, defined up to conjugation by an
element of $PSL_2(\mathbb{R})$, which are quasiconformally conjugate to
$\pi_0$, that is, there is a quasisymmetric homeomorphism $\Phi$ of $S^1$
such that $\pi=\Phi\circ\pi_0\circ\Phi^{-1} $. Notice that this
definition is equivalent to the standard one by a result of [Douady-Earle 1].

\noindent\\
{\it Corollary 9.2.}--- Let $\pi_0,\pi$ be two discrete representation
of a group $G$. Then if $\pi$ lies in the Teichm\"uller space of $\pi_0$,
then the unitary representations
$$ G\overset{\pi_0}{\to}PSL_2(\mathbb{R})\overset{\beta}{\to} U(W_2^{1/2}(S^1)/const ) $$
and
$$ G\overset{\pi}{\to} PSL_2(\mathbb{R})\overset{\beta}{\to}
U(W_2^{1/2}(S^1)/const )$$
are unitarily equivalent.

\noindent\\
{\it Remark 1.}--- The fact that $PSL_2(\mathbb{R})$ acts in
$W_2^{1/2}(S^1)/const $ by unitary operators (with respect to the complex structure given by the Hilbert transform) is well-known [Nag 1]. In fact,
this unitary representation belongs to the discrete series and may be
realized in $L^2$ holomorphic 1-forms in $B^2$.

\noindent\\
{\it Proof.}--- Since $\pi=\Phi\circ\pi_0\circ\Phi^{-1} $ and
$\mathcal{G}_1$ act in $W_2^{1/2}(S^1)/const $, we get an invertible
operator $A$ such that $\beta\circ\pi =A\ \beta\circ\pi_0\ A^{-1} $. By polar
decomposition $A=UP$ where $P$ is positive self-adjoint, $U$ is unitary, $P$
commutes with $\beta\circ\pi_0$ and $U$ intertwines $\beta\circ\pi_0$ and $\beta\circ\pi$, as desired.

The following special case is very important. Let $\pi_0 :G\to
PSL_2(\mathbb{R}) $ be a Fuchsian group corresponding to a Riemann surface
of finite type (that is, a torsion-free lattice in $PSL_2(\mathbb{R})$). Let
$\Sigma= \mathcal{H}^2 /G $ and let $\varphi\in 
\mathcal{M}ap(\Sigma, x_0), x_0\in
\Sigma$. Let $\Phi$ be a quasisymmetric homeomorphism of $S^1$ which is
the boundary value of a quasiconformal homeomorphism $\Psi$  of $(\Sigma,x_0)$,
representing $\varphi$. Then
$$ \pi_0 \circ \varphi^{-1} =\Phi \pi_0 \Phi^{-1} .$$
Let $A_\varphi :W_2^{1/2} /const\to W_2^{1/2} /const $ be an invertible
operator, representing $\varphi$. Let $P_\varphi^2=A_\varphi^* A_\varphi
$.Then $P_\varphi$ commutes with $\beta\circ\pi_0$. We obtained the following

\noindent\\
{\it Theorem 9.2.}--- Let $\pi_0: G\to PSL_2(\mathbb{R}) $ be
a torsion-free lattice. Let $\Sigma =\mathcal{H}^2/G, x_0\in\Sigma,
\varphi\in\mathcal{M}ap (\Sigma,x_0)$, $\Psi$ a quasiconformal homeomorphism
inducing $\varphi$, $\Phi$ the  trace of its lift to $\mathcal{H}^2$ on $S^1$, $A_\varphi $ an invertible
opreator in $W_2^{1/2}(S^1)/const $ given by $A_\varphi
(f)=f\circ\Phi^{-1} $. Then a self-adjoint bounded operator
$$ P_\varphi^2=A_\varphi^* A_\varphi $$
commutes with $\beta\circ\pi_0$. If $P_\varphi^2=\int\lambda \ dE_\lambda $ is the 
spectral decomposition then
$E_\lambda $ commute with $\beta\circ\pi_0$.

\noindent\\
{\it Remarks.}\\
1) If $G$ is cocompact, then we know that $W_2^{1/2}/const\approx
H^1(G,l^2(G))$, so $W_2^{1/2}/const $ is a Hilbert module over the type II
factor defined by $G$ of dimension $\dim_G
W_2^{1/2}/const=L^2b_1(G)=2g-2$.\\
2) In practice, finding $A_\varphi$ is difficult. The reason is that
$\Phi$ is not a diffeomorphism, so the explicit formulae of Chapter 2
do not make sense. Moreover, $\Phi$ is given in a very implicit way as a
boundary value of a quasiconformal map, defined by a quadratic differential
on $\mathcal{H}^2$ which is $G$-invariant! 

\noindent\\
We will now show that for $p>2$ the operator $A_\varphi$ shows very
unusual properties, from the point of view of functional analysis.

\noindent\\
{\it Theorem 9.3.}--- Let $G$ be a fundamental group of a closed
hyperbolic surface $\Sigma^g$. Let $\varphi\in\mathcal{M}ap(\Sigma,x_0)$ be such 
that its image in $\mathcal{M}ap (\Sigma) $ is pseudo-Anosov. Let $A$ be the
operator, representing $\varphi $ in $W_p^{1/p}(S^1),\ p>2 $. Then there is
an element $0\not= v\in W_p^{1/p}(S^1) $ such that
$$ \sum_{k\in \mathbb{Z} }\| A_{\varphi}^k(v) \|^p <\infty .$$
{\it Proof.}--- Let $M$ be a mapping torus of $\Psi $, that is ,
$\mathbb{R}\times \Sigma /\mathbb{Z} $ where $1\in \mathbb{Z} $ acts by
$(t,x)\to (t+1,\Psi (x))$. Then $M$ is hyperbolic [Thurston 2]. We will
view $M$ as a fibration over a circle $\mathbb{R}/\mathbb{Z}$ with
coordinate $t,\ 0\le t<1 $; the fiber over $t$ will be called $\Sigma_t$. We
can trivialize $M\overset{\psi}{\to}\mathbb{R}/\mathbb{Z} $ over $I=[0,1/2]$
so that $(t,x_0),\ 0\le t\le 1/2 $ will be a horizontal curve. Let $g$ be the
hyperbolic metric on $M$ and $g_0$ be some hyperbolic metric on $\Sigma$, then
$g$ and $g_0+dt^2$ are equivalent on $\Sigma\times [0,1/2]\simeq
\psi^{-1}([0,1/2])\subset M$. Lifting to $\widetilde M=\mathcal{H}^3 $, we get a
fibration $\mathcal{H}^3\overset{\widetilde \psi}{\to}\mathbb{R} $ with
$\widetilde\psi^{-1}(t)={\widetilde\Sigma}_t $. Let $G\in PSL_2(\mathbb{C}) $ 
be the mondromy element,
corresponding to $\varphi$. Let $f:\mathcal{H}^3\to \mathbb{R} $ be such that
$$ \int_{\mathcal{H}^3}|\nabla f|^p\ d\mu_h <\infty .$$
We then have
$$ \sum_{k\in\mathbb{Z}} \int_{G^k({\widetilde \psi}^{-1} [0,1/2] )} |\nabla
f|^p\ d\mu_h \le \int_{\mathcal{H}^3} |\nabla f|^p\ d\mu_h <\infty ;$$
on the other hand the left hand side is
$$ \sum_{k\in \mathbb{Z}}\int_{{\widetilde\psi}^{-1}[0,1/2]} |\nabla (f\circ
G^k)|^p\ d\mu_h \ge const\cdot \int_{0}^{\frac{1}{2}}\ dt
\int_{{\widetilde\sum}_t} \sum_{k\in \mathbb{Z}} |\nabla(f\circ G^k) |^p \
dVol(g_0) .$$
It follows that for some $t_0$,
$$ \int_{{\widetilde\Sigma}_{t_0}} \sum_{k\in\mathbb{Z}} |\nabla(f\circ G^k)
|^p \ dVol(g_0) <\infty .$$
Since $g_0$ is a hyperbolic metric, for any function $F$ on $\widetilde\Sigma$
$$ \int_{\widetilde\Sigma} |\nabla F|^p\ dVol(g_0) =const\cdot \|
F|\partial\widetilde\Sigma \|_{W_p^{1/p}(S^1)/const}^p ,$$
actually, we may let the LHS be a definition of the norm in
$W_p^{1/p}(S^1)/const $, making the constant equal one. So
$$ \sum_{k\in\mathbb{Z}}\| f\circ G^k| \partial{\widetilde\Sigma}_{t_0}
\|_{W_p^{1/p}(S^1)/const }^p <\infty .$$
We will now identify $f\circ G^k |\partial{\widetilde\Sigma}_{t_0} $. We have
a boundary map
$$ \partial{\widetilde\Sigma}_{t_0} =S^1\overset{\alpha}{\to}
S^2=\partial\mathcal{H}^3 .$$
We know that $G^k\circ \alpha =\alpha\circ\varphi^{-k} $, so $f\circ
G^k=A_\varphi^k f $ and finally
$$ \sum_{k\in\mathbb{Z}}\|
A_\varphi^k(f|\partial{\widetilde\Sigma}_t)\|_{W_p^{1/p}(S^1)/const }^p
<\infty .$$

Now, for any $u\in W_p^{2/p}(S^2) $ we can take $f$ its harmonic extension.
In particular, any smooth function $u$ will do. Since $\alpha :S^1\to S^2 $
is continuous and nonconstant, we can take $u$ such that $V=u\circ \alpha $
is nonconstant. Then
$$ \sum_{k\in\mathbb{Z}}\| A_\varphi^kv \|_{W_p^{1/p}(S^1)/const}^p <\infty
,$$
as desired.

We remark that $\sum_{k\in\mathbb{Z}}\| A_\varphi^k v\|^p <\infty $ will
hold for all $v$ which are in the image of the bounded operator
$$ W_p^{2/p}(S^2)\to W_p^{1/p}(S^1) ,$$
induced by $\partial\widetilde\Sigma \to \partial \mathcal{H}^3 $.

\noindent\\
{\it Corollary 9.4.}--- Suppose that the space of fixed vectors of  $A_\varphi$
acting in $W_p^{1/p}/const$ possesses a complementary invariant subspace $W$. Then the spectre of $A_\varphi$ in
$W
$ satisfies
$$ \sigma(A_\varphi |W)\cap S^1 \not= \phi .$$
{\it Proof.}--- Suppose the opposite, then $W =W_+\oplus W_-$
such that $A_\varphi^k |W_+ $ and $A_\varphi^{-k}|W_- $ are strict
contractions for some $k>0$. But then
$$ \sum_{k\in\mathbb{Z}} \|A_\varphi^kv\|^p=\infty $$
for all $v\in W_p^{1/p}/const $.

We now turn to a generalization. Let $\widetilde G\subset \mathcal{M}ap_{g,1} $
be a subgroup, which containes $\pi_1(\Sigma_g) $, so that we have an extension
$$ 1\to \pi_1(\Sigma_g)\to \widetilde G\to G\to 1 .$$
Notice that $G\subset \mathcal{M}ap_g $. A well-known problem in hyperbolic
topology is : when there exists a fibration
\renewcommand{\baselinestretch}{1.5}
\large
$$\begin{array}{rcl}
	\Sigma & \to & Q              \\
        &   &  \downarrow       \\
        &   &  T
    \end{array}
$$
\renewcommand{\baselinestretch}{1.2}
\normalsize
with $\pi_1(Q)=\widetilde G $ such that $Q$ is a compact manifold of negative
curvature. In case $T$ is a closed surface, a corollary F.3 to the Theorem
F.1 of [Reznikov 9] provided some necessary condition. This condition is
unfortunately void, as we will show now.

\noindent\\
{\it Theorem 9.5.}--- Let $\Sigma^{g_1}\to Q\to \Sigma^{g_2} $ be a
surface fibration over a surface ($\Sigma^{g_i}$ are hyperbolic and oriented
). Let $\Sigma $ be a section of this fibration. Then
$$ |\Sigma\cap\Sigma |\le 2g_2-2 .$$
{\it Proof.}--- Let $\xi$ be a vertical tangent bundle for $\Sigma$, $e(\xi)$ its
Euler class, then $\Sigma\bigcap\Sigma =(e(\xi),[\Sigma] ) $. We have a natural
homomorphism $\pi_1 (Q)\to \mathcal{M}ap_{g_1,1} $ and a composite
homomorphism
$$ \pi_1 (\Sigma) \to \pi_1(Q)\to \mathcal{M}ap_{g,1} ,$$
which we call $\varphi$. An inclusion $\mathcal{M}ap_{g_1,1}\to\mathcal{G}_1
$ induces an Euler class $\epsilon$ in $H^2(\mathcal{M}ap_{g,1}) $ coming
from the action of $\mathcal{G}_1$ on $S^1$. By [Matsumoto-Morita 1], [Morita 2],  $
\varphi^{-1}\epsilon =e(\xi) $. Moreover, as is well known (and obvious )
$\epsilon $ is a bounded class, in fact, for any homomorphism $\pi_1
(\Sigma^g) \overset{\varphi}{\to} Homeo(S^1) , |(\varphi^*\epsilon, [\Sigma]
)|\le 2g-2 $. This proves the theorem.

\noindent\\
{\it Remarks}. \\
1) If the fibration $Q\to \Sigma^{g_2} $ is holomorphic and the action of
$\pi_1(\Sigma^{g_2}) $ on $H_1(\Sigma^{g_1},\mathbb{R})$ is simple, then a
famous inequality of Arakelov [Arakelov 1] reads $\Sigma\cap\Sigma<0$ for all
holomorphic sections. By Theorem 9.5,
$$ -(2g_2-2) <\Sigma\cap\Sigma <0 .$$
We now have a following result, which seems to be a very strong restriction
on $G$.

\noindent\\
{\it Theorem 9.6.}--- Let $1\to\pi_1(\Sigma^g)\to \widetilde G\to G $ be an
extension. Suppose $\widetilde G$ is a fundamental group of a compact manifold
of negative curvature
$$ -K \le K(Q^n )\le -1 .$$
Then for $p>(n-1)\sqrt{K} $ there is a vector $const\not= v\in
W_p^{1/p}(S^1) $, such that
$$ \sum_{g\in G} \|A_g v\|_{W_p^{1/p}/const}^p <\infty \eqno(*) .$$
{\it Proof.}--- Since the proof is essentially identical to the proof of
Theorem 9.3, we will only indicate the differences. Let $q_0 \in \widetilde Q
$ and let $u:S^{n-1}(T_{q_0}\widetilde Q)\to \mathbb{R} $ be a smooth function.
Composing with a geodesic projection $\widetilde Q\backslash \{0\} \to S^{n-1}
(T_{q_0}\widetilde Q) $ we arrive to a function $f: \widetilde Q\backslash
B(q_0,r)\to \mathbb{R} $ with $\int_{\widetilde Q} |\nabla f|^p \ dVol <\infty $
for $p>(n-1)\sqrt{K} $. Since $\Sigma$ is embedded in $Q$, one has a limit map
$\partial\widetilde\Sigma =S^1\to S^{n-1} =\widetilde Q $ be Theorem 6.1. Let
$v=u\circ \alpha $, where we identified $\partial\widetilde Q$ and
$S^{n-1}(T_{q_0}\widetilde Q ) $. Then $v\in W_p^{1/p}(S^1) $ by Theorem 6.2. As
in Theorem 9.3 we have the inequality $(*)$. Finally, if $v=const$, for any
choice of $u$, then $\alpha$ is almost everywhere a constant map, say to
$z\in S^{n-1} $. Since $\alpha $ is equivariant, it follows that
$\pi_1(\Sigma_g) $ stabilizes $z$. This is obviously impossible.

\section{Spaces $\mathcal{L}_{k,\alpha}^{(n-1)} $ and cohomology with
weights }

In this section we will describe a limit form of Theorem 4.1 when $p=1$, and
discuss $l^{n-1}$-cohomology with weights of cocompact lattices in $SO^{+}(1,n)$.

Let $G$ be a finitely generated group, $w:G\to R_+$ a function such that
$w(g)\to \infty$ as $length(g)\to \infty $. Consider a space $l^p(G,w)$ defined
by $f\in l^p(G,w)$ iff $\sum_g |f(g)|^pw^{-1}(g)<\infty $. Suppose $L_gw/w =O(1)
$ for all $g\in G$, and the same for $R_gw/w$. Then $l^p(G,w) $ becomes a
$G$-bimodule.

\noindent\\
{\it Example 1.}--- If $r(g)=length(g)$ then consider
$w(g)=r^{\alpha}(g),\ \alpha>0$ or $w(g)=r(g)^\alpha \log r(g)\log\log
r(g)\cdots\underbrace{\log\log\cdots\log}_k r(g), \ \alpha>0$. 

{\it 2.}--- Consider
$w(g)=e^{\alpha r(g)},\alpha>0 $.

\noindent\\
Now let $G$ be a cocompact lattice in $SO^{+}(1,n)$, We know by Theorem 5.1,
that $H^1(G,l^p(G))\not= 0$ exactly for $p>n-1$. In particular,
$H^1(G,l^{n-1}(G))=0$. However, by introducing of weights the
situation is changed.

\noindent\\
{\it Theorem 10.1.}--- Let $G$ be a cocompact lattice in $SO^{+}(1,n)$,
then for any $k\ge 1$ and $\alpha>0$,
$$ H^1(G,l^p(G,w))\not= 0 $$
for $p=n-1$ and $w=r(g)\log r(g)\cdots
(\underbrace{\log\log\cdots\log}_k r(g))^\alpha,\ \alpha>1,k\ge 1$.

\noindent\\
{\it Proof.}--- essentially repeats the argument of Proposition 2.1. Let
$u:S^{n-1}\to \mathbb{R} $ be any smooth function and denote again by $u$
its harmonic extension in $B^n$. We have $|\nabla u|_e<const$, therefore
$$ |\nabla u|_h(z)<const\cdot\rho_e(z,S^{n-1})^{-1} $$
Let $\mathcal{F}(h)=u(h^{-1}z_0)$, then a direct computation shows that
$L_g\mathcal{F}-\mathcal{F} \in l^{n-1}(G,w) $ and $\mathcal{F}-const \not=
l^{n-1}(G,w) $ so $l(g)=L_g\mathcal{F}-\mathcal{F} $ is a nontrivial cocycle if
$u$ is one of the coordinate functions on $S^{n-1}$, as in Theorem 2.1.

We would like to compute $H^1(G,l^{n-1}(G,w))$. A construction of Theorem 4.1
produces from any class in $H^1(G,l^p(G,w))$ a function in
$L_w^1(\mathcal{H}^n) $, where the latter space is defined as a space of
locally integrable function $f$ with distributional derivatives such that
$$ \int_{\mathcal{H}^n}|\nabla f|^{n-1} \cdot w^{-1}(z) <\infty \eqno(*) $$
where $w(z)=\rho_h (z_0,z)\log
\rho_h (z_0,z)\cdots(\underbrace{\log\log\cdots\log}_k\rho_h (z_0,z) )^{\alpha}$.

\noindent\\
{\it Definition.}--- A space $\mathcal{L}_{k,\alpha}^{(n-1)} $ is
defined as a Banach space of traces of $L_w^1(\mathcal{H}^n) $ on
$S^{n-1}$. A norm in $\mathcal{L}_{k,\alpha}^{(n-1)} $ is defined as infinum of
integrals $(*)$ taken over the set of all functions $f$ with a given trace.

\noindent\\
{\it Remark.}--- The norm just defined depends on $z_0$. Therefore a
natural action of $SO^{+}(1,n)$ in $\mathcal{L}_{k,\alpha}^{(n-1)} $ is not
isometric.

\noindent\\
We will describe $\mathcal{L}_{k,\alpha}^1 $ as a Zygmund-type space. One can
analogously describe $\mathcal{L}_{k,\alpha}^{(n-1)} $ for $n>2$, of course,
but we will not need it.

\noindent\\
{\it Theorem 10.2.}--- $\mathcal{L}_{k,\alpha}^1 $ consists of all
function $u:S^1\to\mathbb{R} $ for which $(a>0)$
$$ \int_0^a\ dh \int_0^{2\pi} \frac{|u(x+h)-u(x)|}{h^2\log
h\cdots{\underbrace{\log\cdots\log}_k}^\alpha h } <\infty .$$
{\it Proof} is a word-to-word repetition of Uspenski's argument in
[Uspenski 1]. One does not need to use Hardy's inequality, since $p=1$.

\noindent\\
{\it Theorem 10.3.}--- $\mathcal{G}_1 $ acts on $\mathcal{L}_{k,\alpha}^1 $ by
$$ A_\Phi u(x)=u\circ \Phi^{-1} .$$
{\it Corollary 10.4.}--- If $\Phi: S^1\to S^1 $ is quasisymmetric, then as a
function $S^1\to \mathbb{R}^2 $, $\phi \in \mathcal{L}_{k,\alpha}^1$.

\noindent\\
We suggest the reader to compare this result to [Carleson 2] and
[Gardiner-Sullivan 1]

\noindent\\
{\it Proof.}--- Let $\psi : B^2\to B^2 $ be a quasiisometry of the
hyperbolic metric, extending $\Phi$. If $u$ satisfies $(*)$ then $u\circ
\Phi^{-1} $ satisfies $(*)$ as well, whence the result.

\noindent\\
Embedding $\mathcal{M}ap_{g,1} \subset \mathcal{G}_1 $ we obtain a
representation
$$ \mathcal{M}ap_{g,1}\to Aut(\mathcal{L}_{k,\alpha}^1) ,$$
which is a limit case of Theorem 9.1.

\section{Bicohomology and the secondary quantization of the moduli space}

We will now introduce a very important notion of bicohomology spaces which
to an extent linearize 3-dimensional topology.

\noindent\\
{\it Definition.}--- Let $G$ be a finitely generated group. For $p>1$
define
$$ \mathcal{H}_p (G)=H^1 (G_r, H^1(G_l,l^p(G)) ,$$
where $r$ and $l$ stand for the right and left action, respectively.

\noindent\\
{\it Proposition 11.1.}--- A group $Out(G)$ of outer automorphism of
$G$ acts naturally in $\mathcal{H}_p(G) $.

\noindent\\
{\it Proof.}--- By definition, $Out(G)=Aut(G)/(G/Z(G))$. Obviously
$Aut(G) $ acts on $H^1(G_l,l^p(G))$ extending the right action of $G$, so
$Aut(G)/(G/Z(G)) $ will act on $H^1(G_r,(H^1(G_l,l^p(G))) $.

\noindent\\
For a surface group $\pi_1(\Sigma_g) $ we write
$\mathcal{H}_{p,g}=\mathcal{H}_p(G) $.

\noindent\\
{\it Theorem 11.2.}--- There exists a natural representation
$$ \mathcal{M}ap_g\to Aut (\mathcal{H}_{p,g} ). $$
Moreover, for $p>1$, $\mathcal{H}_{p,g}$ is a nontrivial Banach space. For
$p=2$, $\mathcal{H}_{p,g}$ is an infinite-dimensional Hilbert space. There
is a pairing
$$ \mathcal{H}_{p,g}\times \mathcal{H}_{p',g}\to \mathbb{R} \ ,$$
which is $\mathcal{M}ap_g$-invariant. For $p=p'=2$ this pairing is a
nondegenerate symmetric bilinear form. One has therefore a representation
$$ \mathcal{M}ap_g\overset{\psi}{\to} O(\infty,m),\ 0\le m\le \infty ,$$
which we call a secondary quantization of the moduli space of Riemann surfaces.

\noindent\\
The proof of the theorem will occupy the rest of this section.

\noindent\\
For a compact oriented manifold $M$ let $\Omega^{1/p}$ be a space of
measurable $1/p$-powers of densities such that for $\omega\in \Omega^{1/p}$
$$ \int_M |\omega|^p <\infty. $$
Then $\Omega^{1/p} $ is Banach, and for $p=2$, Hilbert. Let $G$ be a
finitely generated group acting in $M$.

\noindent\\
{\it Lemma 11.3.}--- Suppose that any element $g\in G$ has finitely
many repelling points, say $x_1^-,\cdots,x_n^-$ and finitely many attractive
points, say $x_1^+,\cdots,x_m^+$ such that for any set of neighbourhoods
$U_i^-,U_+^+ $ of $x_i^{\pm}$, there is $N$ such that for $k\ge N$,
$g^k(M\backslash \cup U_i^- )\subset \cup U_i^+ $. Suppose there are
$g_1,g_2,g_3,g_4\in G $ such that $\cup U_{i,s}^- \cup U_{i,s}^+ $ are
disjoint for different $s=1,2,3,4$. Then the action of $G$ in $\Omega^{1/p}
$ does not have almost-invariant unit vectors.

\noindent\\
{\it Proof.}--- Suppose the opposite, then there is a sequence
$\omega_j\in \Omega^{1/p} $, $\| \omega_j \|=1 $ and $\| g_s^k\omega_j
-\omega_j \|\underset{j\to\infty}{\to}0 $ for all $s,k$. Choose
$k_s,U_{i,s}^{\pm} $ such that
$$ g_s^{k_s}(M\backslash \cup U_{i,s}^- )\subset \cup U_{i,s}^+ $$
and $\cup U_{i,s}^- $(respectively $\cup U_{i,s}^+ $ ) don't intersect for
different $i$. Let $\omega$ be such that $\| \omega\|=1 $ and
$$ \|g_s^{k_s}(\omega)-\omega \|<(2/3)^{1/p} -(1/3)^{1/p}. $$
For $E\subset M$ define $C(E,\omega)=\int_E |\omega|^p $. We claim that
$$ C(M\backslash \cup U_{s,i}^- \backslash \cup U_{s,i}^+, \omega )<2/3 .$$
Suppose the opposite, then by the invariance of the density $|\omega|^p$,
\renewcommand{\baselinestretch}{1.5}
\large
$$\begin{array}{c}
    C(M\backslash \cup U_{s,i}^- \backslash \cup U_{s,i}^+, \omega\circ
		g_s^{k_s} )\le \\
    =C(g_s^{k_s} (M\backslash \cup U_{s,i}^- \backslash \cup U_{s,i}^+
		),\omega )\le \\
    \le C(g_s^{k_s}(M\backslash \cup U_{s,i}^- ),\omega ) \le 1/3.
    \end{array}
$$
\renewcommand{\baselinestretch}{1.2}
\normalsize
It follows that
\renewcommand{\baselinestretch}{1.5}
\large
$$\begin{array}{c}
    [\int_{M\backslash \cup U_{s,i}^-\backslash \cup U_{s,i}^+ }
		|\omega-\omega\circ g_s^{k_s} |^p]^{1/p} \ge \\
    \ge |[\int_{M\backslash \cup U_{s,i}^- \backslash U_{s,i}^+}| \omega^p |]^{1/p}-
        [\int_{M\backslash \cup U_{s,i}^-\backslash \cup U_{s,i}^+ }
		|\omega\circ g_s^{k_s}|^p] ^{1/p} |\ge (2/3)^{1/p} -(1/3)^{1/p},
    \end{array}
$$
\renewcommand{\baselinestretch}{1.2}
\normalsize
a contradiction.

So $C(\cup U_{s,i}^-,\omega)+C(\cup U_{s,i}^+,\omega)\ge 1/3 $. Since
$ \cup U_{s,i}^{\pm} $ are disjoint for different $s$, we get
$$ 1\ge \sum_{s=1}^4 C(\cup U_{s,i}^-,\omega )+C(\cup U_{s,i}^+,\omega)\ge
	4/3 ,$$
a contradiction. This proves the lemma.

\noindent\\
{\it Corollary 11.4.}--- Let $G\subset SO^{+}(1,n) $ be a cocompact
lattice. Then the natural isometric action of $G$ in $W_p^{(n-1)/p}(S^{n-1})
$ does not have almost-invariant vectors. In particular,
$H^1(G,W_p^{(n-1)/p}(S^{n-1})) $ is Banach for $p>(n-1) $.

\noindent\\
{\it Proof.}--- For $u\in W_p^{(n-1)/p}(S^{n-1})/const $ , let $f$ be
a harmonic extension of $u$ so that
$$ \| u\| =\int_{\mathcal{H}^n}|\nabla f|^p .$$
Since the energy density $|\nabla f|^p d\mu_h $ is invariant under
isometries of $\mathcal{H}^n$, the proof of the Lemma 11.3 applies directly.

\noindent\\
{\it Corollary 11.5.}--- $\mathcal{H}_{p,g} $ is Banach (Hilbert for
$p=2$ ).

\noindent\\
{\it Proof.}--- $H^1(G_l, l^p(G) )=W_p^{1/p} (S^1)/const $. \\

We now describe the pairing
$$ \mathcal{H}_{p,g}\times \mathcal{H}_{p',g} \to \mathbb{R} \ .$$
This is given by the cup-product in cohomology
$$ H^1(G_r,H^1(G,l^p(G)))\times H^1(G_r,H^1(G_l,l^{p'}(G)))\to
	H^2(G_r,H^2(G_l,l^p(G)\otimes l^{p'}(G))) $$
followed by the duality $l^p(G)\times l^{p'}(G)\to \mathbb{R} $ and 
evaluating twice on the fundamental cycle in $H_2(G,\mathbb{R}) $. 
We have also an analytic description, namely a pairing
$$ W_p^{1/p}(S^1)/const \times W_{p'}^{1/p'}(S^1)/const\to 
	\mathbb{R} $$
is given on smooth function by $f,g\to \int_{S^1}fdg $ and then
extended as a bounded bilinear form. This induces a pairing
$$ H^1(G,W_p^{1/p}/const)\times H^1(G,W_{p'}^{1/p'}/const)\to 
	\mathbb{R} \ .$$
{\it Lemma 11.6.}[Korevaar-Schoen 1]--- Let $G$ be a finitely
presented group which is realized as a fundamental group of a
compact Riemannian manifold $M$. Let $\rho : G\to O(\mathcal{H})$
be an orthogonal representation, which does not have almost-invariant
vectors. Let $[l]\in H^1(G,\mathcal{H})$. Let $E$ be a flat vector
bundle with fiber $\mathcal{H}$ over $M$, corresponding to $\rho$. 
Then there is a harmonic 1-form $\omega\in \Omega^1(M,E) $, 
corresponding to $[l]$.

\noindent\\
{\it Proof.}---This is a reformulation of [Korevaar-Schoen 1].

\noindent\\
{\it Corollary 11.7.}--- Let $M$ be K\"ahler. Then if $\rho$ is as 
in the previous lemma, then \\
1) There is a natural complex structure in $H^1(G,\mathcal{H})$,
making it a complex Hilbert space ;\\
2) A pairing 
$$ H^1(G,\mathcal{H})\times H^1(G,\mathcal{H} ) \to \mathbb{R} \ ,$$
given by $[l_1],[l_2]\to ([\omega]^{n-1}([l_1],[l_2]), [M] ) $
where $[\omega] $ is a K\"ahler class, $[M]$ is the fundamental
class and $([l_1],[l_2])\in H^2(G,\mathbb{R}) $ is a cup-product
composed with the scalar product $\mathcal{H}\times\mathcal{H}
\to \mathbb{R} $, is a non-degenerate symplectic structure in
$H^1(G,\mathcal{H}) $.

\noindent\\
{\it Proof } is the same as for finite-dimensional $\mathcal{H}$,
once we have the Hodge theory of the previous lemma.
\\

We now ready to prove that the symmetric pairing
$$ \mathcal{H}_{2,g}\times \mathcal{H}_{2,g}\to \mathbb{R} $$
is nondegenerate. Realize $G$ as a lattice in $SO(1,2)$. Then
$H^1(G,l^p(G))=W_p^{1/p}(S^1)/const $. Let $H$ denote
the Hilbert transform. It is a bounded operator
$$ H: L^p(S^1)/const\to L^p(S^1)/const \quad (p>1) $$
defined as follows: for $u\in L^p(S^1) $ let $f$ be
its harmonic extension and $g$ a conjugate harmonic function,
then $Hu=g|S^1 $. Since
$$ \int_{\mathcal{H}^2}|\nabla f|^p =\int_{\mathcal{H}^2} 
	|\nabla g|^p .$$
$H$ restricts to $W_p^{1/p}(S^1) $ as an isometry.

Now, the symplectic duality $\int f\ dg $ in $W_2^{1/2}(S^1)/const $ is
simply equal to $(Hf,g)$. Moreover, $H$ is
$SO(1,2)$-invariant. Then the Corollary 11.7 implies that the pairing of
Theorem 11.2 is also nondegenerate.

We still have to prove that $\mathcal{H}_{g,p}\not= 0$ and for $p=2$ is
infinite-dimensional. We first describe an element of $\mathcal{H}_{g,p}$
associated to a given realization $G\hookrightarrow SO(1,2)$ as a cocompact
lattice,which we will call a principal state.

Recall that if $M$ is a smooth compact oriented manifold,
$\mathcal{D}iff^1(M)$ a group of orientation-preserving diffeomorphism of
class $C^1$, then one has a cocycle $l$ in $Z^1(\mathcal{D}iff^1(M),C^0(M))$
defined as [Bott 1]
$$ l(\phi)=\log \frac{\phi_*\mu}{\mu} ,$$
where $\mu$ is any smooth density on $M$, and $\phi_*\mu$ a left action. The
class $[l]\in H^1(\mathcal{D}iff^1(M),C^0(M))$ does not depend on $\mu$. For
$r\ge 1$ one similarly gets a class in
$H^1(\mathcal{D}iff^r(M),C^{r-1}(M))$. Now, let $M=S^{n-1}$ and consider a
standard conformal action of $SO^{+}(1,n)$ on $S^{n-1}$. We get a class
$$ [l]_p\in H^1(SO^{+}(1,n),W_p^{(n-1)/p}(S^{n-1})/const) $$
for all $p>1$ simply because $C^\infty (S^{n-1})\subset W_p^{(n-1)/p}
(S^{n-1}) $. We claim $[l]_p\not= 0$ for $n=2$ and $p>n-1$. Since the action
is isometric, it follows from the following lemma (we prove and use it only
for $n=2$).

\noindent\\
{\it Lemma 11.8.}--- Fix $z_0\in B^n $ and let
$r(g)=\rho_h(z_0,g^{-1}z_0)$. Then for any fixed $\mu$, $\| l(g)
\|_{W_p^{(n-1)/p}(S^{n-1})/const } \to \infty $ as $r(g)\to \infty $.

\noindent\\
{\it Proof.}(Only for $n=2$)--- We choose for $\mu$ the harmonic
(Poisson) measure $\mu_0$, associated with $z_0$. Then $l(g)=\log
\frac{g_*\mu_0}{\mu_0} $. For $\beta\in S^{n-1}$, $l(g)(\beta)=B_\beta
(z_0,gz_0)$ where $B_\beta (z_0,\cdot)$ is a Busemann function of $B^n$
corresponding to $\beta\in \partial B^n$ and normalized at $z_0$, that is,
$B_\beta(z_0,z_0)=0$ (see, for example, [Besson-Courtois-Gallot 1]).

We will make the computation only for $n=2$. Let $z_0=0$, $gz_0=w$,
then
$$ B_\beta(0,w)=\log \frac{1-|w|^2}{|w-\beta|^2} .$$
Notice that $\log |\frac{\beta-w}{1-\bar w\beta }| =0 $, since $|\beta|=1$,
so
$$ B_\beta(0,w)=\log(1-|w|^2)-2\log|1-\bar w\beta |=-2\log|1-\bar
w\beta|(\mbox{\it{mod const}}) .$$
Notice that $\log |1-\bar w z| $ is defined and is harmonic in $|z|\le 1$,
so
\renewcommand{\baselinestretch}{1.5}
\large
$$\begin{array}{c}
    \|B_\beta(0,w)\|_{W_p^{1/p}(S^1)/const}^p=2^p\int_{B^2} [\nabla (\log
		|1-\bar w z) |)]_h^p\ d\mu_h =\\
    =2^p \int_{B^2} \frac{|w|^p}{|1-\bar w z|^p} \frac{1}{(1-|z|^2)^{2-p} }\
		dzd\bar z
    \end{array}     \eqno(*)
.$$
\renewcommand{\baselinestretch}{1.2}
\normalsize
{\it Sublemma.}--- An integral $(*)$ grows as $\log (1-|w|) $ as $|w|\to 1$.

\noindent\\
{\it Proof.}--- Computing in polar coordinates, we have
$$ \int_0^1 dr\frac{1}{(1-r^2)^{2-p} } \int_0^{2\pi}
\frac{d\theta}{|1-r|w|e^{i\theta}|^p} .$$
It is elementary to see that the inner integral grows as
$\frac{1}{(1-r|w|)^{p-1}}$, so we arrive at
$$ \int_0^1 dr \frac{1}{(1-r)^{2-p} }\frac{1}{(1-r|w|)^{p-1}}\sim
\int_0^a\frac{ds}{s^{2-p}(A+s)^{p-1}} $$
where $a>0$ is fixed and $A=1-|w|$. Further we have ($s=At$)
$$ \int_0^{a/A}\frac{dt}{t^{2-p}(1+t)^{p-1}}\sim \int_0^{a/A}
	\frac{dt}{t}\sim \log |A| ,$$
which proves the Sublemma.

Finally,
$$ \| B_\beta (0,w)\|_{W_p^{1/p}(S^1)/const} \sim [\log (1-|w|)]^{1/p}, $$
where $\sim$ means that the ratio converges to a constant.

The proof for $n>2$ will be given elsewhere.

Notice that for $p=2$ we have (for $n=2$)
$$ \| l(g) \|_{W_2^{1/2}(S^1)/const }\sim \| g\|^{1/2} $$
where $\|g\|$ is a hyperbolic length of a (pointed) geodesic loop,
representing $g$. This exponent in the RHS is the maximal possible. We will
later prove a general theorem (Theorem III.3.1) showing that for any
orthogonal or unitary representation of $G=\pi_1(\Sigma)$ in a Hilbert space
$\mathcal{H}$ and any cocycle $l\in Z^1(G,\mathcal{H})$,
$$ \| l(g) \|\le const\cdot\mathit{length}(g)^{1/2}\log\log\mathit{length}(g)
$$
as $g$ converges nontangentially to almost all $\theta\in S^1=\partial G $.

Coming back to principal states $[l]_p\in
H^1(G,W_p^{(n-1)/p}(S^{n-1})/const)$, let $E$ be a flat affine bundle over
$M=\mathcal{H}^n/G$ with fiber $W_p^{(n-1)/p}(S^{n-1})$, associated with an
affine action
$$ g\mapsto R_g+l(g) .$$
Notice that
$$ s: z\mapsto log \frac{\mu (z)}{\mu(z_0)} $$
is an $G$-equivariant section of the lift of $E$ on $\widetilde
M=\mathcal{H}^n$, or, equivalently, defines a section of $E$. We claim that
this section  is harmonic. This immediately reduces to a statement
that $B_\beta(z_0,z)$ is harmonic $\mbox {mod const}$ as a function of $z$. In the upper
half-plane model it simply means that $(x,y)\mapsto \log y$ is harmonic $\mbox {mod const}$. The harmonic section just defined does not lift to a harmonic section of the flat affine bundle with fiber $W_p^{(n-1)/p}(S^{n-1})$. 
For $n=2$ we can say more. Let
$$ H:W_p^{1/p}(S^1)/const\to W_p^{1/p}(S^1)/const $$
be the Hilbert transform. It makes $W_p^{1/p}(S^1)/const $ into a complex
Banach space. Then a direct inspection shows that the section of $E$ defined
above is (anti)holomorphic (depending on the choice of a sign of $H$). This
will be used later. Equivalently, $ds$ is an (anti)-holomorphic one-form on
$\mathcal{H}^2/G$, valued in $E$. Again, this holomorphic form does not lift to a $d$ and $\delta$-closed form of a flat bundle with fiber $W_p^{1/p}(S^1)$ even for $p=2$. This latter bundle is a flat bundle with fiber a Hilbert space, but whose monodromy is not orthogonal. The Hodge theory of [Korevaar-Schoen 1] and [Jost 1] does not apply and in fact not every cohomology class is represented by a  $d$ and $\delta$-closed form. We will discuss these subtle obstructions to the Hodge theory in a sequel to this paper [

Reznikov 10].  

As an application of the computation made above, we will complete the proof
of Lemma 5.6 for $p>1$. Let $u\in W_p^{1/p}(S^1)/const $ and let $f:B^2\to
\mathbb{R} $ be a harmonic extension of $u$. We claim that
$$ |f(w)|\le c\cdot [\log (1-|w|)]^{1/p'} .$$
Since the Hilbert transform is invertible in $W_p^{1/p}(S^1)/const$, we can
assume that the Fourier coefficients $\hat u(n)=0 $ for $n<0$, so that $f$
is holomorphic:
\renewcommand{\baselinestretch}{1.5}
\large
$$ \begin{array}{rl}
    &|f(w)|=|\frac{1}{2\pi i} \int_{S^1} \frac{u(\xi)d\xi
		}{\xi-w}|=\frac{1}{2\pi}|\int_0^{2\pi}\frac{u(e^{i\theta})
        e^{i\theta}}{e^{i\theta}-w}d\theta |= \\
    =&\frac{1}{2\pi}|\int_0^{2\pi}\frac{u(e^{i\theta})d\theta}{1-w\cdot
		e^{-i\theta}}|=\frac{1}{2\pi}|\int_0^{2\pi}\frac{u(e^{-i\theta})
		d\theta}{1-w\cdot e^{i\theta} }|= \\
	=&|\frac{1}{2\pi i}\int_0^{2\pi} \frac{[u(e^{-i\theta})\cdot e^{-i\theta}]
		ie^{i\theta}d\theta }{1-w e^{i\theta} }| =|-\frac{1}{2\pi iw}
		\int_0^{2\pi } [u(e^{-i\theta})\cdot e^{-i\theta} ]
		[\log (1-we^{i\theta})]'d\theta |= \\
	=&|-\frac{1}{2\pi iw} < u(e^{-i\theta})\cdot e^{-i\theta} ,
		\log(1-w e^{i\theta}) >|\le \\
	\le & \frac{1}{2\pi |w|}\| u(e^{-i\theta})e^{-i\theta} \|_{W_p^{1/p}(S^1)/const}
		\cdot \| \log(1-w e^{i\theta})\|_{W_{p'}^{1/p'}/const }\le \\
	\le & c\| u\|_{W_p^{1/p}(S^1)/const } |\log (1-|w|)|^{1/p'}.
 \end{array}
$$
\renewcommand{\baselinestretch}{1.2}
\normalsize
It is very plausible that the result is true, for $u\in W_p^{\frac{n-1}{p}}
(S^{n-1})/const $ for $n\ge 3$. Our proof obviously does not work in this case.

We now start to prove that $\mathcal{H}_{2,g}$ is infinite-dimensional. Let
$M_0,M_0'$ be factors generated by the left (respectively, right) actions of
$G$ in $l^2(G) $ [Murray-von Neumann 1]. Notice that $H^1(G_l,l^2(G))$ can
be viewed as a cohomology of a complex
$$
l^2(G)\overset{d_0}{\to}\bigoplus_{i=1}^{2g}l^2(G)\overset{d_1}{\to}l^2(G)
\eqno(*) ,$$
computed from the standard CW-decomposition of $\Sigma^g$ with one
zero-dimensional cell, $2g$ one-dimensional cells and one two-dimensional
cell. Notice that $d_0,d_1$ are given by matrices with entries in
$\mathbb{Z}[G]$, acting on $l^2(G)$ from the left. Letting
$\Delta_l=d_0d_0^*+d_1^*d_1$ we can view $H^1(G,l^2(G))$ as $Ker\Delta_l $.
Notice that $\Delta_l \in M_0 $. It follows that $H^1(G,l^2(G))$ is a module
over $M_0'$.
Now, since $M_0$ is type II, there is a decomposition
$$ W_2^{1/2}(S^1)/const =H^1(G_l,l^2(G))=Ker\ \Delta_l =\bigoplus_{j=1}^m
H_m ,$$
for any $m\ge 1$ where $H_m$ are isomorphic right $G$-modules. Since we know
already that $H^1(G_r,W_2^{1/2}(S^1)/const)\not= 0$, and $H_j$ are
all isomorphic, it follows that $H(G_r,H_j)\not= 0 $ for all $j$,
therefore $\dim H^1(G,W_2^{1/2}(S^1))/const \ge m$. This finally proves
Theorem 11.2.
\noindent\\
 There are natural invariant von Neumann algebras acting in $\mathcal{H}_{2,g}$. Indeed, let $M_1'$ be a double commutant of $M_0'$ in  $H^1(G,l^2(G))=Ker\Delta_l$  and $M_1$ be
a commutant of $M_0'$. We could
define $M_1'$ as a von Neumann algebra , generated by the right action of $G$ in
$H^1(G,l^2(G)) $ and $M_1$ as a commutant of $M_1'$. It follows that
$M_1,M_1'$ do not depend on the choice of the complex $(*)$ and therefore
$\mathcal{M}ap_{g,1}=Aut(G)$ acts in $H^1(G,l^2(G))$ leaving $M_1',M_1$
invariant. Now consider $\mathcal{H}_{2,g}=H^1(G_r,H^1(G_l,l^2(G))$. Then
$\mathcal{H}_{2,g}=Ker\ \Delta_r: 
\bigoplus_{i=1}^{2g}H^1(G_l,l^2(G))\to 
\bigoplus_{i=1}^{2g}H^1(G_l,l^2(G))$ where
a right Laplacian is defined exactly as the left one. It follows that
$\mathcal{H}_{2,g}$ is a module over $M_1$. Let $M_2$ be a double commutant
of $M_1$ and $M_2'$ be its commutant. We have proved a following theorem,
except for the last statement.

\noindent\\
{\it Theorem 11.9.}--- There are infinite-dimensional von Neumann algebras $M_2, M_2'$ acting in
$\mathcal{H}_{2,g}$, which are invariant under the action of
$\mathcal{M}ap_g$. Moreover, there is an involution $\tau$ of
$\mathcal{H}_{2,g}$ which commutes with the $\mathcal{M}ap_g$-action and
permutes $M_2,M_2'$.

\noindent\\
{\it Proof.}--- Everything is already proved except for the last
statement. Notice that there is an involution $\tau:l^2(G)\to l^2(G) $
defined by $\tau f(g)=f(g^{-1})$. A Lyndon-Serre-Hochschild spectral
sequence of the extension $1\to G\to G\times G\to G\to 1 $ shows that
$\mathcal{H}_{2,g}=H^2(G\times G,l^2(G))$. Let $\sigma$ be an involution of
$G\times G$ defined by $\sigma (g,h)=(h,g)$. Then one has $\tau
[(g,h)v]=(\sigma (g,h))\tau (v) $ where $g,h\in G$ and $v\in l^2(G)$. It
follows that $\tau$ induces an involution, which we also call $\tau$, in
$\mathcal{H}_{2,g}$, which obviously commutes with $\mathcal{M}ap_g$-action
and permutes $M_2$ and $M_2'$. This completes the
proof of Theorem 11.9.

Note that since the unitary representation of $G$ in $H^1 (G_l , l^2 (G))=W_2^{1/2}(S^1)/const$ extends to an irreducible representation of $PSL_2(\mathbb{R} )$, the commutator $M_1$ of $G$ in $W_2^{1/2}(S^1)/const $  possesses a faithful trace defined by 
$$tr(a)\cdot Id=\int_{PSL_2(\mathbb{R})/G} gag^{-1}dg.$$

\noindent\\
{\it Proposition 11.10.}--- Let
$\widetilde{\mathcal{H}}_{2,g}$ be a completion of
$M_1$ under the norm $tr\ xx^*$. Then
$\widetilde{\mathcal{H}}_{2,g}$ is a  Hilbert space
and  there is  a  representation
$$ \Tilde\rho: \mathcal{M}ap_g\to 
	Aut(\widetilde{\mathcal{H}}_{2,g}),$$
leaving invariant a nondegenerate form $x\mapsto tr\ x^2 $.

\noindent\\
I don't know at the time of writing if 
$\widetilde{\mathcal{H}}_{2,g}$ is isomorphic to
$\mathcal{H}_{2,g}$ as $\mathcal{M}ap_g$-module. 

We now turn to the holomorphic realization of $\mathcal{H}_{2,g}$. Fix a
realization of $G$ as a cocompact lattice in $SO^{+}(1,2)$, then
$\mathcal{H}_{2,g}=H^1(G,W_2^{1/2}(S^1)/const )$. Recall that $G$ commutes
with the Hilbert transform in $W_2^{1/2}(S^1)/const $. Let
$S=\mathcal{H}^2/G$, then $S$ is a hyperbolic Riemann surface, homeomorphic
to $\Sigma^g$. For any element $w\in H^1(G,W_2^{1/2}(S^1)/const) $ we have by
Lemma 11.3 and Lemma 11.6 a unique harmonic form in a flat Hilbert bundle
$E$ with fiber $W_2^{1/2}(S^1)/const $, associated with the action of $G$.

Uniqueness should be explained. We have a following general fact.

\noindent\\
{\it Lemma 11.11.}--- Let $M$ be a compact Riemannian manifold.
$\rho:\pi_1(M)\to O(H)$ an orthogonal representation in a real Hilbert
space, without fixed vectors, $\omega\in H^1(\pi_1(M),H)$. Then there at
most one harmonic form, $\omega\in\Omega^1(M,E)$, representing $\omega$.

\noindent\\
{\it Proof.}--- If $\omega_1,\omega_2$ are two such forms, then
$\omega_1-\omega_2$ is a derivative of a harmonic section of $M$. But
standard Bochner vanishing theorem shows that such section should be
self-parallel, so $\rho$ has a fixed vector, a contradiction.

Notice that $H$ makes $W_2^{1/2}(S^1)/const $ into a complex Hilbert space.
Then $\frac{1}{2}(\omega-H(\omega\circ J))$, where $J$ is a complex
structure on
$S$, will be a holomorphic 1-form in $E$, whereas
$\frac{1}{2}(\omega+H(\omega\circ J)) $ will be an anti-holomorphic 1-form.
Let $\mathcal{H}_{2,g}^{\pm} $ be the spaces of holomorphic(respectively,
anti-holomorphic) 1-forms in $E$, then
$\mathcal{H}_{2,g}=\mathcal{H}_{2,g}^+\bigoplus \mathcal{H}_{2,g}^- $. Now,
$W_2^{1/2}(S^1)/const $ is identified with exact $L^2$-harmonic 1-forms in
the hyperbolic plane 
$\mathcal{H}^2$, which is isomorphic as a complex Hilbert space ( with
a complex structure, defined by 
the Hodge star operator) to the space of exact
$L^2$-holomorphic 1-form in $\mathcal{H}^2$. So any element in
$\mathcal{H}_{2,g}^+$ defines a holomorphic 1-form on $S$ valued in a bundle
with fibers $L^2$-holomorphic 1-forms on $\mathcal{H}^2$. In other words, let
$G$ act diagonally in $\mathcal{H}^2\times \mathcal{H}^2 $ and
$$ Q=\mathcal{H}^2\times\mathcal{H}^2 /G ,$$
then we have an $L^2$ holomorphic 2-form on $Q$. The space
$\mathcal{H}_{2,g}^+$ therefore is identified with the space of $L^2$
holomorphic 2-forms on $Q$. Similarly, $\mathcal{H}_{2,g}^-$ is identified
with the space of $L^2$ holomorphic 1-form on
$$ Q'=\mathcal{H}^2\times \overline{\mathcal{H}^2}/G ,$$
where $\overline{\mathcal{H}^2}$ is obtained from $\mathcal{H}^2$ by reversing
the complex structure (i.e. $\bar J=-J$). Notice that as complex surfaces,
$Q$ and $Q'$ are not biholomorphic: $Q$ contains a compact curve (the
quotient of the diagonal) whereas $Q'$ does not. We have proved the
following:

\noindent\\
{\it Theorem 11.12}.(Holomorphic realization of quantum moduli space)--- Fix an embedding $G\hookrightarrow
SO^{+}(1,2)$ as a cocompact surface, then $\mathcal{H}_{2,g} $ splits as
$\mathcal{H}_{2,g}^+\oplus \mathcal{H}_{2,g}^- $ where
$\mathcal{H}_{2,g}^+$ (respectively, $\mathcal{H}_{2,g}^-$) is identified
with a space of $L^2$ holomorphic 2-forms on $Q=\mathcal{H}^2\times
\mathcal{H}^2 /G $ (respectively, $Q'=\mathcal{H}^2\times
\overline{\mathcal{H}^2}/G$). Moreover, the splitting is orthogonal with respect
to the canonical symmetric scalar product in $\mathcal{H}_{2,g}$ and the
restriction of this scalar product on $\mathcal{H}_{2,g}^{\pm} $ is
positive (respectively, negative).

\noindent\\
{\it Example.}--- The principal state $[l]_2$ lies in
$\mathcal{H}_{2,g}^-$. We do not know at the time of writing if $\mathcal{H}_{2,g}^+ =0.$

\section{$\mathcal{H}_{p,g}$ as operator spaces and the vacuum vector}

In this section we will develop an algebraic and an analytic theory of
$\mathcal{H}_{p,g}$ as spaces of operators between $W_q^{1/q}(S^1)/const$,
which commute with the action of $G$. We use rather rough estimates of
matrix elements, so the ranges of indices for which the action is
established is certainly not the best possible. We start with a lemma.

\noindent\\
{\it Lemma 12.1.}--- Let $u\in W_p^{1/p}(S^1)$ and $a\in l^q(G)$. Then
$\sum a(g) R_g u\in W_r^{1/r}(S^1) $ where
$\frac{1}{p}+\frac{1}{q}-1=\frac{1}{r} $.

\noindent\\
{\it Remark.}--- $R_g$ means the action of $G$ in $W_p^{1/p}(S^1)$ (a
reminiscent of the actions from the right in $l^2(G) $).

\noindent\\
{\it Proof.}--- Let $f$ be a harmonic extension of $u$, so that
$$ \int_{\mathcal{H}^2}|\nabla f|^p\ d\mu_h <\infty .$$
By Young-A.Weil inequality [Hewitt-Ross 1], $l^p*l^q\subset l^r $, so if $h=\sum
a_gR_g f $, we have
\renewcommand{\baselinestretch}{1.5}
\large
$$\begin{array}{rcl}
	 \int_{\mathcal{H}^2} |\nabla h|^r\ d\mu_h &=&
\int_{\mathcal{H}^2/G}d\mu_h(z) \sum_g |\nabla h(g^{-1}z)|^r \\
    &\le& \int_{\mathcal{H}^2/G}d\mu_h(z) \| \nabla h(g^{-1}z)\|_{l^r}\\
    &\le& c\cdot \int_{\mathcal{H}^2/G}d\mu_h(z)\|a(g)\|_{l^q(G)} \|\nabla
f(g^{-1}z)\|_{l^p} \\
    &\le& c\cdot \|a(g)\|_{l^q(G)} \int_{\mathcal{H}^2} |\nabla f|^p d\mu_h.
    \end{array}
 $$
\renewcommand{\baselinestretch}{1.2}
\normalsize
 The result follows with an estimate
 $$ \| \sum a(g) R_gu\|_{W_r^{1/r}(S^1)/const} \le c\cdot
\|a(g)\|_{l^q(G)}\|u\|_{W_p^{1/p}(S^1)/const} .$$
Now recall that we have a canonical pairing
$$ B: W_r^{1/r}(S^1)/const \times W_{r'}^{1/r'}(S^1)/const\to \mathbb{R} ,$$
so that a formula
$$ (u,v)\mapsto (a\mapsto B(\sum a(g)R_g u,v)) $$
defines a map
$$ W_p^{1/p}(S^1)/const \times W_{r'}^{1/r'}/const\to l^{q'}(G) .$$ Now an element of $W_{r'}^{1/r'}(S^1)/const $
defines an element of
$$ Hom_G(W_p^{1/p}(S^1)/const, l^{q'}(G)) $$
and an induced map
\renewcommand{\baselinestretch}{1.5}
\large
$$ \begin{array}{c}
    Hom(H^1(G,W_p^{1/p}(S^1)/const, H^1(G,l^{q'}(G))=\\
    =Hom(H^1(G,W_p^{1/p}(S^1)/const,W_{q'}^{1/q'}(S^1)/const).
    \end{array}
$$
\renewcommand{\baselinestretch}{1.2}
\normalsize
In other words, we have a map
$$ H^1(G,W_p^{1/p}(S^1)/const)\to Hom_{\mathbb{R}}(W_{r'}^{1/r'}(S^1)/const
\to W_{q'}^{1/q'}(S^1)/const ) ,$$
and it is immediate to check that the image lies in $Hom_G$. So we have a 

\noindent\\
{\it Proposition 12.2.}--- The construction above defines a map
$$ \mathcal{H}_{p,g}\to Hom_G(W_{r'}^{1/r'}(S^1)/const ,
W_{q'}^{1/q'}(S^1)/const ) $$
for $p,q',r' $ satisfying $\frac{1}{p}\ge 1+\frac{1}{q'}-\frac{1}{r'} $,
which is $\mathcal{M}ap_g$-equivariant.

\noindent\\
An induced map in $H^1(G,\cdot)$ produces a bounded
$\mathcal{M}ap_g$-equivariant product
$$ \mathcal{H}_{p,g}\times \mathcal{H}_{r',g}\to \mathcal{H}_{q',g} .$$
We stress again that the range of indices for which 
this product is defined should be
improved. We will see that viewing $\mathcal{H}_{p,g}$ as an operator space
helps to understand $\mathcal{M}ap_g$-action. We turn now to an analytic
description of the above. Let $l\in Z^1(G,W_p^{1/p}(S^1)/const) $. A
construction of Theorem 4.1 produces a smooth map
$$ F: \mathcal{H}^2\to W_p^{1/p}(S^1)/const ,$$
satisfying $F(g^{-1}z)=R_g F(z)+l(g), g\in G$, In particular, $g^*(\nabla
F)(g^{-1}z)=R_g(\nabla F)(z) $. Now let $v\in W_{r'}^{1/r'}(S^1)/const $,
where $r\ge p$. Then we have a scalar function
$$ <F,v>: \mathcal{H}^2\to \mathbb{R} ,$$
where $<\cdot,\cdot>$ is a pairing $W_r^{1/r}(S^1)/const \times
W_{r'}^{1/r'}(S^1)/const \to \mathbb{R} $ defined in Theorem 7.2. Since
$r\ge p$, $W_p^{1/p}(S^1)\subset W_r^{1/r}(S^1)$, so $<F,v>$ is defined.
Without futher assumption one can only say that
$$ |\nabla <F,v>|\le const ,$$
but if we assume $r>p$, say $\frac{1}{p}=1+\frac{1}{q'}-\frac{1}{r'}$, then
$<F,v>$ will satisfy
$$ \sum_g |\nabla <F,v> (gz)|^{q'} <const $$
for all $z\in \mathcal{H}^2$. Integrating over $\mathcal{H}^2/G$, we get
$$ \int_{\mathcal{H}^2}|\nabla <F,v> |^{q'} <\infty ,$$
so there exists $<F,v>|S^1\in W_{q'}^{1/q'}(S^1)$. 
This defines a desired map
$$ H^1(G,W_p^{1/p}(S^1)/const)\to Hom_G (W_{r'}^{1/r'}(S^1)/const,
W_{q'}^{1/q'}(S^1)/const ) .$$
We will use this description now to compute the operator, associated with
the principal state
$$ [l]_p\in H^1(G, W_p^{1/p}(S^1)/const) .$$
{\it Proposition 12.3.}--- For $p,r',q'>1,\ \frac{1}{p}\ge
1+\frac{1}{q'}-\frac{1}{r'} $, an operator in\\
$$Hom_G(W_{r'}^{1/r'}(S^1)/const,W_{q'}^{1/q'}(S^1)/const ),$$ \\
associated to
the principal state $[l]_p$ is proportional to the Hilbert transform
$$ H:W_{r'}^{1/r'}(S^1)/const\to W_{r'}^{1/r'}(S^1)/const ,$$
followed by the 
embedding $W_{r'}^{1/r'}(S^1)\hookrightarrow W_{q'}^{1/q'}(S^1)$.

\noindent\\
{\it Proof.}--- First, we notice that the Hilbert transform acts as an
isometric operator in $W_p^{1/p}(S^1)/const $ for all $p>1$. This follows at
once from the definition of the norm as
$$ \|u\|=\int_{\mathcal{H}^2}|\nabla f|^pd\mu_h ,$$
where $\Delta f=0$ and $f|S^1=u\ (\mbox{mod const}) $. We will prove the proposition by
a direct unimaginative computation. Let
$$ g(z)=\frac{z+z_0}{1+\bar z_0 z},\qquad |z_0|<1,|z|<1 ,$$
so that $g(0)=z_0$. Then the Jacobian of $g$ on the unit circle is
$$ \frac{1-|z_0|^2}{|z-z_0|^2} ,$$
so
$$ l(g)=\log (1-|z_0|^2 )-\log |z-z_0|^2 .$$
Let $\varphi: S^1\to \mathbb{R} $ be smooth. Then
$$ <\varphi, l(g)>=\int_{S^1} \varphi'(\theta)\cdot
[\log(1-|z_0|^2)-\log|e^{i\theta}-re^{i\varphi}|^2 ] d\theta,$$
where $z_0=re^{i\varphi}$. Obviously,
$$ \int_{S^1}\varphi'(\theta)\log(1-|z_0|^2)=0 ,$$
so
\renewcommand{\baselinestretch}{1.5}
\large
$$ \begin{array}{rcl}
    <\varphi,l(g)>&=&-\int_{S^1}\varphi(\theta)\cdot [\log
|e^{i\theta}-re^{i\varphi}|^2]' \\
        &=&-\int\varphi(\theta)\cdot \frac{2r\sin (\theta-\varphi) }{
1+r^2-2r\cos (\theta-\varphi) }\ .
    \end{array}
 $$
\renewcommand{\baselinestretch}{1.2}
\normalsize
As $z_0=re^{i\varphi}\underset{r\to 1}{\to}
 e^{i\varphi}$, this converges to
$$ -v.p.\int\varphi(\theta)\frac{2\sin (\theta-\varphi)}{2-2\cos
(\theta-\varphi) }=v.p.\int \varphi(\theta)\cdot \frac{1}{\mathit{tg}\
\frac{\theta-\varphi}{2}}=\pi H\varphi(\theta) $$
almost everywhere on $S^1$. The Proposition is proved, since smooth
functions are dense in $W_{r'}^{1/r'}(S^1)$.

Notice that since $H$ commutes with the action of $SO^{+}(1,2)$, for any cocycle
$m\in Z^1(G,W_p^{1/p}(S^1)/const)$, $Hm$ is also an cocycle. In
particular, $H[l]_p\in H^1(G,W_p^{1/p}(S^1)/const )$. We wish to compute a
corresponding operator in $Hom(W_{r'}^{1/r'}(S^1)/const \to
W_{q'}^{1/q'}(S^1)/const )$. Let $F$, as above, be a smooth map
$$ F:\mathcal{H}^2\to W_p^{1/p}(S^1)/const ,$$
satisfying $F(g^{-1}z)=R_g F(z)+l_p(g) $. For $v\in W_{r'}^{1/r'}(S^1)$ we
need to find a limit on the boundary of $<HF,v>$. But $H$ repects the
pairing $<\cdot,\cdot>$ and $H^2=-1$, so $<HF,v>=-<F,Hv>$, whose limit on
$S^1$ is $\pi H(-Hv)=\pi v$. We have proved the following lemma.

\noindent\\
{\it Lemma 12.4.}--- For $p,r',q'>1,\ \frac{1}{p}\ge
1+\frac{1}{q'}-\frac{1}{r'}$, an operator in\\
$$Hom_G(W_{r'}^{1/r'}(S^1)/const,W_{q'}^{1/q'}(S^1)/const),$$\\
associated with $\frac{1}{\pi}H[l]_p$, is the identity.

\noindent\\
{\it Theorem 12.5.}A.--- An element $v=H[l]_2\in \mathcal{H}_{2,g}$
does not depend on the choice of the lattice $G\hookrightarrow SO^{+}(1,2)$\\
B.--- The action of $\mathcal{M}ap_g$ in $\mathcal{H}_{2,g}$ fixes $v$.

\noindent\\
{\it Remark.}--- The Theorem is beyond doubt true for all $p>1$ and
not only $p\ge 2$, however I can't prove this at the moment of writing this
paper (July,1999). ({\bf Added January, 2000}). This is in fact true. The proof will appear in [Reznikov 10]).

\noindent\\
The vector $v$ is called a vacuum vector.

\noindent\\
{\it Proof.}--- Consider two embeddings $i_1,i_2: G\to SO^{+}(1,2)$ as
cocompact lattices and let $v_1,v_2$ be corresponding elements. We view
$v_1,v_2$ as elements of $H^1(G_r,H^1(G_l,l^2(G))$. Let $A_1,A_2$ be
associated operators
$$ A_1,A_2: H^1(G_l,l^{r'}(G))\to H^1(G_l,l^{q'}(G)) .$$
We know that $A_1=A_2=id$. It follows that an operator, associated with
$v_1-v_2$ is zero. We are going to show that $v_1-v_2$ is zero. Since
$$ v_1-v_2\in H^1(G_r,V),$$
where $V$ stands for $H^1(G_l,l^2(G))\simeq W_2^{1/2}(S^1)/const $, by a
result of [Korevaar-Schoen 1] cited above (Lemma 11.6) there exists a
harmonic section $F$ of the affine Hilbert bundle over $M=\mathcal{H}^2/G$
with fiber $V$ and nonodromy
$$ g\mapsto R_g(\cdot)+m(g) ,$$
where $m(g)$ is any cocycle, representing $v_1-v_2$. Let $v\in
W_{r'}^{1/r'}(S^1)/const $, then, denoting by $F$ again the lift of this
section on $\widetilde M=\mathcal{H}^2$, we see that $<F,v>$ is a harmonic
function such that
$$ \int_{\mathcal{H}^2}|\nabla (<F,v>)|^{q'}\ d\mu_h <\infty ,$$
and the trace of $<F,v>$ on $S^1$ is constant. It follows that $<F,v>$ is
constant itself, therefore ($v$ is arbitrary!) $F=w=const$ and
$$ m(g)=R_gw-w,$$
so $v_1-v_2=0$. This proves A. Now, if
$\phi\in\mathcal{M}ap_{g,1}=Aut(\pi_1(\Sigma^g))$, simply apply A to $i_1$ and
$i_1\circ \phi $.

We wish to compute $v$. $[l]_2$ is given by a cocycle
$$ g\mapsto -2\log |\beta-w| ,$$
$\beta\in S^1,\ w=g(0)$. This is equal to $2\log |1-\bar w\beta|$. The
latter function is a real part of $2\log (1-\bar wz)$ which is holomorphic
in $|z|\le 1$, so the Hilbert transform is $2Arg(1-\bar w\beta )$. This
means that a cocycle
$$ m(g)(\beta)=2Arg (W-\beta)\ (\mbox{\it{ mod const}}) $$
where $W=1/\bar w , w=g(0)$, represents $v$.

\noindent\\
{\it Theorem 12.6.}--- $H^1(\mathcal{M}ap_{g,1},H^1(G_l,l^p(G)))\not= 0$
for $p\ge 2$.

\noindent\\
{\it Proof.}--- We embed $G$ as a lattice in $SO^{+}(1,2)$ and identify
$H^1(G_l,l^p(G))$ and $W_p^{1/p}(S^1)/const $. We know that
$$ H^0(\mathcal{M}ap_g,H^1(G_r,W_p^{1/p}(S^1)/const)\ni v\not= 0.$$
Notice that $H^0(G_r,W_p^{1/p}(S^1)/const)=0$ since any $G$-invariant
harmonic 1-form in $\mathcal{H}^2$ has infinte $p$-energy. So in the
spectral sequence
$$ E_{i,j}^2: H^i(\mathcal{M}ap_g,
H^j(G_r,W_p^{1/p}(S^1)/const))\Longrightarrow
H^{i+j}(\mathcal{M}ap_{g,1},W_p^{1/p}(S^1)/const )$$
the second differential
$$ d_2: H^0(\mathcal{M}ap_g,H^1(G_r,W_p^{1/p}(S^1)/const )\Longrightarrow
H^2(\mathcal{M}ap_g,H^0(G_r,W_p^{1/p}(S^1)/const) $$
must be zero. Therefore the vacuum vector $v$ survives in $E^\infty$.

\noindent\\
It is plausible that, in fact,
$$ H^1(\mathcal{G}_1,W_p^{1/p}(S^1)/const )\not= 0\qquad (p>1)$$
for the group $\mathcal{G}_1$ of quasisymmetric homeomorphisms. ({\bf Added January, 2000}). This is in fact true. A formula  $$\Phi\to Arg \Phi^{-1} (\beta )- Arg (\beta)\  \mbox{mod const}$$ defines a cocycle of  $\mathcal{G}_1$ in $W_p^{1/p}(S^1)/const$ for any $p>1$. The proof will in [Reznikov 10]).  


\section{Equivariant mapping of the  Teichm\"uller Space, a space of quasifuchsian
representations and a space of all discrete representations into
$\mathcal{H}_{p,g} $}

\noindent
{\it Theorem 13.1.}A.--- A map which associates to a discrete
cocompact representation
$$ G\to SO^{+}(1,2) $$
its principal state
$$ [l]_p\in H^1(G_r,H^1(G_l,l^p(G)) $$
is an $\mathcal{M}ap_g$-equivariant map of the Teichm\"uller space
$\mathbf{T}_{6g-6}$ to $\mathcal{H}_{p,g}$ for all $p>1$. \\
B.--- Let $\varphi: G\to SO^{+}(1,3) $ be a discrete representation. Let
$\alpha_\varphi : S^1\to S^2 $ be the limit map of the boundaries
$$ S^1=\partial\widetilde\Sigma\to \partial \mathcal{H}^3=S^2 ,$$
defined in section 6, associated to $\varphi $. For $p>2$ let
$$ [l]_p\in H^1(SO^{+}(1,3),W_p^{2/p}(S^2)/const )$$
be the principle state. A map
$$ \varphi\mapsto A_\varphi \varphi^*[l]_p \in H^1(G_r,H^1(G_l,l^p(G))) ,$$
defined by first pulling back $[l]_p $ to $\varphi^*[l]_p \in
H^1(G,W_p^{2/p}(S^2)/const )$ and then applying the operator
$$ A_\varphi: W_p^{2/p}(S^2)/const\to W_p^{1/p}(S^1)/const ,$$
induced by $\alpha_\varphi $ and defined in section 6, is an
$\mathcal{M}ap_g$ equivariant map
$$ Hom_{discrete}(G,SO^{+}(1,3))/SO^{+}(1,3)\to \mathcal{H}_{p,g} $$
for all $p>2$.\\
C.--- A restriction of the map, defined in $B$ to
$$ Hom_{quasifuchsian}(G,SO^{+}(1,3))/SO^{+}(1,3) $$
is contained in $\mathcal{H}_{2,g}$.

\noindent\\
{\it Proof.}--- is already contained in section 6--12. We notice that
from the operator viewpoint the map of A sends any realization of $G$ as a
lattice in $SO^{+}(1,2)$ to a Hilbert transform of $W_p^{1/p}(S^1)/const$,
followed by an identification
$$ H^1(G_l,l^p(G))\simeq W_p^{1/p}(S^1)/const ,$$
which depends on the lattice. In other words, fix one lattice embedding
$$ \beta_0: G\to SO^{+}(1,2) .$$
Then any other lattice embedding
$$ \beta: G\to SO^{+}(1,2) $$
can be written as
$$ \beta(g)=\Phi_{\beta_0,\beta}\beta_0(g)\Phi_{\beta_0,\beta}^{-1} ,$$
where $\Phi_{\beta_0,\beta}\in \mathcal{G}_1 $ is a quasisymmetric map. Then
an operator, associated with $\beta$ is
$$ \Phi_{\beta_0,\beta}H\Phi_{\beta_0,\beta}^{-1} \in
Aut(W_p^{1/p}(S^1)/const ) .$$
This gives an $\mathcal{M}ap_g$-equivariant map
$$ \mathbf{T}_{6g-6}\to Aut_G(W_p^{1/p}(S^1)/const) $$
For $p=2$ one gets a map
$$ \mathbf{T}_{6g-6}\to Sp_G (W_2^{1/2}(S^1)/const )$$
because the Hilbert transform and $\mathcal{G}_1$-action are
symplectic(section 7), which can be described as follows. First, one embeds
$\mathbf{T}_{6g-6}$ in the universal Teichm\"uller space
$$ \mathbf{T}  =\mathcal{G}_1/SO^{+}(1,2) .$$
Then using the representation
$$ \mathcal{G}_1\to Sp(W_2^{1/2}(S^1)/const )$$
defined in section 7, one defines an embedding to $Sp/U $:
$$\mathbf{T}\to Sp(W_2^{1/2}(S^1)/const )/U $$
where $U$ is a group of operator in $Sp$ which commutes with $H$ seen as a
complex structure in $W_2^{1/2}(S^1)/const $. Finally, one uses the Cartan
embedding
$$ Sp/U\to Sp . $$

\noindent\\
{\it Theorem 13.2.}(Linearization of pseudoAnosov automorphisms )---
 Let $\phi\in \mathcal{M}ap_{g,1}
=Aut(\pi_1(\Sigma^g)) $ is a pseudoAnosov automorphism. Then for any $p>1$
there exists a nontrivial element $S_p\in \mathcal{H}_{p,g} $ with the
following properties:\\
1) for $p_1<p_2$, $S_{p_2}$ is an image of 
$S_{p_1}$, under the natural map
$\mathcal{H}_{p_1,g}\to \mathcal{H}_{p_2,g} $; \\
2) $S_p$ is invariant under $\bar{\phi} \in \mathcal{M}ap_g $; \\
3) there is a cocycle $\tilde l_p\in Z^1(G,W_p^{1/p}(S^1)/const) $ ,
representing $S_p$, such that for any $g\in G$
$$ \sum_{n\in \mathbb{Z}}\| \tilde l_p (g)\circ \Phi^n
\|_{W_p^{1/p}(S^1)/const} <\infty $$
where $\Phi: S^1\to S^1 $ is a quasisymmetric homeomorphism, associated with
$\phi$ (or, in other words,
$$ \sum_{n\in\mathbb{Z}}\|A_\varphi^m\tilde l_p(g) \|<\infty $$
where $A_\varphi \in Aut(H^1(G_l,l^p(G))$ is induced by $\phi$)

\noindent\\
{\it Proof.}--- is an immediate corollary of [Thurston 2]
(see also an exposition in [Otal 1] ), which shows that the mapping torus of any homeomorphism
$\Psi :\Sigma\to\Sigma$, representing $\varphi$ is a hyperbolic 3-manifold, Theorem
13.1, Theorem 9.1 and Theorem 9.3.
\\

It is plausible that such $S_p$ is unique up to a multiplier. Knowing $S_p$ is
essentially equivalent to knowing the hyperbolic volumes of all ideal simplices
with vertices on the limit curve $S^1\to S^2$.


\chapter{A theory of groups acting on the circle}

Our first main result in this Chapter is Theorem 1.7 which says, roughly,
that Kazhdan group cannot act on the circle. This general theorem draws a line
after many years of study and various special results concerning the actions
of lattices in Lie groups, see [Witte 1], [Farb-Shalen 1], [Ghys 1]. One can see here a historic parallel with a similar, but easier, general theorem of [Alperin 1] and
[Watatani 1] concerning Kazhdan groups acting on trees, which also followed
a study of the actions of lattices. Our technique is absolutely different
from the cited papers and uses a fundamental cocycle, introduced and studied
in section 1. We also use standard facts from Kazhdan groups theory [de la
Harpe-Valette 1].

In Sections 2,3 we quantize equivariant maps between boundaries of universal
covers, studied in Chapter I, Section 6. Our main tool is a harmonic map theory into
infinite-dimensional spaces, as developed in [Korevaar-Schoen 1], see also
[Jost 1]. In Section 4 we review some facts about Banach-Lie groups and
regulators. In Section 5 we describe a series of higher characteristic
classes of subgroups of $\mathcal{D}iff^{1,\alpha}(S^1)$. There are two
construction given. One uses an extension to a restricted linear group of a
Hilbert space of classes originally defined in [Feigin-Tsygan 1] for
infinite Jacobian matrices. Another construction 
uses the action of a restricted
symplectic group $Sp(W_2^{1/2}(S^1)/const )$ on the infinite-dimensional
Siegel half-plane. In both construction we use an embedding of
$\mathcal{D}iff^{1,\alpha}$ into a restricted linear group, by the unitary
and symplectic representation of $\mathcal{D}iff$, respectively. Using the
geometry of the Siegel half-plane, we prove that our classes have polynomial
growth.

There is a striking similarity between the theory of this Chapter and a
theory of symplectomorphism group, see Chapter IV, [Reznikov 2] and
[Reznikov 4]. Note that the extended mapping class group action is not
$C^{1,\alpha}$ smooth, so the results of this Chapter do not apply to this group. On the
other hand, $\mathcal{M}ap_g$ does act symplectically on a smooth compact
symplectic manifold.

\section{Fundamental cocycle}

By $\mathcal{D}iff^{1,\alpha}(S^1)$ we denote a group of
orientation-preserving diffeomorphisms with derivative in the H\"older space
$C^\alpha(S^1)$, which consists of functions $f$ such that
$$ |f(x)-f(y)|<c|x-y|^\alpha .$$

There is a series of unitary representations of
$\mathcal{D}iff^{1,\alpha}(S^1)$ in $L_{\mathbb{C}}^2(S^1,d\theta)$ given by
$$ (\pi(g)(f))(x)=f(g^{-1}x)\cdot[(g^{-1})'(x)]^{\frac{1}{2}+i\beta},\quad
\beta\in \mathbb{R} .$$

We will mostly consider $\beta=0$, in which case one has an orthogonal
representation in $L_{\mathbb{R}}^2(S^1,d\theta)$. An invariant meaning is,
of course a representation in half-densities on $S^1$. Now consider a
Hilbert transform $H$ as an operator in $L_{\mathbb{R}}^2(S^1,d\theta)$
given by a usual formula
$$ Hf(\varphi)=\frac{1}{\pi}v.p.\int_{S^1}\frac{f(\theta)}{\mathit{tg}\
\frac{\varphi-\theta}{2}}d\theta .$$

We wish to consider $[H,\pi(g^{-1})]$. This is a bounded operator in
$L^2(S^1,d\theta)$ given by an integral kernel which we are going to
compute. Notice that
$$ \frac{1}{\mathit{tg}\
\frac{\varphi-\theta}{2}}=\frac{2}{\varphi-\theta}+\textit{smooth kernel} .$$
A computation of [Pressley-Segal 1] shows that
$$
H[\pi(g)f](\varphi)=\frac{2}{\pi}v.p.\int_{S^1}\frac{d\theta}{\varphi-\theta}
f(g^{-1}(\theta))[(g^{-1}(\theta))']^{1/2}
+\textit{smooth kernel}\circ\pi(g) ,$$
so
\renewcommand{\baselinestretch}{1.5}
\large
$$\begin{array}{l}
(\pi(g^{-1})H\pi(g)f)(\varphi)=[g'(\varphi)]^{1/2}\cdot 
\frac{2}{\pi} v.p.\int_{S^1}
\frac{d\theta f(g^{-1}(\theta))
[(g^{-1}(\theta))']^{1/2}}{g(\varphi)-\theta}\\
+\pi(g^{-1})\circ \textit{smooth kernel}\circ \pi(g)
\end{array}.$$
\renewcommand{\baselinestretch}{1.2}
\normalsize
Letting $\theta=g(\eta)$ we have
\renewcommand{\baselinestretch}{1.5}
\large
$$\begin{array}{l}
(\pi(g^{-1})H\pi(g)f)(\varphi)=[g'(\varphi)]^{1/2}\frac{2}{\pi} v.p.
\int_{S^1}\frac{f(\eta)\cdot[g'(\eta)]^{1/2}}{g(\varphi)-g(\eta)}d\eta +\\
+\pi(g^{-1})\circ \textit{smooth kernel}\circ \pi(g) = \\
=\frac{2}{\pi}v.p.
\int_{S^1}\frac{[g'(\varphi)g'(\eta)]^{1/2}}{g(\varphi)-g(\eta)}f(\eta)d\eta
+\pi(g^{-1})\circ \textit{smooth kernel}\circ \pi(g)
\end{array}$$
\renewcommand{\baselinestretch}{1.2}
\normalsize
Finally,
\renewcommand{\baselinestretch}{1.5}
\large
$$\begin{array}{l}
[(\pi(g^{-1})H\pi(g)-H)](\varphi)=\\
=\frac{1}{\pi}\int_{S^1}\frac{[g'(\varphi)g'(\eta)]^{1/2}(\varphi-\eta)
-(g(\varphi)-g(\eta))}{(g(\varphi)-g(\eta))(\varphi-\eta)}f(\eta)d\eta+
\pi(g^{-1})\circ\textit{smooth kernel}\circ\pi(g)+\\
+\textit{smooth kernel}.
\end{array}\eqno(1.1) $$
\renewcommand{\baselinestretch}{1.2}
\normalsize

For a Hilbert space $\mathcal{H}$ and $p\ge 1$ we denote by $J_p (\mathcal{H} )$ a Shatten class of operators such that a sum of the $p$-th powers of their singular numbers converges. By $J_{p+} (\mathcal{H} )$ we mean the intersection of all
 $J_q (\mathcal{H} )$ with $q>p$.

Now recall that $g\in \mathcal{D}iff^{1,\alpha}(S^1) $. A following
proposition sharpens that of [Pressley-Segal 1] for
$\mathcal{D}iff^\infty(S^1)$:

\noindent\\
{\it Propositin 1.1.}\\
A. For $\alpha>1/2, \quad \pi(g^{-1})H\pi(g)-H\in J_2(L^2(S^1,d\theta)) $.\\
B. For $\alpha>0, \quad \pi(g^{-1})H\pi(g)-H\in
J_{1/\alpha +}(L^2(S^1),d\theta) .$

\noindent\\
{\it Proof.}--- As $\varphi-\eta\to 0$,
$$ \frac{[g'(\varphi)g'(\eta)]^{1/2}(\varphi-\eta)-(g(\varphi)-g(\eta))}{
(g(\varphi)-g(\eta))(\varphi-\eta)}<const\cdot (\varphi-\eta)^{\alpha-1} ,$$
so the kernel in (1.1) is in $L^2(S^1\times S^1, d\theta\otimes d\theta )$
for $\alpha>1/2$. This proves A.

To prove B we notice that by [Pietsch 1], the estimate on the kernel implies
that the operator lies in $\mathcal{J}_{1/\alpha +}$. Strictly speaking, the
conditions of [Pietsch 1] require $C^\infty$ smoothness off the diagonal,
whereas we have only the H\"older continuity, but the result stays true.
\\

Now notice that $GL(L^2(S^1,d\theta))$ acts in $J_p$ by conjugation. We
deduce the following

\noindent\\
{\it Proposition 1.2.}--- A map
$$ l:g\mapsto \pi(g)H\pi(g^{-1})-H $$
is a 1-cocycle of $\mathcal{D}iff^{1,\alpha}(S^1)$ in
$J_p(L^2(S^1,d\theta))$ for $p>1/\alpha $.
In particular, $l$ is a 1-cocycle of $\mathcal{D}iff^{1,\alpha}(S^1)$ in
$J_2$ for $\alpha>1/2$.

We will call $l$ a fundamental cocycle of $\mathcal{D}iff^{1,\alpha}(S^1)$.

\noindent\\
Now let $G$ be a subgroup of $\mathcal{D}iff^{1,\alpha}(S^1)$. We obtain a
class in $H^1(G,J_p(L^2(S^1,d\theta)) $ by restricting $l$ on $G$. We are
going to show that this class is never zero, except for completely
pathological actions of $G$ on $S^1$.

\noindent\\
{\it Proposition 1.3.}--- Let $G$ be a subgroup of
$\mathcal{D}iff^{1,\alpha}(S^1)$, $0<\alpha<1$. Suppose $p>1/\alpha$. If
$[l]\in H^1(G,J_p)$ zero, then the unitary action of $G$ in
$L_{\mathbb{C}}^2(S^1,d\theta)$ is reducible. Moreover, if $H^1(G,J_p)=0$
then $L_{\mathbb{C}}^2(S^1,d\theta)$ a direct sum of countably many closed 
invariant subspaces.

\noindent\\
{\it Proof.}--- If $[l]=0$ then there is $A\in J_p$ such that
$$ \pi(g)H\pi(g^{-1})-H=\pi(g)A\pi(g^{-1})-A $$
so that $[\pi(g),H-A]=0$. Since $H$ has two different eigenvalues with
infinitely-dimensional eigenspaces, $H-A\not= const\cdot Id $, so the action
of $G$ in $L_{\mathbb{C}}^2(S^1,d\theta) $ is reducible.

Next, consider an operator $R$ in $L^2(S^1,d\theta)$ with a kernel
$$ K(\varphi,\eta)=\frac{1}{|\mathit{tg}\ \frac{\varphi-\eta}{2}|} .$$

One sees immediately that $R$ is a self-adjoint unbounded operator. Repeating
the computation above, we deduce that
$\pi(g)R\pi(g^{-1})-R\in J_p $, so $\tilde
l(g)=\pi(g)R\pi(g^{-1})-R $ is another cocycle. If this
cocycle is trivial, then we get an unbounded self-adjoint operator
$R-A$ which commutes with the action of $G$. An application of the
spectral theorem shows that $L^2(S^1,d\theta)$ is a countable sum of
invariant subspaces.

\noindent\\
{\it Corollary 1.4.}--- A restriction of $l,\tilde l$ on $SO^{+}(1,2)$ is
not zero, for all $\alpha>0$.

\noindent\\
{\it Proof.}--- $SO^{+}(1,2)$ act in $L^2_{\mathbb C}(S^1,d\theta)$ as a
representation of principal series, which are irreducible.
\\

We now specialize for $\alpha=1/2$ and $p=2$. Since $[\tilde l]\in
H^1(SO^{+}(1,2),J_2)$ is nonzero, $\| \tilde l(g)\|_{J_2} $ is unbounded as a
function of $g$ [de la Harpe-Valette 1]. In fact, one has the following

\noindent\\
{\it Proposition 1.5.}--- Let $\pi: SO^{+}(1,2)\to U(H)$ be a unitary
representation and let $l: SO^{+}(1,2)\to H$ be a continuous cocycle. Suppose
$[l]\not= 0$. Then\\
A. For any cocompact lattice $G\subset SO^{+}(1,2)$, $[l]|G\not= 0$.\\
B. $\| l(g^n)\|$ is unbounded as $n\to\infty $ for any hyperbolic $g$. \\
C. $\|l(\gamma ^n)\|$ is unbounded as $n\to \infty $ for any parabolic $\gamma \not=1$.

\noindent\\
{\it Proof.}--- Let $V\subset SO^{+}(1,2) $ be compact 
and such that $V\cdot
G=SO^{+}(1,2)$. For $v\in V,g\in G$ we have
$$l(vg)=\pi(v)l(g)+l(v) ,$$
so $\| l(vg)\|\le \| l(g)\|+\|l(v)\|$. If $l|G$ is bounded, then so is $l$.
This proves A.  Next, let $P$ be the image of $SO^{+}(1,2)/K$ under Cartan
embedding, where $K$ is a maximal compact subgroup. By the same reason as
above, $l|P$ is unbounded. Let $S^1\subset P$ be a nontrivial orbit of $K$
in $P\approx \mathcal{H}^2 $. Notice that $P$ is closed under raising into
an integral power and there is a compact $V\subset SO^{+}(1,2)$ such that
$$ P\subseteq \bigcup_{n\ge 1}(S^1)^n\cdot V $$
where $(S^1)^n$ is an image of $S^1$ under raising to $n$-th power. We
deduce that $l|\cup_{n\ge 1}(S^1)^n $ is unbounded. Let $\gamma\in S^1$.
Then any element in $(S^1)^n$ is of the form $k\gamma^nk^{-1}$, $k\in K$, so
$$ \|l(k\gamma^n k^{-1} )\| \le \|l(k)\|+\|l(k^{-1})\|+\| l(\gamma^n) \| $$
So $\| l(\gamma^n)\| $ is unbounded. 
Since $\gamma$ can be any hyperbolic element, B
follows. Notice that we proved that $\| l(g_k)\|$ is unbounded for any
sequence $g_k\in P$, which escapes all compact sets. Now let $g\in SO^{+}(1,2)$
be parabolic $\not=1$, and let $\tau$ be the involution fixing $K$. Then
$\tau (g^n)\cdot g^{-n} \in P$ and escapes all compact sets, so $\| l[(\tau
g^n)\cdot g^{-n} ]\| $ is unbounded. It follows that either $\| l(\tau
g^n)\|$ or $\| l(g^{-n})\| $ is unbounded. But all parabolics are conjugate
in $SO^{+}(1,2)$, so $C$ follows.

\noindent\\
{\it Proposition 1.6.}--- Let $G\subset
\mathcal{D}iff^{1,\alpha}(S^1),\quad \alpha>1/2 $. Suppose that $G$ contains
an element $g$ which is conjugate in
$\mathcal{D}iff^{1,\alpha}(S^1)$ to a hyperbolic or a nontrivial parabolic fractional-linear transformation. Then $[l]|G\not= 0$ in $H^1(G,J_2)$.

\noindent\\
{\it Proof.}--- Any such $g$ is conjugate in
$\mathcal{D}iff^{1,\alpha}(S^1)$ to an element $g'\in SO^{+}(1,2) $ for which
$\| l({g'}^n) \|$ is unbounded, so $\| l(g^n)\| $ is unbounded as well.
\\

We are ready to formulate the main result of this section.

\noindent\\
{\it Theorem 1.7.}--- Let $G\subset \mathcal{D}iff^{1,\alpha}(S^1),\
\alpha>1/2 $. Suppose that either \\
1) a natural unitary action ($\beta=0$) of $G$ in $L^2(S^1,d\theta) $ given
by
$$ \pi(g)(f)(\varphi) =f(g^{-1}(\varphi) )\cdot [(g^{-1}(\varphi))']^{1/2} ,$$
is irreducible or is a direct sum of finitely many irreducible factors, or\\
2) $G$ contains an element, conjugate in $\mathcal{D}iff^{1,\alpha}(S^1)$
to a hyperbolic fractional-linear transformation, or\\
3) $G$ contains an element, conjugate in $\mathcal{D}iff^{1,\alpha}(S^1)$
to a parabolic ($\not= 1$) fractional-linear transformation, or\\
4) $$ \sup_{g\in G}{\int\int}_{S^1} \left[
\frac{\sqrt{g'(\varphi)g'(\eta)}(\varphi-\eta) - (g(\varphi)-g(\eta)) }{
(g(\varphi)-g(\eta))(\varphi-\eta) }\right]^2\ d\varphi d\eta =\infty $$
Then $G$ is not Kazhdan.

\noindent\\
{\it Proof} follows from the formula (II.1.1), Proposition 1.3,
Proposition 1.5 and Proposition 1.6.

\section{Construction of $N=2$ quantum fields with lattice symmetry}

It is possible that the physical time-space is discrete. Correspondingly, in
the axiomatic quantum field theory it is possible that the fields must yield
invariance not under the whole Poincar\'e group, but only under a lattice in
it. See [Michailov 1], [Belavin 1] in this respect. We are going to
construct mathematical objects, which yield such invariance on one hand, and
quantize the equivariant measurable maps considered in I.6.3, on the other.

\noindent\\
{\it Theorem 2.1.}--- Let $G$ be a cocompact lattice in $SO^{+}(1,2)$. Let
$\mathcal{H}=L_{\mathbb{R}}^2(S^1,d\theta) $ with the orthogonal action
$\pi$,
corresponding to $\beta\in \mathbb{R}$. Then there exists a measurable map
to the space of bounded operators $$ S^1\overset{\rho}{\to} \mathcal{B}(H) $$
with the following properties.\\
1) Equivariance: for $s\in S^1$ and $g\in G$
$$ \rho(gs)=\pi(g)\rho(s)\pi(g^{-1}) $$
almost everywhere on $S^1$. \\
2)One has $$\int_{S^1}(\rho(s)-H)\psi(s)ds\in J_2 $$ for $\psi\in C^\infty (S^1) $. \\
3)There exists $J\in J_2(\mathcal{H})$ such that $\rho (s)$ is a weak nontangential limit
$$\rho (s)= \lim_{g\to s} \pi(g)(H+J) \pi (g^{-1})$$
as $g\in G$ converges nontangentially to $s\in S^1 =\partial G$ a.e. on $S^1.$

\noindent\\
{\it Proof.}--- As a Hilbert space with orthogonal $G$-action,
$J_2=L^2(S^1\times S^1,d\theta\otimes d\theta ) $. By the proof of Lemma
I.11.3, $G$ does not have almost invariant vectors in $J_2$. Let
$\Sigma=\mathcal{H}^2/G $ and let $E$ be a flat affine vector bundle over
$\Sigma $ with a fiber $J_2$ and monodromy
$$ g\mapsto Ad\pi(g) +l(g) .$$
Then by a result of [Korevaar-Schoen 1], and [Jost 1] (lemma I.11.6), there
exists a harmonic map
$$ \tilde f: \mathcal{H}^2\to J_2 $$
satisfying
$$ \tilde f(gx)=\pi(g)\tilde f(x)\pi(g^{-1})+l(g) $$
Consider $f(x)=\tilde f(x) +H $. Then
$$ f(gx)=\pi(g) f(x)\pi(g^{-1}), $$
in particular, $\| f(x)\| $ is bounded in operator norm. An operator version
of Fatou theorem [Naboko 1 and references therein ] shows that $f$ has nontangential limit
values a.e. on $S^1$, say $\rho(s)$. Obviously, $\rho$ is $G$-invariant. On
the other hand, $\tilde f$ is a Bloch harmonic $J_2$-valued function, that
is,
$$ \sup_{x\in \mathcal{H}^2} \| \nabla \tilde f\|_{J_2} <\infty .$$
It follows that $\| \tilde f(w) \|_{J_2}<c\cdot \log(1-|w|)$,
$w\in B^2= \mathcal{H}^2 $. This implies by a standard argument
(see e.g. [Gorba\v cuk 1] that $\tilde f$ has a limit on $S^1$ as an element
of $\mathcal{D}'(S^1,J_2) $. So for $\psi\in C^\infty(S^1) $,
$$ \int_{S^1}(\rho-H)\psi \in J_2 ,$$
which proves the Theorem.

\noindent\\
{\it Remarks.}--- 1)As was mentioned above, the invariant meaning of the
representation $\pi$ is that $L^2(S^1,d\theta) $ should be regarded as a
space of half-densities. Correspondingly, an integral operator is defined by
a kernel which is a half-density on $S^1\times S^1$ of the type
$K(\varphi,\eta )(d\varphi d\eta)^{1/2} $. If $K(\varphi,\eta)$ is smooth
and has a zero of second order on the diagonal $\Delta\subset S^1\times S^1$,
then one has an invariant definition of its residue or second derivative,
which is a quadratic differential. A direct computation which we leave to the reader  shows that for $g\in
\mathcal{D}iff^{\infty}(S^1) $\\:
1) $l(g)=\pi(g)H\pi(g^{-1})-H$ is given by a kernel which has a zero of
second order on $\Delta$;\\
2) a corresponding residue $S(g)$ is the Schwartzian of $g$.

This shows that $l(g)$ is a quantization of the Schwartzian cocycle. The
operator field $\rho(s)$ of Theorem 1.8 seems therefore to be related to
objects axiomatized, but not constructed, in [Belavin-Polyakov-Zamolodchikov
1].
2)The Theorem and the proof stay valid for any representation
$$\varphi : G\to \mathcal{D}iff^{1,\alpha }(S^1 ),$$
$\alpha >1/2
$, such that the action on $S^1 \times S^1$ satisfies the very mild conditions of Lemma I.11.3.

\section{Construction of $N=3$ quantum fields with lattice symmetry}

A theory developed have for $\mathcal{D}iff(S^1) $ does not generalize to
$\mathcal{D}iff(S^n), n\ge 2$. The reason is that the action of
$\mathcal{D}iff(S^1)$ on $S^1$ is conformal. There are two ways to
generalize various aspects of the theory to higher dimensions, by either
considering $SO^{+}(1,n)$ acting on $S^{n-1}$ or, very surprisingly, a group of
symplectomorphisms of a compact symplectic manifold $M$ (see Chapter IV).
Here we consider the action of $SO^{+}(1,3)\simeq PSL_2 (\mathbb{C})$ on $S^2$.
We set $d(x,y)$ to be a spherical distance in $S^2$. Let $d\theta$ denote
the spherical measure and let $\mathcal{H}=L^2(S^2,d\theta )$. For $g\in
SO^{+}(1,3)$ let $\mu_g(x)$ denote a conformal factor, that is $\mu_g^2(x)$ is a
Jacobian of $g$ with respect to $d\theta $. A formula
$$ \pi(g)f(x)=f(g^{-1}(x))\cdot \mu_{g^{-1}}^{1+i\beta}(x),\quad
\beta\in\mathbb{R} ,$$
defines a unitary representation of $SO^{+}(1,3)$ in $\mathcal{H}$. Now we introduce
an operator $H$ with the kernel
$$ K(\varphi,\theta)=\frac{1}{d^2(\varphi,\eta)} .$$
This operator is self-adjoint and unbounded. Our goal is to compute
$$\pi(g)H\pi(g^{-1})-H =l(g). $$

\noindent\\
{\it Proposition 3.1.}--- $l(g)\in J_2 $ for all $g\in SO^{+}(1,3) $ and
$\beta=0 $.

\noindent\\
{\it Proof.}--- A direct computation. One needs to show that as
$d(x,y)\to 0$,
$$ d^2(g(x),g(y))-\mu_g(x)\mu_g(y)d^2(x,y) $$
is of order $d^4(x,y)$. In other words, for a fractional-linear
trangformation $g$ of $\mathbb{C}$ one needs to show that as $x\to y$, $Im\
x,Im\ y>0,\ g(x)=x$,
$$ \left| \frac{g(x)-g(y)}{g(x)-\overline
{g(y)}}\right|^2-|g'(x)||g'(y)|\frac{Im\ y}{Im\ g(y)} \left|
\frac{x-y}{x-\bar y} \right|^2 $$
is of order $|x-y|^4$. This verifies the result for hyperbolic metric
instead of spherical metric, which is of course equivalent. One
computes directly using Taylor series for holomorphic function $g$.

Now arguing as in section 2 we arrive at the following result.

\noindent\\
{\it Theorem 3.2.}--- Let $G$ be a cocompact lattice in $SO^{+}(1,3)$. Let
$\mathcal{H}=L_{\mathbb{R}}^2(S^2,d\theta) $ with orthogonal action of $G$
corresponding to $\beta=0$. Then there exists a harmonic map
$$\mathcal{H}^3\overset{\psi}{\to} J_2(\mathcal{H}) $$
with the property that $z\mapsto \psi(z)+H $ is
equivariant:
$$ \psi (gz)+H= \pi(g)(\psi (z)+H) \pi(g^{-1}) $$
for all $g\in G$ and $z\in \mathcal{H}^3 $.

Since $H$ is unbounded, the boundary value of $\psi(z)+H$ does not exist as
a measurable map to the space of bounded operators. It is possible that
there is a more clever choice of a conformally natural singular integral
operator which is bounded, but I don't know how to do it. Note in this
respect that there is a very different realization of an orthogonal representation of
$SO^{+}(1,3)$ in the space of  functions on $S^2$ , discovered in [Reznikov 1].
Namely, look at the natural action of $SO^{+}(1,3) $ on smooth half
co-densities, that is, sections of $\sqrt{\varLambda^2 T S^2} $. Using the
spherical metric, we can identify this space with $C^\infty (S^2)$. Then the
above-mentioned action leaves invariant a nonnegative quadratic form
$$ Q(f)=\int_{S^2}((\Delta f)^2-2|\nabla f|^2)darea $$
whose kernel consists of constants and linear functions. It is possible that
there are $G$-equivariant quantum fields valued in operators acting in the
associated Hilbert space.

\section{Banach-Lie groups and regulators: an overview}

A Banach-Lie group is a Banach manifold with a compatible group structure.
Usual Lie theory largely extends to this case. In particular, if
$\mathcal{G}$ is a Banach-Lie group and $\mathfrak{g}$ its Banach-Lie
algebra, then a continuous n-cocycle on $\mathfrak{g}$ defines a
left-invariant closed form on $\mathcal{G}$, so that one has a homomorphism
$$ H_{cont}^n(\mathfrak{g},\mathbb{K})\to H_{top}^n(\mathcal{G},\mathbb{K})$$
where $H_{top}^* $ is a cohomology of a topological space. In [Reznikov 2]
we defined $\mathbb{K}$-homotopy groups of a Lie algebra, so that there is a
map
$$ \pi_i(\mathcal{G})\otimes \mathbb{K}\to \pi_i (\mathfrak{g}) $$
which in the case $\mathcal{G}=SL_n(C^\infty(M))$, $M$ a compact manifold,
$n>>1$, reduces to the Chern character
$$ K_i^{top}(M)\to
HC_i(C^\infty(M))=\Omega^i(M)/d\Omega^{i-1}(M)\oplus
H^{i-2}(M,\mathbb{K})\oplus\cdots $$
($\mathcal{G}$ is not a Banach-Lie group but a Frech\'et-Lie group in this
case). More interesting is a secondary class (=regulator) map. Define
an algebraic $K$-theory of $\mathcal{G}$ as
$$ K_i^{alg}(\mathcal{G})=\pi_i ((B\mathcal{G}^{\delta})^+ ) $$
and the augmented $K$-theory as a kernel of the map $K_i^{alg}\to K_i^{top}$:
$$ 0\to \overline K_i^{alg}(\mathcal{G})\to K_i^{alg}(\mathcal{G})\to
\pi_i(B\mathcal{G})=\pi_{i-1}(\mathcal{G}).$$
Then the regulator map is a homomorphism
$$ r: \overline K_i^{alg}(\mathcal{G})\to coker
(\pi_i(\mathcal{G})\otimes\mathbb{K}\to \pi_i (\mathfrak{g})). $$
Lifting this map to cohomology, that is , constructing a map
$$H^{\ast}_{cont}(\mathfrak{g},\mathbb{K})\to H^{\ast}(\mathcal{G}^{\delta},\mathbb{K})$$
meets obstructions described in the van Est spectral sequence. If $\mathcal{K}\subset \mathcal{G}$ is a closed subgroup such that $\mathcal{G} / \mathcal{K}$ is contractible, then these obstructions vanish and one gets a map
$$H^{\ast}_{cont}(\mathfrak{g},\mathfrak{k})\to H^{\ast}(\mathcal{G}^{\delta})$$
given explicitly by a Dupont-type construction [Dupont 1].
This is essentially the same as geometric construction of secondary classes
of flat $\mathcal{G}$-bundles, described in [Reznikov 3]. In case
$\mathcal{G}=SL_n(C^\infty (M))$ this gives a usual regulator map in
algebraic $K$-theory. However, for various diffeomorphism groups one
construct new interesting classes. For symplectomorphism groups two series
of classes, mentioned in the Introduction to Chapter 4, were constructed in
[Reznikov 2] and [Reznikov 4], and a new class associated to a Lagrangian
submanifold, will be constructed in Chapter 4. The symmetric spaces for
$Sympl(M)$, used in [Reznikov 2] are sort of continuous direct products of
finite-dimensional Siegel upper half-planes. On the other hand, a symmetric
space which we will use in this chapter to construct classes in
$H^*(\mathcal{D}iff^{1,\alpha}(S^1))$ is an infinite-dimensional Siegel
half-plane. The trouble is, however, that, for a compact manifold $Y$, (say,
$S^1$) a group of diffeomorphisms of finite smoothness, like
$\mathcal{D}iff^k(Y)$, is not a Banach-Lie group: the multiplication from
the right is not a diffeomorphism (the multiplication from the left is).
This is neatly explained in [Adams-Ratiu-Schmid 1]. Luckily, to construct secondary classes we
only use the fact that the multiplication from the left is a diffeomorphism.

\section{Charateristic classes of foliated circle bundles}

As is well known, the continuous cohomology of $\mathcal{D}iff^\infty (S^1)$
is generated by the Euler class and by the integrated Godbillon-Vey class
[Geldfand-Fuks 1], [Fuks 1 and references therein]. Moreover, the square of
the Euler class is zero. This already shows that the degree of smoothness is
crucial. For if one considers the action of the extended mapping class group
$$ \mathcal{M}ap_{g,1}\hookrightarrow \mathcal{G}_1\hookrightarrow Homeo(S^1)
,$$
then the pull-back of the Euler class has nonzero powers to a degree which
goes to infinity with $g$ [Miller 1], [Morita 1], [Mumford 1]. It appears
that the scarcity of the cohomology of $\mathcal{D}iff(S^1)$ is a
consequence of an (artificial) restriction of excessive degree of
smoothness. Notice that the proofs in [Fuks 1] depend hopelessly on
$C^\infty$- smoothness. We will give two constructions of a series of new
classes in $H^*(\mathcal{D}iff^{1,\alpha}(S^1))$, $0<\alpha<1$ using both
the unitary representation in $L^2(S^1,d\theta)$ and the symplectic
representation in $Sp(W_2^{1/2}(S^1)/const)$. As in the case of the powers of
the Euler class, a nonvanishing of these classes is an obstruction to 
smoothability, i.e. to a conjugation to a subgroup of
$\mathcal{D}iff^\infty(S^1)$. We will also prove that our classes are of
polynomial growth if $\alpha>1/2$. A related result (but not the argument) for $C^\infty$
Gelfand-Fuks cohomology in all dimensions is to be found in
[Connes-Gromov-Moscovici 1]. Both in spirit and technology, the construction
of the classes in $H^*(\mathcal{D}iff^{1,\alpha}(S^1))$ resembles our
construction of a series of classes in $H_{cont}^k(Sympl(M),\mathbb{R})$,
$k=2,6,10,\cdots $, where $M$ is a compact symplectic manifold and
$Sympl(M)$ is its symplectomorphism group [Reznikov 4].

We start with the construction using the unitary representation. By
Proposition 1.1, $\pi(g)H\pi(g^{-1})-H \in J_p $ where $g\in
\mathcal{D}iff^{1,\alpha}(S^1),\quad p>1/\alpha $, $\pi$ is a unitary
action in $L_{\mathbb{C}}^2(S^1,d\theta )$, and $H$ is a complexification of
the Hilbert transform. That is $H(e^{in\theta})=sgn(n)\cdot e^{in\theta} $.
The group of $\Phi\in GL(\mathcal{H})$,
$\mathcal{H}=L_{\mathbb{C}}^2(S^1,d\theta)$ such that $\Phi H\Phi^{-1}-H\in
J_p$ will be denoted $GL_{J_p}(\mathcal{H})$, following [Pressley-Segal 1].
The unitary subgroup of $GL_{J_p}(\mathcal{H})$ is denoted
$U_{J_p}(\mathcal{H})$. Let $\mathcal{H}_+,\mathcal{H}_-$ be the eigenspaces
of $H$ with eigenvalues $+1$ and $-1$ respectively. By $Gr_{J_p}(\mathcal{H})$ we
denote the restricted Grassmanian $U_{J_p}/U(\mathcal{H}_+)\times
U(\mathcal{H}_-) $. Then $Gr_{J_p}(\mathcal{H})$ is a Banach manifold,
modelled by the Banach space $J_p$. The Banach-Lie group
$GL_{J_p}(\mathcal{H})$ acts smoothly on $Gr_{J_p}(\mathcal{H})$. On the
other hand, though $\mathcal{D}iff^{1,\alpha}(S^1)$ is a group and a Banach
manifold, it is not a Banach-Lie group [Adams-Ratiu-Schmid 1].  However, multiplication from
the left $L_g(h)=gh$ is a diffeomorphism (but not a multiplication from the
right). The embedding
$$ \mathcal{D}iff^{1,\alpha}(S^1)\to U_{J_p}\to GL_{J_p}(\mathcal{H})$$
is not continuous. However, an induced action of
$\mathcal{D}iff^{1,\alpha}(S^1)$ on $Gr_{J_p}(\mathcal{H})$ is smooth
[Pressley-Segal 1].

We will introduce a series of $U_{J_p}$-invariant differential forms on
$Gr_{J_p}(\mathcal{H})$. These forms induce cohomology classes in the Lie
algebra cohomology $H^*(Lie(U_{J_p}))$, extending the classes introduced in
[Feigin-Tsygan 1] for the Lie algebra of Jacobian matrices. Notice that a
tangent space to the origin of $Gr_{J_p}(\mathcal{H})$ can be identified with
matrices of the form
$$
C=\begin{pmatrix}
    0 & B \\
    A & 0
 \end{pmatrix}
$$
where $A\in J_p(\mathcal{H}_+,\mathcal{H}_-)$ and $B\in
J_p(\mathcal{H}_-,\mathcal{H}_+)$. Let $C_1,\cdots,C_{2k},\quad (k \mbox{
odd} )$ be a collection of such matrices. Define
$$ \mu_k(C_1,\cdots,C_{2k})=\sum_{\sigma\in S_{2k}}sgn(\sigma) P_k(\rho
(C_{\sigma(1)},C_{\sigma(2)}),\cdots, \rho(C_{\sigma(2k-1)},C_{\sigma(2k)})$$
where $P_k$ is the $k$-th invariant symmetric functions of $k$ matrices,
which is a polarization of $tr\ A^k$ (not an elementary symmetric
polynomial, as in [Fuks 1] ).  Now, $\rho(C_1,C_2)$ is defined as follows: let
$\pi(C)$ is the left upper corner of $C$, i.e. an operator in
$\mathcal{B}(\mathcal{H}_+)$. Then $\rho(C_1,C_2)=\pi
([C_1,C_2])-[\pi(C_1),\pi(C_2)]$. The "meaning" of $\pi$ is that of a connection
of a principal bundle on something like the classifying space of the Lie
algebra $Lie(GL_{J_p})$, and of $\rho$ is that of the curvature of this
connection. Then $\mu_k$ becomes a characteristic class, somewhat analogous to the characteristic classes  in the standard Chern-Weil theory. Notice that $\mu_k$ is defined for all $k\ge
[1/\alpha]+1$. In [Feigin-Tsygan 1], $\rho(C_1,C_2)\in
\mathfrak{gl}(\infty,\mathbb{K}) $ and $\mu_k$ is defined for all $k$. The
form $\mu_2$ defines the famous ``Japanese cocycle'', [Verdier 1].

\noindent\\
{\it Lemma 5.1.}--- $\mu_k$ is $U_{J_p}$-invariant and closed.

\noindent\\
{\it Proof.}--- The invariance is obvious. The proof of closedness is
standard and left to the reader, see the remarks above and [Feigin-Tsygan 1].

\vspace{3mm}
Pulling back to $\mathcal{D}iff^{1,\alpha}(S^1)$ (this is possible by the
remarks made above) we obtain a left-invariant closed differential form on
$\mathcal{D}iff^{1,\alpha}(S^1)$. Pulling back to the universal cover
$\widetilde{\mathcal{D}iff}^{1,\alpha}(S^1)$, we obtain a left-invariant closed
differential form $\tilde\mu_k$ on $\widetilde{\mathcal{D}iff}^{1,\alpha}(S^1)$.
A following theorem follows.

\noindent\\
{\it Theorem 5.2.}--- The secondary characteristic class,
corresponding to $\tilde\mu_k$ is a well-defined class $r(\tilde\mu_k)$ in
$H^{2k}([\widetilde{\mathcal{D}iff}^{1,\alpha}]^\delta,\mathbb{R} )$.

\noindent\\
{\it Proof.}--- $\widetilde{\mathcal{D}iff}^{1,\alpha}(S^1)$ is
contractible.
\\

Notice that for $\alpha>1/2$ the class $\mu_1\in
H^2(\mathcal{D}iff^{1,\alpha}(S^1))$ is defined, which is just the
integrated Godbillon-Vey class.

Our second construction uses the symplectic action. For simplicity, we only
treat the case $\alpha>1/2$. Recall (Corollary I.7.3) that $\mathcal{G}_1$
acts symplectically in $V=W_2^{1/2}(S^1)/const$. Restricting on
$\mathcal{D}iff^{1,\alpha}(S^1)$, we obtain a representation
$$ \mathcal{D}iff^{1,\alpha}(S^1)\overset{\pi}{\to}Sp(V).$$
Let $H$ be the Hilbert transform in $V$, normalized such that $H^2 =-1$.
Denote by $Sp_{J_p}$ a subgroup of $A\in Sp(V)$ such that $[A,H]\in J_p$.
Denote $U=U(V)$ the unitary group of such $A$ that $[A,H]=0$. Denote
$$ X=Sp_{J_p}/U $$
a restricted Siegel half-plane. This is a Banach contractible manifold
[Palais 1]. For $p=2$ this is a Hilbert manifold with canonical
$Sp_{J_2}$-invariant Riemannian metric of nonpositive curvature. The metric
is defined as follows. The tangent space $T_H(X)$ is identified with
operators $A$ such that $A\in Lie(Sp_{J_2})$ and $AH=-HA$. It follows that
$A\in J_2$, and $A=A^*$. Then the metric is defined as $trA^2$. This
definition is dimension-free and so the proof that the curvature is
nonpositive follows from the explicit formulae, as in finitely-dimentional
case.

\noindent\\
{\it Lemma 5.3.}--- For $\alpha>1/2$, $\pi
(\mathcal{D}iff^{1,\alpha}(S^1))\subset Sp_{J_2}(V) .$

\noindent\\
{\it Proof.}--- We will use the computation of [Segal 1]. Let $g\in
\mathcal{D}iff^{1,\alpha}(S^1)$. We need to show that
$$ S=\sum_{n,m>0} \frac{m}{n} \left| 
	\int_{S^1} e^{i(ng(\theta)+m\theta)}\ d\theta \right|^2
<\infty .$$
As in [Segal 1] we have, using a trick of Kazhdan,
$$ S=\sum_{N=1}^{\infty} \sum_{m=1}^{N-1} \frac{m}{n} \left|\int_{S^1}
e^{iN\varphi }\cdot [g_{\beta}^{-1} ]'(\varphi) d\varphi \right|^2 ,$$
where $\beta=\frac{n}{N},\ n=N-m,\ g_\beta (\theta)=\beta g(\theta)
+(1-\beta ) \theta, \ \theta\in S^1=\mathbb{R}/2\pi\mathbb{Z} $. For $0\le
\beta \le 1$, $g_{\beta}^{-1}$ are uniformly in
$\mathcal{D}iff^{1,\alpha}(S^1) $ with $\alpha>1/2$, so
$$ \int_{S^1}e^{iN\varphi}[g_\beta^{-1}]' (\varphi)\ d\varphi \le const\cdot
N^{-\alpha} \cdot c_N $$
with $\sum_{N=1}^{\infty} c_N^2<\infty $. Since
$\sum_{m=1}^{N-1}\frac{m}{n}\sim \log N $, we have
$$ S\le const\cdot \sum_{N=1}^{\infty} N\log N N^{-2\alpha}\cdot c_N^2
<\infty .$$

Now let $k$ be odd, $A_1,\cdots, A_{2k}\in T_H X $ and
$$ \nu_k (A_1,A_2,\cdots,A_{2k})= $$
{\it Lemma 5.4.}--- $\nu_k$ is closed and $Sp_{J_2}$-invariant.

\noindent\\
{\it Proof.}--- is identical to the finite-dimensional case [Borel 1].

\noindent\\
{\it Theorem 5.5.}--- The secondary characteristic class,
corresponding to $\nu_k$ defines an element $r(\nu_k)$ in
$H_{cont}^{2k}(Sp_{J_2}(V))$ and in
$H^{2k}([\mathcal{D}iff^{1,\alpha}(S^1)]^{\delta},\mathbb{R} ) $, $\alpha >1/2$. All these
classes are of polynomial growth.

\noindent\\
{\it Proof.}--- Only the last statement needs a proof. For
$x_0,\cdots, x_s \in X$ denote a geodesic span $\sigma (x_0,\cdots,x_s )$ in
the following inductive way: $\sigma (x_0,x_1)$ is a geodesic segment
joining $x_0$ and $x_1$ and $\sigma(x_0,\cdots,x_s)$ is a union of geodesic
segments joining $x_0$ and points of $\sigma (x_1,\cdots,x_s)$. By standard
comparison theorems $Vol_s(\sigma (x_0,\cdots,x_s))\le const\cdot
[\max_{0\le i\le j\le s}\rho (x_i,x_j)]^s $, where $\rho (\cdot,\cdot) $ is
the distance function (this is where we use  non-positive curvature). By
[Dupont 1], $r(\nu_k ) $ can be represented by a cocycle
$$ g_1,\cdots,g_{2k}\mapsto \int_{\sigma
(x_0,g_1x_0,g_1g_2x_0,\cdots,g_1,g_2\cdots g_{2k}x_0)} \nu_k $$
where $g_i\in Sp_{J_2} $ and $x_0\in X$ is fixed. Since $\nu_k$ is uniformly
bounded, the result follows.

We will give an independent proof of polynomial growth of $\mu_2\in
H^2(\mathcal{D}iff^{1,1}(S^1))$. Let $Var(S^1)$ be a space of functions of
bounded variation on $S^1$ mod constants. Then for $f_1,f_2\in Var(S^1)$,
$$ \int_{S^1} f_1\cdot d\ f_2\le c\| f_1\|_{Var}\cdot \| f_2\|_{Var} .$$
Now, $Homeo(S^1)$ acts isometrically in $Var(S^1)$ and there is a cocycle
$\psi\in H^1(\mathcal{D}iff^{1,1}(S^1),Var) $ given by $g\mapsto \log (g^{-1})'$.
By an formula of Thurston, $\mu_2$ can be represented as
$$ \int_{S_1}\psi(g_1)\ d\psi(g_2) .$$
The result now follows from Lemma I.1.1. For $\mu_2$ as a class in
$H^2(\mathcal{D}iff^{\infty}(S^1))$ see also [Connes-Gromov-Moscovici 1].

\section{Examples}

A typical example of a group in $\mathcal{D}iff^{1,\alpha}(S^1) $ is a
following one. Let $K\subset \mathbb{R}$ be a subfield (i.e., a 
number field). By
$S^1(K)$ denote $K$-rational points of $S^1\subset
\mathbb{R}^2$. Define $G_{K}$ as a group of $C^1$-diffeomorphism
$g$ such that there are points $x_0,\cdots,x_n=x_0\in S^1(K)$ in
this order such that $g_k=g_{[x_k,x_{k+1}]}$ is a restriction of an element
of $PSL_2(K)$. The $C^1$-condition simply means that
$g_k'(x_{k+1})=g_{k+1}'(x_{k+1})$. Then automatically $G_{K}\subset
\mathcal{D}iff^{1,1}(S^1)$. Groups of this type, or rather their obvious analogues which act by piecewise-affine transformations on $S^1$ viewed as $\mathbb{R} / \mathbb{Z}$  appeared in [Thompson 1],
[Greenberg-Sergiesku 1,2], [Brown-Georghegan 1], etc. where various properties were
studied. The ``proper'' Thompson group can be smoothed , that is, embedded in
$\mathcal{D}iff(S^1)$ [Ghys-Sergiesku 1] so that the Theorem 1.7 applies. However, it also acts on
a tree so it is not Kazhdan already by a result of [Alperin 1], [Watatani 1]. Generally,
subgroups of $\mathcal{D}iff^{1,\alpha}(S^1)$ like described above, do not
have any obvious action on a tree and one needs our Theorem 1.7 to show that they are
not Kazhdan. A parallel theorem for symplectomorphism groups will be given
in Chapter IV. Notice also that the proof that our characteristic classes constructed in Section 5 are in polynomial cohomology agrees with a recent result on the growth of the Dehn function of the Thompson group [Guba 1].

\chapter{Geometry of unitary cocycles}

In this Chapter we return to the asymptotic geometry of finitely generated
groups. If $G$ is not Kazhdan, then an orthogonal cocycle $l\in
Z^1(G,\mathcal{H})$ should be viewed as a way to linearize the geometry of
$G$. Our first result is a convexity theorem 2.1 which says that the
embedding of $G$ into the Hilbert space $\mathcal{H}$ given by $l$ coarsely
respects the geometry in a sense that inner points of big "domains" in $G$
are mapped inside the convex hull of the image of boundary points.

We have seen in Chapter I that primitive functions $\mathcal{F}: G\to
\mathbb{R} $ of cocycles in $Z^1(G,l^p(G))$ of a surface group satisfy
$$ |\mathcal{F}(g)|< c\cdot length (g)^{1/p'} $$
Here, we start a general study of cocycle growth. We show in Theorem 3.1 that for any
orthogonal cocycle $l: G\to \mathcal{H} $,
$$ \| l(g)\|< c(\theta) [length (g)\log\log\ length (g) ]^{1/2} $$
for almost all $\theta\in \partial G\simeq S^1 $ and $g\to \theta $
nontangentially. We use in proof an adjusted version of Makarov's law of
iterated logarithm. The result extends to all complex hyperbolic cocompact
lattices of any dimension.

Using another deep result of Makarov, we show the following in Proposition  3.3. Let $G$ be a
surface group, $\beta: G\to \mathbb{Z} $ a surjective homomorphism and
$G_0=Ker\beta $. Then the conical limit set of $G_0$ has Hausdorff dimension
1, in particular, the exponent $\delta (G_0)=1 $. We do not know if this set has a full Lebesgue measure ( it is certainly a doable problem).
 
Notice that the proof of Lemma I.11.8 shows that the estimate on $\| l(g)\|$
is essentially sharp. It also shows that this estimate does not hold in
other Banach spaces. However, imposing various restricitons on a Banach
space, one still hopes to get an estimate, reflecting a fine structure of
$G$.

\section{Smooth and combinatorial harmonic sections}

Let $G$ be a finitely generated group. $\pi: G\to O(\mathcal{H})$ an
orthogonal representation without almost invariant vectors and $l:
G\to\mathcal{H}$ a nontrivial cocycle. If $M$ is a compact Riemannian
manifold with $\pi_1(M)=G$ (so that $G$ is finitely presented) then one
forms a flat affine bundle $E$ over $M$ with fiber $\mathcal{H}$ and
monodromy
$$ g\mapsto (v\mapsto \pi(g)v+l(g) )$$

A result of [Korevaar-Schoen 1] and [Jost 1] (Lemma I.11.6) states that
there is a harmonic section $f$ of $E$. If $M$ is K\"ahler then there is
another cocycle $m: G\to \mathcal{H} $ so that a complex affine bundle with
fiber $E\otimes \mathbb{C} $ and monodromy
$$ g\mapsto (v+iw \mapsto \pi(g)v+i\pi(g)w +l(g)+im(g)) $$
admits a holomorphic section. Our first result is a combinatorial version of
this theorem.

Let $\{ \gamma_i \} $ be a finite set of generators for $G$. Let $V$ be a
space of "sections", that is, $G$-equivariant maps
$$ f: G\to \mathcal{H} .$$
This simply means that $f(g^{-1}x)=\pi(g) f(x)+l(g) $. Obviously, such map
is determined by $f(1)\in \mathcal{H} $. Therefore, $V\approx \mathcal{H} $. A
combinatorial Laplacian is defined as
$$ \triangle f(x)=\sum_i f(\gamma_i x)+f(\gamma_i^{-1} x)-2f(x) .$$

\noindent
{\it Proposition 2.1.}--- There exists an equivariant $f: G\to
\mathcal{H} $ with $\triangle f=0 $.

\noindent\\
{\it Proof.}---
Let $v=f(1)$, then $f(x^{-1})=xv+l(x)$. Therefore
\renewcommand{\baselinestretch}{1.5}
\large
$$ \begin{array}{rcl}
    \triangle f(x^{-1})&=& \sum f(\gamma_i x^{-1})+f(\gamma_i^{-1}x^{-1})-2f(x^{-1} ) \\
    &=&\sum x\gamma_i^{-1}v+l(x\gamma_i^{-1})+x\gamma_iv+l(x\gamma_i)-2xv-2l(x) \\
    &=&\sum x(\gamma_i^{-1}v+\gamma_iv-2v)+\sum 
		xl(\gamma_i^{-1})+l(x)+xl(\gamma_i)+l(x)-2l(x)\\
    &=&x\sum (\gamma_i^{-1}+\gamma_i-2)v+x\sum [l(\gamma_i^{-1})+l(\gamma_i)],
    \end{array}
 $$
\renewcommand{\baselinestretch}{1.2}
\normalsize
 so that we need only to solve an equation
 $$ \sum (\gamma_i^{-1}+\gamma_i-2)v=-\sum [l(\gamma_i^{-1})+l(\gamma_i)] .$$
 Notice that $\widetilde \triangle : \mathcal{H}\to \mathcal{H} $ defined by
$v\mapsto \sum (\gamma_i^{-1}+\gamma_i-2)v $ is selfadjoint. 
Moreover, since $\widetilde
\triangle =-\sum (\pi (\gamma_i)-1)^* (\pi(\gamma_i)-1) $, 
$\widetilde \triangle $ is
nonpositive and if $0\in spec (\widetilde \triangle )$, then $\pi: G\to
O(\mathcal{H}) $ has almost invariant vectors. Therefore, 
$\widetilde \triangle
$ is invertible and the result follows.

\section{A convexity theorem}

We keep the notation of 3.1. Any cocycle $l: G\to \mathcal{H} $ can be seen
as an embedding of $G$ in the Hilbert space. If $\| l(g)\|\to \infty $ as
$length(g)\to \infty $, then this embedding is uniform in the sense that $\|
l(g)-l(h) \|\to \infty $ as $\rho (g,h)\to \infty $ for any word
left-invariant metric on $G$. For instance, Proposition I.2.1 implies that any
group $G$ acting discretely (but possibly not cocompactly) on an Hadamard
manifold of pinched negative curvature, admits a uniform embedding into
$l^p(G),\ p>1 $. We are, however, interested in a finer geometry of the cocycle 
embeddings.

For a finite $A\subset G$ and $C>0$, a $C$-interior $int_C(A) $ is defined
as $\{ x| \rho(x,y)<C \Rightarrow y\in A\} $. A $C$-boundary $\partial_C
(A)$ is defined as $A\backslash int_C(A)$.

\noindent\\
{\it Theorem 2.2.}--- Let $\pi: G\to \mathcal{O}(\mathcal{H})$ be an
orthogonal representation without almost-invariant vectors. Let $l: G\to 
\mathcal{H} $ be a cocycle for $\pi$. Then there are constants 
$C_1,C_2(l)>0$ such that for any finite $A\subset G$ and any $x\in A$,
$$ dist_{\mathcal{H}}( l(x)-\overline{conv} (l(\partial_{C_1}A))) \le C_2.
\eqno(*) $$

\noindent
{\it Proof.}--- Let $f:G\to \mathcal{H}$ be an equivariant harmonic
map of Proposition 1.1. Since $\| f(x^{-1})-l(x)\|=\| f(1)\| =const $, we
can replace $(\ast)$ by a condition
$$ dist_{\mathcal{H}}(f(x)-\overline{conv} f(\partial_{C_1}(A) )\le C_2' ,$$
where however, one uses a right-invariant word metric on $G$ in definition
of $\partial_C(A)$. This result follows from the maximum principle of
harmonic functions. Indeed, let $x\in int_{C_1}(A)$ be such that
$dist_{\mathcal{H}}(f(x)-\overline{conv} f(\partial_{C_1}(A) )) $ is maximal
possible (and $>C_2$)( a choice of $C_1,C_2$ will be made later). Let $v$ be
a unit vector, such that
$$ (f(x)-y,v)=dist_{\mathcal{H}}(f(x)-\overline{conv}f(\partial_{C_1}(A)) $$
for some $y\in \overline{conv}f(\partial_{C_1}(A))$. Let $h(z)=(f(z)-y,v)$.
Then $h(x)>C_2$ and $h(\partial_{C_1}(A))\subset (-\infty,0] $. Moreover,
$\widetilde\triangle h=0 $ and $h(z)\le h(x) $ for $z\in int_{C_1}(A) $. It
follows that 
$h(\gamma_i x)=h(x)$ for all $i$. 
Replacing $x$ by $\gamma_i x$ and
continuing until we hit $\partial_{C_1}A$, we arrive to a contradiction 
with $C_1=2$, $C_2=2 \| f(1) \|+1 $.

\section{Cocycle growth for a surface group}

In this section we continue, for general representations, a subject started 
in I.5.2. Recall that, for any group $G$, any primitive function $\mathcal{F}: G\to \mathbb{R} $ 
of a class in $H^1(G,l^p(G))$ satisties
$$ | \mathcal{F}(g)|\le const\cdot length (g) $$
at least of $G$ is finitely presented. However, if $G=\pi_1(\Sigma) $, a 
surface group, then one has much finer estimate, established in Theorem I.5.2:
$$ | \mathcal{F}(g) |\le const\cdot length (g)^{1/p'} .$$

\noindent\\
{\it Theorem 3.1.}--- Let $G=\pi_1(\Sigma) $ be a surface group. Let 
$\pi : G\to O(\mathcal{H})$ be an orthogonal representation 
without almost-invariant 
vectors and let $l:G\to \mathcal{H} $ be a cocycle. Then for 
almost all $\theta\in S^1\approx \partial G$,
$$ \| l(g)\| \le const(\theta)[length(g)\cdot\log\log\ length(g)]^{1/2} 
\eqno(*) $$
as $g$ converges nontangentially to $\theta$. Here "almost all" corresponds 
to a Lebesgue measure on $\partial G$, identified with $S^1$ under some 
lattice embedding $G\hookrightarrow SO^{+}(1,2)$.

\noindent\\
{\it Remark.}--- Nontangential convergence of points of $B^2$ to 
$\theta\in \partial B^2$ is an invariant of quasi-conformal homeomorphism 
[???]. Therefore $(\ast)$ is $\mathcal{M}ap_{g,1}$-invariant. Let $A\subset 
S^1$ be an exceptional set where $(\ast)$ does not hold. It follows that the 
Lebesgue measure:
$$ meas\ \varphi (A)=0 $$
for all $\varphi\in \mathcal{M}ap_{g,1}$, considered as a quasisymmetric 
homeomorphism of $S^1$.

\noindent\\
{\it Proof.}--- Complexifying, we find a holomorphic section of an 
affine bundle $E_{\mathbb{C}}$ as in section 1. Lifting to $\mathcal{H}^2$,
we obtain an equivariant holomorphic map (we replace $\mathcal{H}$ by
$\mathcal{H} \otimes \mathbb{C} $)
$$ \widetilde f: \mathcal{H}^2\to \mathcal{H} .$$
Notice that $\widetilde f$ is a Bloch function, that is, $\| \nabla \widetilde
f\|\le const $. The result now follows from a version of the Makarov law
[Makarov 1] of iterated logarithm for Hilbert-space-valued Bloch functions.

\noindent\\
{\it Proposition.}--- Let $\psi: B^2\to \mathcal{H} $ be holomorphic
and $\| \nabla \psi\|_h \le const $. Then for almost all $\theta\in S^1$,
$$ \limsup_{z\to \theta} \frac{\| \psi(z) \|}{ \sqrt{\log
(1-|z|)\log\log\log(1-|z|) } } <\infty .$$
{\it Proof.}---
 We will simply note which changes should be made in a proof for
complex-valued functions [Pommerenke 1]. The Hardy identity [Pommerenke 1,
page 174] holds in the following form. Let $S$ be a Riemannian surface,
$z_0\in S $, $g:S\to \mathcal{H} $ a holomorphic function, $(x,y)$ normal
coordinates in the neighbourhood of $z_0$. Let $n$ be a positive integer.
Then
$$
\frac{\partial}{\partial x}(g,g)^{n+1} =(n+1)(g,g)[(g'_x,g)+(g,g'_x)]
,$$
\renewcommand{\baselinestretch}{1.5}
\large
$$ \begin{array}{rcl}
	\frac{\partial^2}{\partial x^2}(g,g)^{n+1}=n(n+1)(g,g)^{n-1}[g_x',g)+(g,g_x')]^2+&& \\
	    
		(n+1)(g,g)^n[2(g_x',g_x')+(g_x'',g)+(g,g_x'')]
	\end{array}
$$
\renewcommand{\baselinestretch}{1.2}
\normalsize
and the same for $\frac{\partial^2}{\partial y^2} $. Summing up, we have
\renewcommand{\baselinestretch}{1.5}
\large
$$ \begin{array}{rcl}
	\triangle (g,g)^{n+1}&=&(\frac{\partial^2}{\partial
		x^2}+\frac{\partial^2}{\partial y^2})(g,g)^{n+1} \\
    &=&n(n+1)(g,g)^{n-1}\cdot 4|(g',g)|^2+(n+1)(g,g)^n\cdot 2(g',g') ,
    \end{array}
$$
\renewcommand{\baselinestretch}{1.2}
\normalsize
because $\triangle g=0 $ and $g_y'=\sqrt{-1}g_x'$. If $S$ is a unit disc
then in polar coordinates $z=re^{it}$
$$
	\triangle= \frac{\partial^2}{\partial
	r^2}+\frac{1}{r}\frac{\partial}{\partial
	r}+\frac{1}{r^2}\frac{\partial^2}{\partial
	t^2}
    =\frac{1}{r}\frac{\partial}{\partial r}(r\frac{\partial}{\partial
	r})+\frac{1}{r^2}\frac{\partial^2}{\partial t^2}
    =\frac{1}{r^2}[(r\frac{\partial}{\partial
	r})^2+\frac{\partial^2}{\partial t^2} ].
 $$
 So
 $$ 
 	\frac{1}{r^2}((r\frac{\partial}{\partial 
	r})^2+\frac{\partial^2}{\partial 
	t^2})(g,g)^{n+1}=4n(n+1)(g,g)^{n-1}|(g',g)|^2+2(n+1)(g,g)^n|g'|^2.
 $$
 Integrating over $0\le t\le 2\pi $ and using Cauchy-Schwartz inequality, we
arrive at the inequality of [Pommerenke 1, Theorem 8.9]. The rest of the
proof will go unchanged once we know the Hardy-Littlewood maximal theorem
for $(g,g)^n$, which is used in [Pommerenke 1, page 187]. Let
$$
    g^*(s,\xi)=\max_{0\le r\le 1-e^{-s}} |g(r\xi)|, e\le s<\infty, \xi\in 
S^1 .$$
Since $g:B^2\to\mathcal{H} $ is holomorphic, it is also harmonic and yields 
the Poisson formula. Then a proof of the Hardy-Littlewood maximal theorem 
given in [Koosis 1] applies, since it reduces it to the Hardy-Littlewood 
inequality for the maximal function of $|g|$.

\noindent\\
{\it Remark 3.2.}--- Theorem 3.1 holds for complex hyperbolic 
cocompact lattices. This is because Makarov's law of iterated logarithm 
holds for the complex hyperbolic space, as we can see by passing to totally 
geodesic spaces of complex dimension 1. It is plausible that a version of 
Theorem 3.1 holds for real hyperbolic lattices (but not quaternionic and 
Cayley, as these are Kazhdan, see a new proof in Chapter VI). On the other hand, another deep result of [Makarov 2] saying that Bloch functions are nontangentially bounded for a limit set of Hausdorff dimension one, fails for Hilbert space valued functions. In fact, we have shown in Chapter I  that there are unitary cocycles on a surface group  such that $ \| l(g)\| \to \infty$ as $length(g)\to \infty$.

If $G$ is any finitely generated group, and we are given an orthogonal 
representation $\pi: G\to O(\mathcal{H}) $ and a cocycle $l\in 
Z^1(G,\mathcal{H}) $ with a control on $\| l(g)\| $ from below, then for any 
embedding of the surface group $\pi_1(\Sigma) $ into $G$ we immediately have 
a comparison inequality between the word lengthes of elements of 
$\pi_1(\Sigma)$ in $\pi_1(\Sigma)$ and $G$. To get a nontrivial result, we 
need a low bound on $\| l(g)\| $ better then $[length(g)\log\log\ 
length(g)]^{1/2} $. To find such groups and cocycles seems to be a very 
attractive problem.

We will now use similar ideas to estimate the Hausdorff dimension of limit 
sets of some infinite index subgroups of $G=\pi_1(\Sigma)$.

\noindent\\
{\it Theorem 3.2.}--- Let $\beta: G\to\mathbb{Z}$ be a surjective
homomorphism and let $G_0=Ker\beta$. Let $A$ be a conical limit set of
$G_0$. Then $\dim A=1 $. 

\noindent\\
{\it Proof.}--- Let $[\beta ]\in H^1(G,\mathbb{Z})$ be an induced
class. Realize $G$ as a cocompact lattice in $SO^{+}(1,2)$ so that
$S=\mathcal{H}^2/G$ is a hyperbolic surface. Let $\omega$ be a holomorphic
1-form on $S$ such that $Re[\omega]=[\beta]$ and let $\widetilde \omega $ be a
lift of $\omega$ on $\mathcal{H}^2$. Let $f:\mathcal{H}^2\to \mathbb{C}$ be
holomorphic with $df=\omega $. Then $f$ is a Bloch function. By a result of
[Makarov 2], there is a set $B\subset S^1$ with $\dim B=1$ such that
$$ \limsup_{z\to \theta} |f(z)|<\infty $$
for any $\theta\in B$ and nontangential convergence. Notice that
$f(gz)=f(z)+([\omega],[g])$ where $g\in G$ and $[g]$ is an image of $g$ in
$H_1(G,\mathbb{Z})$. Now it is clear that $B\subseteq A$, so $\dim A=1$.

\noindent\\
{\it Remark.}--- This result does not contradict a theorem of
[Sullivan 1] and [Tukia 1] because $G_0$ is infinitely generated.

\noindent\\
In the opposite direction we have the following . Let $\Sigma_1, \Sigma_2$ be two closed surfaces
and let $\psi: \Sigma_1\to \Sigma_2$ be a smooth ramified covering. Let
$G_i=\pi_1(\Sigma_i)$ and let $G_0 =Ker\psi_* : G_1\to G_2 $. Let
$G_1\hookrightarrow SO^{+}(1,2) $ be a realization of $G_1$ as a lattice. Then
for any $z\in B^2$,
$$ \sum_{g\in G_0}| 1-gz|<\infty .$$
In other words, either $\delta (G_0) <1 $ or $\delta (G_0)=1$ and $G_0$ is
of convergence type. In the latter case, the Patterson measure of the
conical limit set of $G_0$ is zero. To see this, notice that we can find hyperbolic structures on $\Sigma_i,\ i=1,2 
$ so that $\psi$ is holomorphic. Let $\widetilde \psi $ be a lift of $\psi $ as 
a map $\widetilde\psi: B^2\to B^2 $. Since $\widetilde \psi$ is a bounded 
holomorphic function, $\widetilde \psi$ has limit values almost everywhere on 
$S^1$. By I.6.5., $| \widetilde \psi |S^1| =1 $ almost everywhere. So $\widetilde 
\psi$ is an inner function. Let $C\subset B^2$ be a countable set of zeros 
of $\widetilde\psi$. We claim that $C$ is a finite union of orbits of $G_0$. 
First, it is clear that $C$ is $G_0$-invariant. Let $Q\subset B^2 $ be compact
which contains a fundamental domain for $G_1$. Then $\widetilde\psi(Q)$ is
compact so there is a finite set $R\subset G_2$ such that $g(0)\notin
\widetilde\psi(Q)$ if $g\notin R$. Let $T\subset G_1$ be finite and such that
$\psi_*(T)\supseteq R $. Let $Q_1=\bigcup_{g\in T^{-1}}gQ $ so that $Q_1$ is
compact and therefore $C\cap Q $ is finite. Let $x\in C$, then $x=gy$ with
$y\in Q$. So $0=\widetilde \psi (x) =\psi_*(g)\widetilde \psi(y) $, i.e.,
$\psi_*(g^{-1})(0)=\widetilde\psi(y)\in \widetilde \psi(Q) $. This means
$\psi_*(g^{-1})\in R$ so $g^{-1}\in TG_0$, and $g\in G_0T^{-1} $, say
$g=g_0t^{-1}$, $g_0\in G_0$, $t\in T$. Then $t^{-1}y\in C$ and $t^{-1}y\in
Q_1$, so there are finitely many options for $t^{-1}y$.

We deduce that there are $x_1,\cdots,x_n$ such that $C=\bigcup_{i=1}^n G_0
x_i $. The decomposition formula for inner functions implies that
$$ \widetilde\psi(z)=c\cdot \prod_{\stackrel{g_0\in G_0}{1\le i\le n}
}\frac{\overline{g_0x_i}}{g_0x_i} \frac{z-g_0x_i}{1-\overline{g_0x_i}z} $$
which gives an explicit formula for holomorphic maps between hyperbolic
Riemann surfaces (one still needs to find $x_i$). By a well-known result
on zeros of a bouded holomorphic function [Koosis 1, IV: B, Theorem 1],
$$ \sum_{g_0\in G_0}(1-|g_0x_i| )< \infty .$$
The rest follows from [Nicholls 1].

\chapter{A theory of groups of symplectomorphisms}

We already have noticed an intriguing similarity between groups acting on the 
circle and groups acting symplectically on a compact sympletic manifold.
The two leading topics studied in Chapter 2, namely, (non-) Kazhdan groups 
acting on $S^1$ and characteristic classes, have exact analogues 
for {\it Sympl(M)}.
In fact, a theory of characteristic classes parallel to II.5 , has already been
presented in [Reznikov 2] and  [Reznikov 4]. In the second cited paper, we 
noticed
that the K\"{a}hler action of 
{\it Sympl(M)} on the twistor variety allows us to 
define a  series of  classes in $H^{2k}_{cont}(Sympl(M),\mathbb{R})$, $k$ 
odd , which are
highly non-trivial. In the first cited paper,we introduced bi-invariant 
forms on 
{\it Sympl(M)} and the classes in 
$H^{odd}_{top}(Sympl(M)) $ and
$H^{odd}$(${Sympl(M)}^\delta,\mathbb{R}/A)$ (cohomology of a topological space 
and a discrete group) where $A$ is a group of periods of the above-mentioned 
forms.
Here we present a fundamental class in $H^{1}$$({Sympl(M)},L^2(M))$ whose 
nontriviality on a subgroup $G\subset{Sympl(M)}$ implies that $G$ is not
Kazhdan, similarly to the situation in  ${\mathcal{D}}iff^{1,\alpha}(S^1)$. 
>From the nature of our class it is clear that its vanishing imposes severe 
restriction
on the symplectic action, roughly, the transformations of $G$ should satisfy 
a certain PDE . We give an explicit formula for our class in the case of a 
flat torus.

We then introduce a characteristic class in $H^{n+1}({Sympl}^\delta (M^{2n}),
\mathbb{R})$ associated with a compact Lagrangian immersed submanifold. This
class is a sympletic counterpart, and a generalization, of the Thurston-Bott class 
[Bott 1]. We use this class to give a formula for the volume of compact 
negatively curved manifold through Euclidean volumes of ``Busemann bodies" 
(the images of the manifold under Busemann functions).

\section{Deformation quantization : an overview}

Let $F$ be a field and $A|F$ a (commutative) algebra . A deformation of A is 
an algebra structure over $F[[\hbar]]$ of $A[[\hbar]]$
extending that $A$, so that if $x,y\in A$,
$$
x\ast y=x\cdot y+{b}_1(x,y)\hbar +{b}_2(x,y) \hbar^2 +\cdots
$$
where $x\cdot y$ is a multiplication in $A$ and $ x\ast y$ is a deformed 
multiplication .

If $F=\mathbb{R}$, $A=C^{\infty}(M)$, where $M$ is a symplectic manifold , 
then a deformation quantization is a deformation of $A$ with
${b}_1(f,g)=\{f,g\}$, a Poisson bracket . A deformation quantization 
always exists by a result of [Moyal 1], [Vey 1],
[Bayen-Flato-Fronsdal-Lichnerowicz-Sternheimer 1]
[Fedosov 1].
For any algebra $A|F$ one defines a Hochschild collomology ${HH}^k(A)=
{Ext}^{k}_{A\otimes A}(A,A) $. There is a natural Lie superalgebra 
structure in
${HH}^{\ast}(A) $ [Gerstenhaber 1]. There exists a simple   explicit 
complex , computing ${HH}^{k}(A) $ with ${C}^k(A) = {Hom}_{F}(
{\bigotimes}_{i=1}^{k} A,A) $. In particular , $b_1$ above is a cocycle 
(for any deformation ). If $F=\mathbb{R}$ and A is a topological algebra, one
modifies the definitions to obtain topological Hochschild cohomology . 
If $M$ is a smooth manifold and $A=C^{\infty}(M) $ with a pointwise
multiplication ,
then
$$
{HH}^k(A)=\Gamma(M,{\Lambda}^k T M  ),
$$
a space of poly-vector fields . The Lie superalgebra structure coincides 
with a   classical bracket of poly-vector fields.

We will need an explicit form of the cocycle condition for a 2-cocycle 
$b: A\otimes A \rightarrow \mathbb{R} $ :
$$
xb(y,z)-b(xy,z)+b(x,yz)-b(x,y)z=0
.$$

\section{A fundamental cocycle in $ H^1(Sympl(M), L^2(M))$}

Let $(M^{2n},\omega) $ be a compact symplectic manifold . Fix a deformation
quantization
$$
f\ast g=f\cdot g+\{f,g\}\hbar +
	{\sum_{i=2}^{\infty}c_i(f,g)\cdot\hbar^i}
$$

Let $\Phi: M\rightarrow M $ be symplectic and let
\renewcommand{\baselinestretch}{1.5}
\large
$$
\begin{array}{rcl}
f \tilde{\ast} g &= & (f\circ {\Phi}^{-1}\ast g\circ{\Phi}^{-1})\circ {\Phi}\\
&= & f\cdot g +\{f,g\}\hbar +
	\sum_{i=2}^{\infty}c'_i(f,g)\cdot\hbar^i
\end{array}
$$
\renewcommand{\baselinestretch}{1.2}
\normalsize
{\it Lemma 2.1.}--- Let $A|F $ be an algebra and let
$$
f\ast g=f\cdot g+c_1(f,g)\hbar + \cdots +c_{k-1}(f,g)\hbar^{k-1}+
	{\sum_{i=k}^{\infty}c_i(f,g)\cdot\hbar^i}
$$
and
$$
f \tilde{\ast} g =  f\cdot g +c_1(f,g)\hbar + \cdots +c_{k-1}(f,g)\hbar^{k-1}+
	\sum_{i=k}^{\infty}c'_i(f,g)\cdot\hbar^i
$$
be two deformations, which coincide up to the order $\hbar^{k-1}$. Then
$$
c_i - c'_i : A\otimes A \rightarrow A
$$
 is a Hochschild cocycle .

\noindent\\
{\it Proof.}---
$$
(f\ast g) \ast p - (f \tilde{\ast} g )\tilde{\ast} p = c_k(f,g)\cdot
p+c_k(f\cdot g,p)-c'_k(f,g)\cdot p-c'_k(f\cdot g,p)\qquad (\mbox{\it mod }\hbar^{k+1})
$$
Similarly,
$$
f\ast {(g\ast p)}  - f\tilde{\ast} {(g \tilde{\ast} p )} = f\cdot c_k(g,p)
+c_k(f,g\cdot p)-fc'_k(g,p)-c'_k(f,g\cdot p)\qquad (\mbox{\it mod } \hbar
^{k+1})
$$
So for $ c=c_k - c'_k $,
$$
f\cdot c(g,p)+c(f,g\cdot p)-c(f,g)p-c(f\cdot g ,p)=0
,$$
which means that $ c$ is a 2-cocycle.

\noindent\\
{\it Lemma 2.2.}---  A formula
$$
\Phi \mapsto \left[ (f,g) \mapsto c_2(f\circ{\Phi}^{-1},g\circ{\Phi}^{-1} )
\circ{\Phi}-c_2(f,g)\right]
$$
defines a smooth cocycle of {\it Sympl(M)} in the space $ Z^{2} (C^{\infty}(M)$,
$C^{\infty}(M)) $ of Hochschild 2-cocycles for $C^{\infty}(M) $.

\noindent\\
{\it Proof.} Follows from Lemma 2.1.\\

Passing to Hochschild cohomology, we obtain a 1-cocycle of $Sympl(M)$ in\\
 $${HH}^2(C^{\infty}(M))=\Gamma (M,{\Lambda}^2 TM).$$
Using the symplectic structure ,we identify $\Gamma (M,{\Lambda}^2 TM)$ with ${\Omega}^2(M) $, a space of 2-forms on $M$.
Multilying by ${\omega}^{n-1}$ we obtain a cocycle
$$
\mu \in {H}^{1}(Sympl(M),C^{\infty}(M)).
$$

\section{Computation for a flat torus and the main theorem}

If $M$ is a coadjoint orbit of a compact Lie group, one can find an explicit formula for the deformation quantization $f\ast g$.
A classical case $M=T^{2n}$, a flat torus, is due to H.Weyl.

\noindent\\
{\it Proposition 3.1.}---
One has a following deformatiom quantization on $T^{2n}$ :
$$
f\ast g =\displaystyle \sum_{k=0}^{\infty} \frac{1}{k!}{\left(-\frac{i\hbar}{2} \right)}^k{\sigma}^{i_1 j_1}\cdots {\sigma}^{i_k j_k}
                {\frac{{\partial}^k f}{{{\partial {y_{i_1}}}}\cdots {\partial {y_{i_k}}}}} 
		{\frac{{\partial}^k g}{{{\partial {y_{j_1}}}}\cdots {\partial {y_{j_k}}}}}
$$
where ${\sigma}^{ij} $ are elements of the matrix, inverse to the matrix $({\sigma}_{ij}) $ of a (constant) symplectic form , and the ``repeated
indices" summation agreement is applied .

Now, since our definition of a fundamental cocycle is completely explicit, 
one can derive an explicit formula for $\mu $ in this case.
We give an answer for $ T^2 $ (the formula for $T^{2n} $ is completely analogous). The computation is tedious (takes several pages) but straightforward and is left to
reader . Here is the formula for $T^2$ :
$$
\Phi \mapsto {\frac{{\partial}^2{\Phi}_2}{\partial {y}_1^2}} 
	{\frac{{\partial}^2{\Phi}_1}{\partial {y}_2^2}} - 
	{\frac{{\partial}^2{\Phi}_1}{\partial {y}_1^2}}
                               {\frac{{\partial}^2{\Phi}_2}{\partial {y}_2^2}}
$$
where $\Phi = ({\Phi}_1 ,{\Phi}_2 ) $ a symplectomorphism of the form $T^2$. Summing up, we have

\noindent\\
{\it Theorem 3.2.}--- Let $M^{2n}$ be a compact symplectic manifold , let $ Sympl(M) $ its symplectomorphism group , acting orthogonally on
a Hilbert space $L^2(M) $. There exists a cocycle
$$
\mu \in Z^1(Sympl(M),L^2(M)),
$$
defined canonically by a given deformation quantization of $C^{\infty}(M) $ with the following properties :\\
A. Let
$$
\begin{array}{rcl}
f \tilde{\ast} g &=&  (f\circ {\Phi}^{-1}\ast g\circ{\Phi}^{-1})\circ {\Phi} =  f\cdot g +\{f,g\}\hbar 	+c'_2(f,g)\cdot\hbar^2+\cdots ,\\
f\ast g &= & f\cdot g+\{f,g\}\hbar +{c_2}(f,g)\cdot\hbar^2+\cdots
\end{array}
$$
and let us indentify the class of the Hochschild 
cocycle $c'_2-{c_2} $ with a section $\nu$ of ${\Lambda}^2 TM $.
Let $\hat{\nu}$ be a 2-form obtained from $\nu$ by lifting the indices using the symplectic form . Then
$$
\mu(\Phi) \cdot {\omega}^n = \hat{\nu}. {\omega}^{n-1}
.$$
B. 
$\mu(\Phi) $ depends only on the second jet of $\Phi$.\\
C. 
For $M=T^2 $ and the Weyl deformation quantization , $\Phi =({\Phi}_1,{\Phi}_2)$,
$$
\mu(\Phi) = {\frac{{\partial}^2{\Phi}_2}{\partial {y}_1^2}} 
	{\frac{{\partial}^2{\Phi}_1}{\partial {y}_2^2}} - 
	{\frac{{\partial}^2{\Phi}_1}{\partial {y}_1^2}}
                               {\frac{{\partial}^2{\Phi}_2}{\partial {y}_2^2}}
.$$
D. 
If $G$ is a Kazhdan subgroup of $Sympl(M)$, then
$$
{\|{\mu(\Phi)}\|}_{L^2}<const  \qquad (\Phi \in G).
$$
{\it Examples.}---\\
1) 
$M=T^{2n}$, $G=Sp(2n,Z)$   \quad (Kazhdan for $n\geq 2$).  Then $\mu $ is identically zero.
\noindent\\
2)Let $ \Gamma$ be a surface group, and let $M$ be a component of
$$Hom(\Gamma,SO(3))/SO(3),$$
consisting of representations with nontrivial Stiefel-Whitney class. 
Then $M$ is a compact symplectic manifold and ${Map}_g$ acts
symplectically on $M$. 
We do not know if part D of Theorem 3.2 holds in this case and if ${Map}_g$ is Kazhdan or not. There
is a ``Teichm\"uller structure" on $M$ defined by a holomorphic map of the 
Teichm\"uller space into the twistor variety of $M$,
described in [Reznikov 4], see also Chapter 5.

\noindent\\
{\it Remark.}---
The case of two-dimensional $M^2$ is much easier,
simply because ${SL}_2(\mathbb R)$ is not Kazhdan . If $Sympl(M,x_0)$
is a subgroup fixing $x_0\in M$ then one gets a nontrivial unitary cocycle 
on $Sympl(M,x_0)$ by pulling back from ${SL}_2(\mathbb R)$
under the tangent representation. Using the 
measurable transfer (= Shapiro's lemma) one constructs a cocycle of $Sympl(M
)$.
See [Zimmer 1] for details .

\section{Invariant forms on the space of Lagrangian immersions and new regulators for symplectomorphism groups}

In this section we will ``symplectify" the Thurston-Bott 
class in the cohomology of diffeomorphism groups. Let $M$ be any
(possibly noncompact) symplectic manifold, 
and let $L_0\hookrightarrow M$ be a Lagrangian immersion of  a compact oriented
manifold $L_0$. Let $Lag(L_0,M)$ be a space of Lagrangian immersions
of $L$ into $M$ which can be jointed to $L_0$
by an exact Lagrangian homotopy. This means the following. 
If $f_t \rightarrow M $ is a smooth family of Lagrangian immersions,
than $\displaystyle{\frac{d}{dt} {f_t |}_{t=0}} $ is a vector field along $L_0$.
Projecting to the normal bundle $NL_0$ and
accounting that $NL_0$ is canonically isomorphic to $T^{\ast}L_0$, we get a 1-form on $L_0$ which is immediately seen to be
 closed. A Lagrangian homotopy $f_t$ is exact, if this form is exact for all t. 
There is therefore a well-defined function $F$
 (mod const ) on $L$ which can be seen as a tangent vector of the deformation .

\noindent\\
{\it Definition.}---
A canonical (n+1)-form $\nu$ on $Lag(L,M)$ is defined by
$$
\nu(F_0 \cdots F_n)=\int_L F_0 d{F_1 \cdots F_n} = {Vol}_{n+1} (\widetilde {Q})\qquad (\ast),
$$
where $\widetilde {Q} $ is any chain in ${\mathbb R}^{n+1} $ spanning 
$ (F_0,\cdots ,F_n)(L)$. 

\noindent\\
{\it Proposition 4.1}---$\nu$ is closed. 

\noindent\\
{\it Proof} is an exercise for reader.
\\
 
Let ${Sympl}_0 (M) $ be a group of Poissonian transformations of $M$. 
Then $Lag(L,M) $ is invariant under $Sympl_0 (M)$.

\noindent\\
{\it Proposition 4.2.}---
$\nu $ is $Sympl_0 (M)$-invariant.

\noindent\\
{\it Proof} is obvious.
\\

A standard theory of regulators [Reznikov 3], [Reznikov 2] implies that, 
first, one has an induced class in $H^{n+1}(\mathfrak{g},\mathbb R)$.
where $ \mathfrak{g}$ = $Lie (Sympl_0(M))$= 
$C^{\infty} (M)/const $ given by $(\ast),$ where now $F_i \in C^{\infty}(M) $ and second,
a class in
$$
Hom(\pi_{n+1}({BSympl_0^\delta (M)}^{+},\mathbb{R}/A))  \qquad (n+1 \geq 5),
$$
where $A$ is a group of periods of $\nu $ on maps 
$\Sigma^{n+1} \rightarrow Sympl_0(M) $ of homology spheres to $Sympl_0(M) $.
This class often lifts to a class in $H^{n+1}(Sympl_0^\delta(M),\mathbb R)$ 
under suitable conditions on topology of  $Sympl_0(M) $
(see the discussion in the papers cited above ).

As an example, let $Q$ be a compact oriented simply connected manifold, 
$ M=T^{\ast}Q $ and $L_0=Q$, a zero section.
Then we obtain a class $[\nu]$ in $H^{n+1}(Sympl_0(T^{\ast}Q),\mathbb R)$. 
Notice that the restriction of this class on
$\mathcal{D}iff(Q) \hookrightarrow Sympl_0(T^{\ast}Q) $ is zero, 
as $\mathcal{D}iff(Q)$ fixes the zero section. However , our class is
an extension of Thurston-Bott class [Bott 1] in 
$\mathcal{D}iff(Q)$ by means of the following construction. Let
$G\subset  Sympl_0(T^{\ast}Q) $ be a subgroup of symplectomorphisms of the form
$$
p_{x} \mapsto {\phi}^{\ast}p_x + df(x),
$$
where $f \in C^{\infty}(Q)$, $ \phi \in $ $\mathcal{D}iff(Q)$, 
$x\in Q $, $p_x \in T^{\ast}_x Q$. Then $G$ is an extension
$$
0 \rightarrow C^{\infty}(Q)/const \rightarrow G \rightarrow {\mathcal D}iff(Q) \rightarrow 1.
$$
Any 1-cocycle $ \psi \in Z^1({\mathcal D}  iff(Q), C^{\infty}(Q)/const )$ induces a spliting of this exact sequence:
$$
S_{\psi} : {\mathcal D}  iff(Q) \rightarrow G
.$$
Now let $\mu$ be a smooth density on $Q$ then $\psi = \frac{{\phi}_{\ast}\mu}{\mu}$ 
is a 1-cocycle, so it defines such a splitting.
A pull-back $ S_{\psi}^{\ast}([\nu]|G) $ of our class on 
${\mathcal D}iff(Q)$ is precisely the Thurston-Bott class.

We sum up :

\noindent\\
{\it Theorem 4.3.}---\\
A.  A formula
$$
\nu (F_0 , \cdots , F_n ) = \int_L F_0dF_1\cdots dF_n = Vol_{n+1}(\widetilde Q)
$$
defines an $Sympl_0(M)$-invariant closed $(n+1)$-form in $Lag(L,M)$. It induces a class
 $[\nu]\in H^{n+1}(Lie(Sympl_0(M),\mathbb R)$ and a regulator
$$
[\nu] : \pi_{n+1}(BSympl_0^+(M)) \rightarrow \mathbb R , \qquad  n+1\geq 5,
$$
which lifts to a class
$$
[\nu] \in  H^{n+1}(Sympl_0^\delta (M),\mathbb R)
$$
if ${\widetilde H}_i(Lag(L,M),\mathbb R)) = 0 , \quad 0\leq i \leq n+1 $.\\
B. In particular , if $Q$ is a smooth oriented simply-connected closed manifold, then
$$
[\nu] \in  H^{n+1}(Sympl^\delta (T^{\ast}Q),\mathbb R)
$$
pulls back to the Thurston-Bott class under any splittting
$$
{\mathcal D} iff (Q) \rightarrow C^{\infty}(Q)/const \rtimes {\mathcal D} iff(Q)
,$$
coming from a smooth density on $Q$.

\section{A volume formula for negatively-curved manifolds}

This  section is ideologically influenced by [Hamenst\"adt 1] and discussions with G.Besson (Grenoble, 1996). 
Let $N^n $ be an Hadamard manifold. Let  $CN$
be the space of oriented geodesic of $N$, which is a symplectic manifold
of dimension $2n-2$. Any
point $x\in N$ defines a Lagrangian sphere $S_x\subset CN $ of geodesics passing through x.

\noindent\\
{\it Lemma 5.1.}---
A pull-back $ S^{\ast}\nu $ of the form $\nu \in {\Omega}^n(CN) $ to $N$ is the Riemannian volume form on $N$ times
a constant .

\noindent\\
{\it Proof.}---
An exercise in Jacobi fields.\\

Now if $G$ acts discretely and cocompactly on $N$, we have\\
$$ [S^{\ast}\nu, \mbox{fundamental class of} \ N/G]=c\cdot Vol(N/G).$$

\noindent\\
{\it Corollary 5.2.}---
$[\nu] \neq 0 $ in $ H^n{(Sympl}^\delta (N),\mathbb R) .$\\

Now we assume that the curvature of $N$ is strictly negative and moreover, 
the induced action of $G$ on the sphere at
infinity $ S_{\infty} $ is of class $C^{1,\frac{n-1}{n}} $. 
For $n=2$ this is always the case [Hurder-Katok 1], whereas for
$ n\geq 3$ seems to require a pinching of the curvature. Notice that the map
$$
s_{+} : CN\rightarrow S_{\infty},
$$
sending any geodesic $\gamma(t)$ to $\gamma(\infty)$, 
is a Lagrangian fibration. Therefore if we fix a Lagrangian
section of $s_+$, we will have a symplectomorphism 
$CN \simeq T^{\ast} (S_{\infty})$. Fix $p_0 \in N$, then
$S_{p_0}$ is such a section. Notice that an induced homomorphism 
$G \rightarrow Sympl(T^{\ast}S_{\infty})$
is given by
,$$
g \mapsto (z \mapsto \pi(g)z+dF(p_0,g^{-1}p_0,\theta)),
$$
where $g \in G,\ z\in T_{\theta}^*S_{\infty},\ 
\pi: G\rightarrow {{\mathcal D}iff}^{1,\frac{n-1}{n}}(S_{\infty})\rightarrow Sympl(T^* S_{\infty})$
is induced by the action of $G$ on $S_{\infty}$ and $ B(p_0,g^{-1}p_0,\theta))$ is the Buseman function. 
Our assumption
imply that $B(p_0,p_1, \cdot) \in C^{\frac{n-1}{n}}(S_{\infty})\subset {W_n}^{\frac{n-1}{n}}(S_{\infty})$. Recall that for
$ F_1,\cdots , F_n \in {W_n}^{\frac{n-1}{n}} $ we have an n-form
$$
\int_{S_{\infty}} F_1 dF_2 \cdots dF_n = \int_{B^n} du_1 \cdots du_n
,$$
where $u_i$ is a harmonic extension of $F_i$.

We derive a

\noindent\\
{\it Corollary 5.3.}---
Let $N^{n}/G$ be  a compact negatively curved manifold such that 
the induced action of $G$ on $S_\infty$ is of class
$C^{1,\frac{n-1}{n}} $. If the fundamental class of $G$ is
$$
\displaystyle \sum_i [{g_1}^{(i)} \cdots {g_n}^{(i)}],
$$
then the following volume formula holds :
$$
Vol(N/G)=c(n)\cdot\displaystyle \sum_i \int_{S_{\infty}}{F_1}^{(i)} d{F_2}^{(i)} \cdots d{F_n}^{(i)}
,$$
where ${F_k}^{(i)}(\theta) = B\left(p_0,{\left({g_k}^{(i)}\right)}^{-1}p_0,\theta\right) .$

One can say that a volume of a negatively curved manifold is a sum of Euclidean volumes 
of Busemann bodies in ${\mathbb R}^n$
bounded by $ (F_1,\cdots ,F_n)(S_{\infty}) $.

Replacing the Busemann cocycle by a Jacobian cocycle $ \frac{g\ast \mu}{\mu}$,
where $\mu$ is a smooth density on $S_{\infty}$,
we arrive to a similar  formula for Godbillon-Vey-Thurston-Bott 
invariant of $N/G$,under the same regularity assumptions.
This seems to have been also accomplished in a preprint [Hurder 1] cited in [Hurder-Katok 1],
though I was unable to obtain this
paper from its author. The case $n=2$ is ,however, covered in [Hurder-Katok 1].

\chapter{A theory of groups of volume-preserving diffeomorphisms 
	and the nonlinear superrigidity alternative}

In this chapter, we shift the focus from linear functional analytic techniques 
to nonlinear PDE, notably harmonic maps into nonlocally
compact spaces, a theory recently developed in [Korevaar-Schoen 1] and [Jost 1]. 
The main idea is to use twistor varietes, which
were in a center of the characteristic classes construction of [Reznikov 4], for a deeper study of volume-preserving actions of groups.
We introduce an invariant of a volume-preserving action, which we call $\Lambda$, which is a sort of a $\log L^2$-version of a $sup$-displacement
studied in [Zimmer 2]. Our first main result, Theorem 2.3, states that if $G$ is a Kazhdan group acting on a compact manifold
$M$ preserving volume, then either $\Lambda >0$ or $G$ fixes a $\log L^2$-metric. 
A much weaker analogue of this result for the special case
of lattices in Lie groups and $sup$-displacement was known before [Zimmer 2, Theorem 4.8].

We then apply our technique to  a major open problem in the field, 
that of the nonlinear superrigidity of volume-preserving actions of lattices in Lie groups. From a nonlinear version of Margulis theorem  given in [Zimmer 3] one knows that a volume preserving action of a lattice in a semisimple Lie group of rank $\ge 2$ on a low dimensional (with respect to the group) manifold fixes a {\bf measurable} Riemannian metric. Since measurable metrics do not define a geometry on a manifold, one wishes ,of course, to prove a much stronger result: that the action preserves a smooth metric. Zimmer noticed [Zimmer 2 and references therein] that such strong result would follow if one is able to find an invariant metric whose dilations with respect to any smooth metric are in the class $L^2_{loc}$. The central question of how to find such a ``bounded'' invariant metric was left completely open. We present a completely new approach to the problem which leads to 
Theorem 3.1. It  states that if a cocompact
lattice acts on $M$ preserving volume, then either it {\bf nearly} preserves a $\log L^2$-metric, 
or a sort of $G$-structure. This theorem, though constituting a clear progress in solution of the main problem is still less than what one wants in two respects: first, we deal with $\log L^2$-metrics, not $L^2$-metrics, second, we leave open a very delicate situation when an action nearly preserves a $\log L^2$-metric, but does not exactly preserve such a metric. This situation is purely infinite-dimensional (if an action on a finitely dimensional space of nonpositive curvature nearly preserves a point, it actually preserves a point at infinity). As already
said, we use a heavy machinery : harmonic maps into twistor varietes and vanishing 
results of [Mok-Siu-Yeung 1] and [Corlette 1]. These results will also be applied in the next Chapter to study quaternionic K\"ahler groups.

As is well-known, an original Kolmogorov's definition of entropy  used extremum 
over all partitions and only became  computable after
it had been realized by Kolmogorov and Sinai that certain partitions realize entropy. In a way of a pleasant similarity, we show how to compute our invariant 
$\Lambda$ for $G= \mathbb Z$ in case $G$
leaves a geodesic in the twistor space invariant, like a hyperbolic element of 
$SL(n,\mathbb Z)$ acting on $T^n$. This clearly shows an
advantage of $\log L^2$-displacement over {\it sup}-displacement.

\section{$\log L^2$-twistor spaces}

Twistor varietes $(C^{\infty})$ were used in [Reznikov 4] to define secondary characteristic classes for volume-preserving and symplectic
actions. More specifically, we have defined, for a compact oriented manifold $M$ equipped with a volume form $\nu$, a series of classes
in $H^{\ast}_{cont}({{\mathcal{D}}iff}_{\nu}(M))$ of dimension $ 5,9,13,\cdots$ (where ${{\mathcal{D}}iff}_{\nu}(M))$ is a group of  volume-preserving
diffeomorphisms). Likewise, for a compact symplectic manifold $M$ we have defined classes in $H^{\ast}_{cont}(Sympl(M))$ of dimensions
of $2,6,10,\cdots$. For purposes of the present paper, we will need to work with a $\log L^2$-version of twistor varietes, defined below.

\noindent\\
{\it Remark 1.1.}---
I would like to use an opportunity to note that by some strange reason 
I have overlooked an integrated Euler class in
$H^n_{cont}({{\mathcal{D}}iff}_{\nu}(M^n))$. The definition is exactly 
like that in [Reznikov 4] for classes in dimensions 5, 9, $\cdots$
if one realizes that there exists an $n$-form on the twistor variety for $M$,
which is ${{\mathcal{D}}iff}_{\nu}$-invariant. Alternatively,
if ${{\mathcal{D}}iff}_{\nu}(M,p_0)$ is a subgroup fixing a point $p_0$,
then one pulls back the Euler class of ${SL}_n(\mathbb R)$ under
the tangent representation
$$
{{\mathcal{D}}iff}_{\nu}(M,p_0) \rightarrow {SL}_n(\mathbb R)
,$$
and then applies a measurable transfer (see the above cited paper). 
The just defined class viewed as a class in
$H^n({{\mathcal{D}}iff}_{\nu}^s(M))$ is bounded. This follows
from the fact that the Euler class is bounded [Sullivan 2] exactly in the same manner as in [Reznikov 4].\\

We now define the $\log L^2$-twistor variety $X$ for $(M,\nu)$. First, 
one defines a bundle ${\mathcal{P}}$ of  metrics with volume form
$\nu$ as an $SL(n)/SO(n)$-bundle, associated with a principal $SL(n)$-bundle,
defined by $\nu$. Fix a smooth section (=a Riemannian
metric with volume form $\nu$) $g_0$ of this bundle. 
For any other measurable section $g$ of ${\mathcal{P}}$ define
$$
{\rho}^2(g_0,g)=\int_M {\rho}_x^2(g_0,g)d{\nu}, \eqno(*)
$$
where ${\rho}_x$ is a distance in ${\mathcal{P}}_x$ induced by (fixed once forever) 
$SL(n)$-invariant metric on $SL(n)/SO(n)$.
Now the twistor variety $X$ consists of $\log L^2$-metrics, that is,
$$
\rho(g_0,g)<\infty.
$$
Alternatively, let $A_x$ be a self-adjoint (with respect to ${(g_0)}_x$) operator such that $g_x=g_0(A_x\cdot,\cdot)$. Then $(\ast)$
can be written as
$$
\int_M {\|\log A_x\|}^2 d\nu < \infty
.$$
A crucial fact about ${\mathcal{P}}$ is a following

\noindent\\
{\it Proposition 1.2.}--- ${\mathcal{P}}$ is a complete Hilbert Riemannian manifold 
with nonpositive curvature operator.
The action of  ${{\mathcal{D}}iff}_{\nu}(M)$ on ${\mathcal{P}}$ is isometric.

\noindent\\
{\it Proof.}---
We will only define a metric, leaving all routine checks to the reader. 
A tangent space at $g_0$ consists of
$L^2$-sections of $S^2{T}^{\ast}M$, with trace identically zero. 
If $A$ is such a section (so that $A_x$ is $g_0$-self-adjoint for all
$x\in M,$) then we define a square of the length of $A$ as
$$
\int_M tr A^2 d\nu
.$$
This metric is invariant under $SO(n)$-valued gauge transformations. Now we define 
a $\log L^2$ $SL(n)$-gauge group as a group of
 measurable sections of $Aut(TM)$ such that with respect to $g_0$,
$$
\int_M {\|\log(A^{\ast}A)\|}^2 d\nu < \infty
.$$
Then ${\mathcal{P}}$ is a homogeneous space under the action of this group. 
We define a metric on ${\mathcal{P}}$ as a unique invariant metric,
which agrees at $g_0$ with the metric just defined.

Now let $(M^{2n},\omega)$ be a compact symplectic manifold. 
Let ${\mathfrak{T}}$ be the twistor bundle, that is, an $Sp(2n)/U(n)$-bundle,
associated with the principal $Sp(2n)$-bundle, defined by $\omega$. 
A smooth section of ${\mathfrak{T}}$ is exactly a tamed 
 almost-complex structure. One then defines a space $Z$ of $\log L^2$-sections of 
${\mathfrak{T}}$ as above  (the $C^{\infty}$-version was used
in [Reznikov 4]).

\noindent\\
{\it Prosition 1.3.}--- The spaces $X$ and $Z$ are Alexandrov and Busemann nonpositively curved.

\noindent\\
{\it Proofs} are standard.

\section{A new invariant of smooth volume-preserving dinamical systems}

Let $(M,\nu)$ be a compact oriented manifold with volume form $\nu$. 
Let $G$ be a finitely generated group which acts on $M$ by smooth
transformations, preserving $\nu$. We are going to define a new dinamical 
invariant which we call $\Lambda$. This is a nonnegative
real number. Though it depends on the choice of a system of generators of $G$,
the crucial fact of whether $\Lambda >0$ or $\Lambda =0$
does not. This relates our $\Lambda$ to Kolmogorov's entropy [Kolmogorov 1]. 
The invariant $\Lambda$ is highly nontrivial already for
$G=\mathbb Z$, that is, as a new  invariant of a volume-preserving diffeomorphism. 
It is also an invariant under conjugation in
 ${{\mathcal{D}}iff}_{\nu}(M)$. A central result of this section is Theorem 2 
below stating that if $G$ is a Kazhdan group then either $\Lambda >0$
 or $G$ fixes a $\log L^2$-Riemannian metric (again this connects $\Lambda$ to the Kolmogorov's entropy).

Let  $g_1,\cdots,g_n$ be a system of generators for $G$. 
Let $X$ be the twistor variety for $(M,\nu)$. Let $\rho$ be the distance function for $X$,
introduced in Section 1. We define $\Lambda$ as the displacement of $G$-action:
$$
\Lambda = \underset{z\in X} \inf \underset{i} \max\, \rho(g_iz,z)
.$$
{\it Proposition 2.1.}---
	 $\Lambda$ is invariant under conjugation in  ${{\mathcal{D}}iff}_{\nu}(M)$.

\noindent\\
{\it Proof.}--- $\rho$ is  ${{\mathcal{D}}iff}_{\nu}$-invariant.

\noindent\\
{\it Proposition 2.2.}---
Let $M=(T^n,can)$ and let $G=\mathbb Z$ act by iterations of a hyperbolic 
element of $SL(n,\mathbb Z)$. Then  $\Lambda >0$.

\noindent\\
{\it Proof.}--- The proof is based on an observation about Alexandrov 
non-positively curved spaces and a trick from [Reznikov 4].

\noindent\\
{\it Lemma.}---
Let $X$ be an Alexandrov non-positively curved space and let 
$\phi : X \rightarrow X$ be an isometry which leaves invariant a geodesic
$\gamma$ of $X$. Then the displacement of $\phi$ is realized on the points of $\gamma$, 
that is, for $y \in \gamma$,
$$
\rho(y,\phi y)=\underset{x\in X} \min\, \rho(x,\phi x).
$$
{\it Proof.}---
For $x\in X$ let $y\in \gamma$ be a point which realizes the distance from $x$ to $\gamma$. 
Then $\rho(y,\phi y)\leq \rho(x,\phi x)$.\\

Now let $X$ be the twistor space of $T^n$ and let $Y \subset X$ be the space of metrics, 
invariant under shifts (we view $T^n$ as
a Lie group). Then $Y$ is totally geodesic in $X$, because it is a manifold of 
fixed points of a family of isometries. As a Riemannian
manifold, $Y\simeq SL(n)/SO(n)$. Any hyperbolic matrix $\phi$ by 
definition leaves invariant a geodesic in $Y$. The result follows.

A main result in the theory of invariant $\Lambda$ is as follows.

\noindent\\
{\it Theorem 2.3.}--- Let $G$ be a Kazhdan group acting on a compact 
oriented manifold $(M,\nu)$ preserving a
volume form $\nu$. Then either $\Lambda >0$ or $G$ fixes a $\log L^2$-metric on $M$.\\

{\it Proof.}--- Consider an isometric action on $X$. If the displacement 
function $\underset{i}\sup\, \rho(g_{i}z,z)$
is not bounded away from zero, then either there is a fixed point $z_0 \in X$ for $G$,
or $G$ is not Kazhdan, by a result of
[Kovevaar-Schoen 1]. The result follows.

\section{Non-linear superrigidity alternative}

{\it Theorem 3.1.}---Let $G$ be either a  semisimple Lie group of rank$\geq 2$, or $Sp(n,1)$ or $Iso(\mathbb{C}a\mathbb{H}^2)$. 
Let $\Gamma \subset G$ be a
cocompact lattice. Let $(M^n,\nu)$ be a compact oriented manifold,
on which $\Gamma$ acts preserving the volume form $\nu$.
Then either 

\vspace{2mm}\noindent
a) $\Gamma$ preserves a $\log L^2$- metric on $M$, or

\vspace{2mm}\noindent
b) there exists a sequence $g_0,g_1,\cdots$ of smooth Riemannian metrics 
on $M$ with volume form $\nu$ such that
$$
\int_M {\| \log A_i\|}_{g_0}^2 d\nu \rightarrow \infty,
$$
where $g_i = g_0(A_i\cdot,\cdot)$ and
$$
0<{const}_1<\underset{j} \sup \int_M {\| \log B_{ij}\|}_{g_i}^2 d\nu 
	< {const}_2,  \qquad (i \rightarrow \infty),
$$
where ${\gamma}^{\ast}_{j}g_i=g_i(B_{ij}\cdot,\cdot)$, $\{{\gamma}_j\}$ is a 
fixed finite set of generators for $\Gamma$, or

\vspace{2mm}\noindent
c) there is a nonconstant totally geodesic $\Gamma$-invariant map
$$
\Psi : G/K \rightarrow X,
$$
where $K$ is a maximal compact subgroup of $G$.

\noindent\\
{\it Remarks}.-

\vspace{2mm}\noindent
1) \quad In case b) we say that $\Gamma$ nearly fixes a $\log L^2$-metric on $M$.

\vspace{2mm}\noindent
2) \quad the case c) implies, for $G$ simple, that $\dim G/K \leq \dim SL(n)/SO(n)$, a so-called Zimmer conjecture.

\vspace{2mm}\noindent
3) \quad for $G=SL(m, \mathbb R), m\geq 3$ and $n=m$, one deduces in case c) an existence of 
a measurable frame field
$\hat {e}(x), \hat{e}=(e_1,\cdots , e_n)$, such that for almost all $x\in M$,
$$
\pi (\gamma)_*[\hat{e}(x)]=\gamma \hat{e}(\pi(\gamma)x)
$$
where $\gamma\in\Gamma$ and $\pi(\gamma)$ is an action of $\gamma$ on $M$.

\vspace{2mm}\noindent
4) \quad Conversely, a standard action of $SL(n,\mathbb Z)$ on $T^n$ does 
not satisfy a) (which is well-known) and b). To see this,
we notice that $SL(n,\mathbb Z)$ leaves invariant a totally geodesic 
space $Y$ introduced in the proof of Proposition 2.2. The
argument of this proof implies that it is enough to show that the displacement 
function of the action of $\Gamma$ on $Y$
diverges to $\infty$ as one escapes all compact sets of $Y$. This follows from 
the fact that $Y$ is a Riemannian symmetric
space of non-compact type and $\Gamma$ does not fix a point at infinity of $Y$.

\vspace{2mm}\noindent
5) \quad The statement of Theorem constitutes a definite progress in 
the nonlinear superrigidity problem. There is still
a mystery in the option b) where one would prefer a statement that $\Gamma$ 
fixes a ``point at infinity" of the space of metrics $X$,
perhaps a measurable distribution of $k$-dimension planes, $k\leq n$. 
At the time of writing this chapter (August, 1999) I am
unable to make such a reduction.

\noindent\\
{\it Proof.} follows a long-established tradition [Siu 1], [Corlette 1],
[Mok-Siu-Yeung 1], see also a treatment of [Jost-Yau 1],
in a new infinite-dimensional target context. If neither a) or b) holds then,
accounting that $\Gamma$ is Kazhdan, we deduce that the displacement function of $\Gamma$ 
tends to infinity as one escapes all bounded sets in $X$. Let
$F\rightarrow \Gamma \setminus G/K $ be a flat fibration with fiber $X$,
corresponding to the action of $\Gamma$ in $X$. A theorem
of [Kovevaar-Schoen 1], or [Jost 1] implies that there is a harmonic section of $F$. 
By Propositon 1.2 and main theorem of [Corlette 1] and [Mok-Siu-Yeung 1],
this section must be totally geodesic. The result follows.\\

In the case of symplectic action of lattice $\Gamma$ on a compact symplectic manifold $(M,\omega)$ we have a comletely similar
theorem, as follows.

\noindent\\
{\it Theorem 3.2.}--- Let $G$ be either a semi-simple Lie group of rank$\geq 2$, 
or $Sp(n,1)$ or $Iso(\mathbb{C}a\mathbb{H}^2)$, $ \Gamma$ a cocompact lattice in $G$ which
acts symplectically on a compact symplectic manifold $(M^{2n},\omega)$. Then either 

\vspace{2mm}\noindent
a) \quad $\Gamma$ fixes a $\log L^2$ tamed almost-complex structure $J$, 
 or

\vspace{2mm}\noindent
b) \quad there exists a sequence of tamed smooth almost-complex structures 
$J_i \in Z$ with $\rho(J_0,J_i) \rightarrow \infty$
and
$$
0<{const}_2<\underset{j} \sup\, \rho({\gamma}_j J_i,J_i)<{const}_1, \qquad\mbox{or}
$$
c) \quad there is a $\Gamma$-invariant totally geodesic map
$$
\Psi : G/K \rightarrow Z
.$$
{\it Proof.} is exactly as above.\\

In case c) and $G$ simple it follows that $\dim G/K \leq \dim Sp(2n)/U(n)$. 
If $M=(T^{2n},can),\ \ G=Sp(2n,\mathbb R)$ and case c)
one deduces an existence of a measurable symplectic frame $\hat{e}(x)=(e_1,\cdots, e_{2n}(x))$,
such that for $\gamma \in \Gamma,$
$$
{\pi(\gamma)}_{\ast}[\hat{e}(x)]=\gamma \hat{e}(\pi(\gamma)(x))
.$$

\chapter{K\"{a}hler and quaternionic K\"{a}hler groups}

In a letter to the author [Deligne 1] P.Deligne asked if one can 
extend the author's theorem on rationality of secondary characteristic
classes of a flat bundle over a projective variety to quasiprojective varieties. 
In 1994 the author was able to answer this  question
positively for the special case of noncompact ball quotients using an analytic 
technique of [Gromov-Schoen 1] and the scheme of the
original proof for projective varietes. Here we present a full answer to Deligne's 
question, Theorem 1.1,  using an analytic technique of [Jost-Zuo 1], who
produced harmonic maps of infinite energy but controlled growth.

We then turn to a well-known open problem of finding restriction on topology 
of compact quaternionic K\"{a}hler manifolds. In case of
positive scalar curvature the situation is well-understood,
but in case of negative scalar curvature the twistor spaces of [Solomon 1]
are not K\"{a}hler and its technique fails. The only result known was a 
theorem of [Corlette 1] stating that the  fundamental group
does not have infinite linear representations unless the manifold is locally symmetric. 
Our result, Theorem 2.2, states that the
fundamental group is Kazhdan. This is of course,
a severe restriction (Kazhdan groups are rare). As a by-product of our technique,
we obtain a  new proof of a  classical theorem,
stating that the lattices in semisimple Lie groups of rank $\ge 2$, $Sp(n,1)$ and $Iso(\mathbb{C}a\mathbb{H}^2)$ are Kazhdan. We also show using I.1 that
the classes in second cohomology space of a K\"ahler 
non-Kazhdan group, constructed in [Reznikov 6] and shown there nontrivial, are of
polynomial growth. This again is very rare for ``just a group",
as polynomial growth in cohomology is connected to a polynomial
isoperimetric inequality in the Cayley graph, which needs a special reason to hold. 
This means K\"{a}hler groups are rare, too.

\section{Rationality of secondary classes of flat bundler over quasiprojective varietes}

A rationality theorem for secondary classes of flat bundles over 
compact K\"{a}hler manifolds (previously known as Bloch conjecture
[Bloch 1]) has been  proved in 1993 in [Reznikov 3] and [Reznikov 5]. 
In a letter to the author [Deligne 1] P.Deligne asked if one can prove
such a statement for local system with logarithmic singularities 
over a quasiprojective variety. The answer happens to be yes.

\noindent\\
{\it Theorem 1.1.}---
Let $X$ be a quasiprojective variety, $\rho : {\pi}_1(X) \rightarrow SL(n,\mathbb C)$ 
a representation. Let $b_i(\rho)$ be the imaginary part and
$ChS_i(\rho)$ the $\mathbb R / \mathbb Z$-part of the secondary 
class $c_i(\rho)\in H^{2i-1}(X,\mathbb C / \mathbb Z)$ of the flat bundle with
monodromy $\rho$. Then\\
A. \quad $b_i(\rho)=0 \quad (i\geq 2)$ (the Vanishing Theorem).\\
B. \quad $ChS_i(\rho)\in H^{2i-1}(X,\mathbb Q / \mathbb Z)$ (the Rationality Theorem).

\noindent\\
{\it Proof.}---
For any smooth manifold, A implies B, as explained in the above cited papers. 
So we only prove A. Again it is explained in the above cited papers that
we may assume $\rho$ to be irreducible. Then by a recent result [Jost-Zuo 1] 
an associated $SL(n,\mathbb C)/SU(n)$ flat bundle over $X$
possesses a pluriharmonic section $s$ which satisfies the Sampson degeneration condition. 
This means the following. The derivative
$Ds_x$, $x\in X$ can be viewed as a 
$\mathbb R$-linear map to the space $P$ of Hermitian matrices. Let
${(Ds_x)}^{\pm}_{\mathbb C}(Y)=(Ds_x(Y)\pm \sqrt{-1}Ds_x(\sqrt{-1}Y))$ be 
a map of $TX $ to $P\otimes \mathbb C$. Then the image of
${(Ds_x)}^{\pm}_{\mathbb C}$ consists of commuting matrices. 
Now a first proof of the Main Theorem in [Reznikov 5] applies
word-to-word and the result follows.
\\



\section{Kazhdan property $T$ for K\"{a}hler and quaternionic K\"{a}hler groups}

There are two ways to geometrize group theory. One approach (a time geometry in the terminology of [Reznikov 7]) is 
to consider finitely generated groups which act on a (usually compact) space with
some structure (a volume form, a symplectic form,
a tree, a circle, a conformal structure, etc). An amazing phenomenon,
amply demonstrated
in the previous chapters is that these groups tend to be not Kazhdan. 
Another approach (a space geometry) is to consider groups which are fundamental groups
 of a compact (or closed to compact) manifold with some structure 
(like K\"{a}hler). It happens that these groups tend to be Kazhdan.
Therefore these two families of ``geometric" groups are essentically disjoint. 
A following result is a main theorem of [Reznikov 6].

\noindent\\
{\it Theorem.}--- Let $G$ be a fundamental group of a compact K\"{a}hler manifold. 
If $G$ is not Kazhdan, then
$H^2(G,\mathbb R) \neq 0$. Moreover, if $H$ is not Kazhdan and $\psi: G \rightarrow H $ 
is surjective then
$0\neq{\psi}^{\ast}:H^2(H,\mathbb R)\rightarrow H^2(G,\mathbb R).$\\

I would like to notice an important property, which I overlooked in [Reznikov 6].

\noindent\\
{\it Proposition 2.1.}---Under the conditions of the Theorem, there is a nontrivial 
class of polynomial growth in
$H^2 (G,\mathbb R).$

\noindent\\
{\it Proof.}---
There is a unitary representation $\rho:G\rightarrow U({\mathcal{H}})$ and a 
class $l\in H^1(G,{\mathcal{H}})$ such that a class $\gamma$
in $H^2(G,\mathbb R)$ given by $\langle l,l \rangle$ is nonzero,
where $\langle \cdot,\cdot \rangle$ is an imaginary part of the
scalar product in ${\mathcal{H}}$. This is proved in [Reznikov 6]. Now the result follows from Lemma I 1.1.\\

It is extremely rare for a finitely generated group to have nonzero polynomial cohomology.

We now turn to quaternionic K\"{a}hler manifolds.
If a scalar curvature is positive, then the topology is very well understood [Solomon 1].
On the contrary, if the scalar curvature is negative,
the only result known is that the fundamental group satisfy the geometric superrigidity
[Corlette 1]. This means if $ {\pi}_1(X)$ admits a Zariski-dense representation 
in an algebraic Lie group, then  $ {\pi}_1(X)$ is a lattice, and
$X$ a symmetric space of a known type. However, it is a rare occasion 
for a group to have any finite dimensional linear representation with infinite image.
Using a combination of ideas of [Corlette 1] and [Reznikov 6] which is 
based on [Korevaar-Schoen 1] we now prove a very stong structure
 theorem.

\noindent\\
{\it Theorem 2.2.}---Let $X$ be a quaternionic K\"{a}hler manifold of negative scalar curvature. 
Then  $ {\pi}_1(X)$ is Kazhdan.

\noindent\\
{\it Proof.}---Suppose not. Then by [Korevaar-Schoen 1] 
there exists an affine flat Hilbert bundle $E$ over $X$ with a
nonparallel harmonic section. By a vanishing result of [Corlette 1] 
this section must be totally geodesic. Then $X$ must be covered
by a flat torus, a contradiction.

\noindent\\
{\it Remark.}--- The same argument provides a new proof of 
the classical theorem [Kazhdan 1], [Kostant 1] , that the (cocompact) lattices in
semisimple Lie groups of rank$\geq 2$, $Sp(n,1)$ and $Iso(\mathbb{C}a\mathbb{H}^2)$ are Kazhdan. One uses a vanishing 
result of [Mok-Siu-Yeung 1] ( see also a treatment in [Jost-Yau 1]) for lattices in semisimple Lie groups of rank$\geq 2$,  and the above-mentioned result of [Corlette 1] for $Sp(n,1)$ and  $Iso(\mathbb{C}a\mathbb{H}^2)$. Once established for cocompact lattices, the result follows for all lattices because a Lie group and a lattice in it are or are not Kazhdan at the same time.



\chapter*{REFERENCES}
\small

[Adams-Ratiu-Schmid 1] M.Adams, T.Ratiu, R.Schmid.---{\it The Lie grouip structure of diffeomorphism groups and Fourier integral operators with applications},
		in: Infinite dimensional Lie groups with applications, V.Kac, Editor, Springer, 1985.\\

\noindent
[Alperin  1] R.Alperin.---{\it Locally compact groups acting on trees and property T },
		Mh.Math, {\bf 93} (1982), 261-265.\\

\noindent
[Beurling-Ahlfors 1] A.Beurling,  L.Ahlfors.---{\it The boundary correspondence under quasiconformal mappings},
		Acta Math., {\bf 96} (1956), 125--142.\\

\noindent
[Arakelov 1] \\

\noindent
[Bayen-Flato-Fronsdal-Lichnerowicz-Sternheimer ] F.Bayen, M.Flato, C.Fronsdal, A.Lichnerowicz, 
		D.Sternheimer.--- {\it Deformation theory and quantization}, Ann. Phys.  {\bf III} (1978), 61--151.\\

\noindent
[Belavin 1] A.A.Belavin.--- {\it Discrete groups and integrability of quantum 
		systems}, Funct. Anal. Appl, {\bf 14} (1980), 18-26 (Russian).\\

\noindent
[Belavin-Polyakov-Zamolodchikov 1] A.A.Belavin, A.M.Polyakov, A.B.Zamolodchikov.---
		{\it Infinite conformal symmetry in two-dimensional quantum field 
			theory}, Nucl. Phys., {\bf B241} (1984), 333-380.\\

\noindent
[Brown 1] K.S.Brown.--- 
		{\it Cohomology of Groups},
			Springer.\\

\noindent
[Bloch 1] S.Bloch.--- 
		{\it Applications of the dilogarithm functions in algebraic K-theory and algebraic geomery},
			in :Proc. Int. Symp. Alg. Geom, Kyoto, Kinokumiya, 1977,  103--114.\\

\noindent
[Bojarski 1] B.Bojarski.--- 
		{\it Generalized solutions of a system of differential equations of the first order with discontinuous coefficients},
			Math. Sbornik, {\bf 43}  (1957), 451--503 (Russian).\\

\noindent
[Bott 1] R.Bott.--- 
		{\it On the characteristic classes of groups of diffeomorphisms},
			Enseign. Math.  {\bf 23} (1977), 209-220.\\

\noindent
[Besson-Courtois-Gallot 1] G.Besson, G.Courtois, S.Gallot.---
		{\it Entropies et rigidit\'es des espaces localement sym\'etriques de courbeure 
			strictement n\'egative}, GAFA  {\bf 5} (1995), 731-799.\\

\noindent
[Brown-Georghegan 1] K.S.Brown, R.Georghegan.---
		{\it An infinite dimensional torsion free $FP_{\infty}$ group}, Invent. Math.  {\bf 77} (1984), 367--381.\\

\noindent
[Connes-Gromov-Moscovici 1] A.Connes, M.Gromov, H.Moscovici.--- 
		{\it Group cohomology with Lipschitz control and higher signatures}, 
			GAFA  {\bf 3} (1993), 1-78.\\

\noindent
[Cannon-Thurston 1] Cannon, W.Thurston.--- 
		{\it Equivariant Peano curves}, 
			Preprint , 1986.\\

\noindent
[Carleson 1] L.Carleson.--- 
		{\it The extension problem for quasiconformal mappings}, in:
			Contributions to Analysis, AP, 1974, 39--47.\\

\noindent
[Carleson 2]\\

\noindent
[Connes-Moscovici 1]
A.Connes, H.Moscovici.---
		{\it Cyclic cohomology, the Novikov conjecture and hyperbolic groups},
			Topology {\bf 29} (1990), 345-388.\\

\noindent
[Corlette 1] K.Corlette.---
		{\it Archimedian superrigidity and hyperbolic geometry}, 
Ann. Math. {\bf 135} (1990), 165--182.\\

\noindent
[Deligne 1] P.Deligne.---{\it A letter to the author}, 1994.\\

\noindent
[Duady-Earle 1] A.Duady, C.J.Earle.--- 
		{\it Conformally natural extensions of homeomorphisms of the circle},
			Acta Math.,  {\bf 157} (1986), 23--48\\

\noindent
[Dupont 1] J.L.Dupont.--- 
		{\it Simplicial de Rham cohomology and characteristic classes of flat bundles},
			Topology, {\bf 18} (1979), 295--304.\\

\noindent
[Edmunds-Opic 1]
D.E.Edmunds, B.Opic.--- 
		{\it Weighted Poincar\'e and Friedrichs inequalities},
			J. London Math. Soc.,  {\bf 47} (1993), 79-96.\\

\noindent
[Edmunds-Triebel 1]
D.E.Edmunds, H.Triebel.---
		{\it Logarithmic Sobolev spaces and their applications to sectral theory},
			Proc. London Math. Soc. {\bf 71} (1995), 333-371.\\

\noindent
[Farb-Shalen 1] B.Farb, P.Shalen.---
		{\it Real-analytic action of lattices},
			Invent. Math.  {\bf 135} (1999), 273--296.\\

\noindent

\noindent
[Fedosov 1] Fedosov.---
	{\it Index theorems},
		 {\bf  } (), (Russian).\\

\noindent
[Feigin-Tsygan 1]
B.L.Feigin, B.L.Tsygan.---
		{\it Cohomology of Lie algebras of generalized-Jacobean matrices},
			Funct. Anal. Appl.  {\bf 17} (1983), 86-87(Russian).\\

\noindent
[Fuks 1]
D.Fuks.---
	{\it Cohomology of Infinite-Dimensional Lie Algebras},
		Consultants Bureau, New York and London, 1986.\\

\noindent
[Furstenberg 1] H.Furstenberg.--- 
		{\it A Poisson formula for semi-simple Lie groups},
			Ann. Math., {\bf 77} (1963), 335--386.\\

\noindent
[Furstenberg 2] H.Furstenberg.---{\it Boundary theory and stochastic procecesses on homogeneous spaces},
		Proc. Symp. Pure Math.  {\bf 26} (1973), 193--233.\\

\noindent
[Garnett 1] L.Garnett.---{\it Foliations, the ergodic theorem and Brownian motion},
			J. Funct. Anal., {\bf 51}(1983), 285--311.\\

\noindent
[Gelfand-Fuks 1]
I.M.Gelfand, D.B.Fuks.---
	{\it Cohomology of Lie algebras of vector fields on the circle},
		Funct. Anal. Appl. {\bf 2} (1968), 92-93(Russian).\\

\noindent
[Gardiner-Sullivan 1] F.P.Gardiner, D.Sullivan.---
	{\it Symmetric structures on a closed curve},
	Amer. J. Math. {\bf 114} (1999), 683--736.\\

\noindent

[Ghys 1]
E.Ghys.---
	{\it Actions de r\'eseaux sur le cercle},
	Invent. Math. {\bf 137} (1999), 199-231.\\

\noindent
[Ghys-Sergiescu 1]
E.Ghys, V.Sergiescu.---
	{\it Sur un groupe remarquable de diff\'eomorphismes du cercle},
	Comment. Math. Helv.  {\bf 62} (1987), 185--239.\\

\noindent
[Goldshtein-Kuzminov-Shvedov 1]V.M.Goldshtein, V.I.Kuzminov, I.A.Shvedov ,{\it On a problem of Dodziuk},
		Trudy mat. Inst. Steklov,  {\bf 193} (1992), 72-75.\\

\noindent
[Gromov-Schoen 1] M.Gromov, R.Schoen---
	{\it Harmonic maps into singular spaces and $p$-adic rigidity for lattices in groups of rank one},
		Publ. Math. IHES,  {\bf 76} (1992), 165--246.\\

\noindent
[Gehring 1]F.Gehring.--- 
		{\it The $L^p$-integrability of the partial derivatives of a quasiconformal mapping},
			Acta Math., {\bf 130}(1973), 265--277.\\

\noindent
[Guichardet 1] A.Guichardet.--- 
		{\it Cohomologie des groupes topologiques et des alg\'ebres de Lie},
			Cedic/Fernand Nathan, 1980.\\

\noindent
[Greenberg-Sergiesku  1] P.Greenberg, V.Sergiesku.--- {\it An algebraic extension of the braid group},
			Comment. Math. Helv. {\bf 66} (1991), 109--138.\\

\noindent
[Guba 1]
V.S.Guba.---
	{\it Polynomial upper bounds for the Dehn function of R.Thompson group $F$},
		Journ. Group Theory  {\bf 1} (1998), 203--211.\\

\noindent
[de la Harpe-Valette 1]
P.de la Harpe, A.Valette.---
	{\it La propri\'et\'e (T) de Kazhdan pour les groupes localement compactes},
		Ast\'erisque {\bf 175} (1989).\\

\noindent
[Hewitt-Ross 1] E.Hewitt, K.A.Ross.---
	{\it Abstract Harmonic Analysis}, Springer, 1970.\\

\noindent
[Hopf 1]
E.Hopf.---
	{\it Statistik der geodatischen linien in manigfaltigkeiten negativer kr\"ummung},
		Ber. Verb. Sachs. Akad. Wiss. Leipzig  {\bf 91} (1939), 261-309.\\

\noindent
[Hurder-Katok]
S.Hurder, A.Katok.---\\

\noindent
[H\"amenstadt 1] U.H\"amenstadt.---{\it A lecture in Leipzig conference "Perspectives in Geometry"}, 1998.\\

\noindent
[Jost 1]J.Jost.---{\it Equillibrium maps between metric spaces}, Calc. Var. {\bf 2} (1994), 173--204.\\

\noindent
[Jost-Yau 1]J.Jost, S.T.Yau.---{\it Harmonic maps and superrigidity}, Proc. Symp. Pure Math., {\bf 54} (1993), 245--280.\\

\noindent
[Jost-Zuo 1] J.Jost, K.Zuo.---{\it Harmonic maps of infinite energy and rigidity results for representations of fundamental groups of quasiprojective varietes},
		J. Diff. Geom.,  {\bf 47} (1997), 469--503.

\noindent
[Kaimanovich-Vershik 1]V.A.Kaimanovich, A.M.Vershik.--- Ann. Probab. {\bf 11} (1983), 457-490.\\

\noindent
[Koosis 1]
P.Koosis.---
	{\it Introduction to $H_p$ Spaces},
		Cambridge UP, 1998.\\

\noindent
[Kudryavcev 1]L.D.Kudryavcev.---{\it Habilit\"ationschrift},
		Steklov Math. Inst., 1956.\\

\noindent
[Kudryavcev 2]L.D.Kudryavcev.---{\it Direct and inverse imbedding theorems. Applications to the solutions of elliptic equations by variational methods} Trudy Steklov Math. Inst. {\bf 55} (1959)\\

\noindent
[Korevaar-Schoen 1]N.Korevaar, R.Schoen.---{\it Global existence theorems for harmonic maps to non-locally compact spaces}, 
		Comm. Geom. Anal., {\bf 5} (1997), 333-387.\\

\noindent
[Korevaar-Schoen 2]N.Korevaar, R.Schoen.---{\it Sobolev spaces and harmonic maps for metric space targets}, 
		Comm. Geom. Anal., {\bf 1} (1997), 561--659.\\
\noindent
[Kostant 1] B.Kostant.---{\it On the existence and irreducubility of certain series of representations}, 
		Lie Groups and Their Representations, Halsted, NY, 1975, 231--329.\\
\noindent
[Lions 1]J.L.Lions.---{\it Th\'eor\`ems de trace et d'interpolation},
		 I, Ann. Schuola Norm. Super. Pisa,  {\bf 13} (1959), 389-403.\\

\noindent
[Lizorkin 1]
P.I.Lizorkin.---{\it Boundary values of functions from "weight" classes},
		Sov. Math. Dokl. {\bf 1} (1960), 589-593.\\

\noindent
[Lizorkin 2]
P.I.Lizorkin.---
	{\it Boundary values of a certain class of functions}, 
		Dokl. Anal. Nauk SSSR {\bf 126} (1959), 703-706(Russian).\\

\noindent
[Makarov 1]
N.G.Makarov.--- 
	{\it On the distortion of boundary sets under conformal mappings},
		Proc. London Math. Soc. {\bf 51} (1985), 369-384.\\

\noindent
[Makarov 2]
N.G.Makarov.---
	{\it On the radial behaviour of Bloch functions},
		Soviet Math. Dokl. {\bf 40} (1990), 505-508.\\

\noindent
[Matsumoto-Morita 1] S.Matsumoto, S.Morita.---
	{\it Bounded cohomology of certain groups of homeomorphisms},
		PAMS {\bf 94} (1985), 539-544.\\

\noindent
[Morita 1] S.Morita.---{\it Characteristic classes of surface bundles},
		Invent. Math.,  {\bf 90} (1987), 551--577\\

\noindent
[Morita 2] S.Morita.---{\it Characteristic classes of surface bundles and bounded cohomology},
		in: A F\^ete of Topology, AP, 1988, 233--257.\\

\noindent
[Mikhailov 1]\\

\noindent
[Miller 1]
E.Y.Miller.---{\it The homology of the mapping class group},
		J. Diff. Geom. {\bf 24} (1986), 1--14.
\\

\noindent
[Mishchenko 1] A.S.Mishchenko.--- 
	{\it Infinite-dimensional representation of discrete groups and higher 
    	siguatures},
		Math. USSR. Izv.{ \bf 8} (1974), 85-111.\\

\noindent
[Mishchenko 2] A.S.Mishchenko.---
	{\it Hermitian K-theory, the theory of characteristic classes and
    	methods of functional analysis},
		Russian Math. Surveys  {\bf 31} (1976), 71-138.\\

\noindent
[Mok-Siu-Yeung 1]Mok, Siu, Yeung.---
	{\it },
		Invent. Math. {\bf  } (1993),  .\\

\noindent
[Mumford 1]
D.Mumford.---
	{\it Towards an enumerative geometry of the moduli space of curves},
		in: Arithmetic and geometry, Progress in Math.  {\bf 36}, Birkh\"auser, 1983, 271--328.\\

\noindent
[Murray-von Nuemann 1]F.J.Murray, J. von Neumann.---{\it On rings of operators}, Ann. Math {\bf 37} (1936), 116--129, TAMS {\bf 41} (1937), 208--248, Ann. Math., {\bf 41} (1940), 94--161,  {\bf 44} (1943), 716--808.\\

\noindent
[Naboko 1]S.N.Naboko.---
	{\it Nontantengial boundary values of operator-valued R-functions},
		Leningrad Math. Journ.   {\bf 1} (1990), 1255-1278.\\

\noindent
[Nag 1]S.Nag.---\\

\noindent
[Nicholls 1]
P.J.Nicholls.---
	{\it A measure on the limit set of a discrete groups in Ergodic Theory},
		in: Symbolic Dynamics and Hyperbolic Spaces, T.Bedford, M.Keane,
			C.Series, eds, Oxford UP, 1991, 259-296.\\

\noindent
[Otal 1] J-P.Otal.---
	{\it Le th\'eor\`eme d'hyperbolization pour les vari\'et\'es fibr\'ees sur le circle }, Pr\'epublication Orsay.\\

\noindent
[Palais 1] R.Palais.---
	{\it On the homotopy type of certain groups of operators},Topology {\bf  3} (1965), 271--279\\

\noindent
[Pietsch 1]
A.Pietsch.---
	{\it Operator Ideals}, VEB,Berlin, 1978.\\

\noindent
[Pommerenke 1]
Ch.Pommerenke.---
	{\it Boundary Behaviour of Conformal Maps},
		Springer, 1992.\\

\noindent
[Pansu 1] P.Pansu.---
	{\it Cohomologie $L^p$ des vari\'et\'es \`a courbure n\'egative, cas du degr\'e un,}, Rend. Sem. Mat. Torino  (1989), 95--120.\\

\noindent
[Pansu 2] P.Pansu.---
	{\it Differential forms and connections adapted to a contact structure, after M.Rumin}, in: Symplectic Geometry, D.Salamon, Editor, Cambridge UP (1993), 183--195.\\

\noindent
[Pressley-Segal 1]
A.Pressley, G.Segal.---
	{\it Loop Groups},
		Clarendon Press, Oxford, 1986.\\

\noindent
[Reimann 1]H.M.Reimann.---
	{\it Functions of bounded mean oscillation and quasiconformal mappings},
		Comment. Math. Helv. {\bf 49} (1974), 260--276.\\

\noindent
[Rempel-Schulze 1]
S.Rempel, B.-W.Schulze.--- 
		{\it Index Theory of Elliptic Boundary Problems},
			Academie-Verlag, Berlin, 
1992.\\

\noindent
[Reznikov 1] A.Reznikov.---
	{\it The space of spheres and conformal geometry},\\
		Riv. Math. Un. Parma {\ bf 17} (1991), 111-130.\\

\noindent
[Reznikov 2] A.Reznikov.---
	{\it Characteristic classes in symplectic topology}, 
		Sel. Math. {\bf 3} (1997), 601--642.\\

\noindent
[Reznikov 3]  A.Reznikov.---
	{\it Rationality of secondary classes},
		J. Diff.Geom. {\bf 43} (1996), 674--682.\\

\noindent
[Reznikov 4] A.Reznikov.---
	{\it Continuous cohomology of volume-preservivg and symplectic diffeomorphisms, 
	measurable transfer and higher asymptotic cycles}, 
		Sel. Math. {\bf 5} (1999), 181--198.\\

\noindent
[Reznikov 5] A.Reznikov.---
	{\it All regulators of flat bundles are torsion},
		Ann. Math.  {\bf 141} (1995), 373--386.\\

\noindent
[Reznikov 6] A.Reznikov.---
	{\it Structure of K\"ahler groups,I: second cohomology},
		Preprint, April, 1998 (Math. DG 9903023).\\

\noindent
[Reznikov 7] A.Reznikov.---
	{\it Analytic Topology},
		in:Proceedings of the European Congress of Mathematics, to appear.\\

\noindent
[Reznikov 8] A.Reznikov.---
	{\it Arithmetic Topology of units, ideal classes and three and a half-manifolds},
		in preparation.\\

\noindent
[Reznikov 9] A.Reznikov.---
	{\it Harmonic maps, hyperbolic cohomology and higher Milnor inequalities},
		Topology {\bf 32} (1993),  899--907.\\

\noindent
[Reznikov 10] A.Reznikov.---
	{\it Analytic Topology II},
		in preparation.\\

\noindent
[Siu 1] Y.T.Siu.---
	{\it The complex-analycity of harmonic maps and strong rigidity of compact K\"ahler manifolds},
		Ann. Math. {\bf 112} (1980), 73--111.\\

\noindent
[Segal 1] G.B.Segal.---
	{\it Unitary representations of some infinite dimensional groups},
		Comm. Math. Phys.  {\bf 80} (1981), 301-342.\\

\noindent
[Solomon 1]
S.Solomon.---
	{\it Quaternionic K\"ahler manifolds},
		Invent. Math. (1981), 301-342\\

\noindent
[Sullivan 1] D.Sullivan.---
	{\it Entropy, Hausdorff measures old and new, and limit sets of geometrically finite Kleinian groups}, 
		Acta Math.  (1989), 259-277\\

\noindent
[Sullivan 2]\\

\noindent
[Sullivan 3]D.Sullivan.--- 
		{\it On the ergodic theory at infinity of arbitrary discrete group of hyperbolic motions},
			in: Riemann Surfaces and Related Topics, I.Kra and B.Maskit, Editors, Up (1981)., 465--496.\\

\noindent
[Thompson 1]J.Thompson.--- , Unpublished.

\noindent
[Triebel 1]
H.Triebel.--- 
		{\it Theory of Function Spaces},
			Birkh\"auser, 1983.\\

\noindent
[Tukia 1]
P.Tukia.---
	{\it The Hausdorff dimension of the limit set of a geometrically finite Kleinian group},
		Acta Math. {\bf 152} (1989), 127-140.\\

\noindent
[Thurston 1]W.Thurston.---{\it The Geometry and Topology of Three-Manifolds}, Princeton Lecture Notes.\\

\noindent
[Thurston 2]W.Thurston.---{\it Hyperbolic structure on 3-manifolds, II: surface groups and 3-manifolds which fiber ove the circle}, Preprint, 6 August 1986.\\

\noindent
[Tukia-V\"ais\"al\"a] P.Tukia, J.V\"ais\"al\"a .--- 
		{\it Quasiconformal extension from dimension $n$ to $n+1$},
			Ann. Math., {\bf 115} (1982), 331--348.\\

\noindent
[Vasharin 1]
A.A.Vasharin.---
	{\it The boundary properties of functions having a finite Dirichlet integral with a weight},
		Dokl. Anal. Nauk SSSR  {\bf 117} (1957), 742-744.\\

\noindent
[Verdier 1]J.-L.Verdier.---
	{\it Les r\'epresentations des alg\'ebres de Lie affines: applications \`a quelques probl\`emes de physique (d'apr\`es E.Date, M.Jimbo, M.Kashivara, T.Miwa},
		S\'eminaire Bourbaki, Expos\'e 596 (1981--1982), 1--13.\\

\noindent
[Vey 1]\\

\noindent
[Uspenski\u\i 1]
S.V.Uspenski\u\i.---
	{\it Imbedding theorems for weighted classes},
		Trudy Math. Inst. Steklov {\bf 60} (1961), 282-303 (Russian)(Amer. Math. Soc. Trans.  {\bf 87} (1970)).\\

\noindent
[Watatani 1]
Y.Watatani.---
	{\it Property (T) of Kazhdan implies property (FA) of Serre},
		Math. Japon. {\bf 27} (1981), 97-103.\\

\noindent
[Zimmer 1]R.Zimmer.---
	{\it Kazhdan groups acting on manifolds},
		Invent. Math., {\bf 75} (1984), 425--436.\\

\noindent
[Zimmer 2] R.Zimmer.---{\it Lattices in semisimple groups and invariant geometric structures on compact manifolds}, Discrete Groups in Geometry and Analysis, Progress in Math.  {\bf 67} (1987),
			152--210.\\
\noindent
[Zimmer 3] R.Zimmer.---

		{\it Strong rigidity for ergodic actions of semisimple Lie groups},
			Annals of Math., {\bf 112} (1980),
			511--529.\\

\end{document}